\documentclass{article}

\usepackage{amsmath,amssymb,amsthm,mathrsfs}
\usepackage{natbib}
\usepackage{booktabs}
\usepackage{algorithm,algorithmic}
\usepackage{enumitem}
\usepackage{fullpage}
\usepackage{xurl}
\usepackage{ragged2e}
\usepackage[colorlinks=true,citecolor=blue,linkcolor=blue,urlcolor=blue]{hyperref}
\usepackage{mathtools}
\mathtoolsset{showonlyrefs}

\newtheorem{theorem}{Theorem}
\newtheorem{lemma}[theorem]{Lemma}
\newtheorem{proposition}[theorem]{Proposition}
\newtheorem{corollary}[theorem]{Corollary}
\theoremstyle{definition}
\newtheorem{definition}[theorem]{Definition}

\theoremstyle{remark}
\newtheorem{remark}[theorem]{Remark}

\newcommand{\cP}{\mathcal{P}}
\newcommand{\cQ}{\mathcal{Q}}
\newcommand{\LL}{\operatorname{LL}}
\newcommand{\eff}{\operatorname{eff}}
\newcommand{\lstar}{\ell^\star}

\renewcommand{\Pr}{\field{P}}

\newcommand{\bh}{\boldsymbol{h}}

\newcommand{\bq}{\boldsymbol{q}}

\newcommand{\bu}{\boldsymbol{u}}

\newcommand{\bC}{\boldsymbol{C}}

\newcommand{\bz}{\boldsymbol{z}}

\newcommand{\bv}{\boldsymbol{v}}

\newcommand{\argmax}{\mathop{\mathrm{argmax}}}

\newcommand{\field}[1]{\mathbb{#1}}

\newcommand{\R}{\field{R}}

\newcommand{\E}{\field{E}}
\newcommand{\Var}{\mathrm{Var}}

\newcommand{\scI}{\mathcal{I}}

\DeclareMathOperator{\Regret}{Regret}

\DeclareMathOperator{\Wealth}{Wealth}

\DeclareMathOperator{\KL}{KL}
\DeclareMathOperator{\entropy}{H}
\DeclareMathOperator{\KLBern}{KL_\text{Bern}}

\newcommand{\indevent}{\mathbf{1}}

\title{Optimal Priors for Betting Tests with Bounded Payoffs:\\ An Exact Expected-Time Characterization}
\author{Francesco Orabona\\
King Abdullah University of Science and Technology (KAUST)\\
Thuwal, 23955-6900, Kingdom of Saudi Arabia\\
\texttt{\small francesco@orabona.com}}
\date{}

\begin{document}
\maketitle

\begin{abstract}
Testing by betting is a powerful framework that recasts hypothesis testing as a betting game.
A popular approach spreads the initial wealth over fixed bets according to a prior distribution, using a universal portfolio algorithm. For suitable priors, this approach achieves the optimal expected rejection time as the error probability goes to zero at a fixed alternative.
Yet, priors with this guarantee can behave quite differently near the null.
This motivates the question: Which prior should we use in testing-by-betting algorithms? We study bounded scalar payoffs and, under technical conditions, prove a sharp formula for expected rejection time, uniform over all bounded alternatives as the error probability tends to zero. The formula shows how the prior's mass around the zero bet determines the cost of adapting to an unknown alternative. We characterize the priors that attain the best possible leading constant when the error probability and alternative vary together. We also show that no valid test can uniformly improve on this cost across ranges of alternatives near the null. In particular, the uniform prior can require arbitrarily more observations than suitable priors that concentrate around the zero bet. A general inversion argument converts our results into sharp expected times to reach a prescribed confidence-sequence width. We also give a concrete prior choice and an FFT-based algorithm with almost linear computational complexity.
\end{abstract}

\section{Introduction}
\label{sec:intro}

Suppose a hypothesis lets us place bets that are fair whenever it is true. The observations need not themselves be the quantity of interest: a bet may be based on a bounded score, a difference between two measurements, or an e-value~\citep{VovkW21} supplied by another procedure. Our wealth over time represents the evidence against that hypothesis: if we gain enough money, we can reject the hypothesis while keeping the chance of a false alarm small. This is the testing-by-betting framework~\citep{ShaferV05,Shafer21,RamdasW25}.

Here, we focus on the practically relevant setting of scalar bounded payoffs. Write the payoff as $X_i-\mu_0$, where $X_i\in[0,1]$ and $\mu_0\in(0,1)$.
Under the null, we assume each payoff has mean zero given the previous observations. Every fixed betting fraction $\beta\in[-1/(1-\mu_0),1/\mu_0]$ then gives wealth
\begin{equation}
\label{eq:fixed-bet-wealth}
\Wealth_n(\beta)
:=\prod_{i=1}^n(1+\beta (X_i-\mu_0))\,.
\end{equation}
This is a nonnegative martingale under the null, so Ville's inequality guarantees that rejection at wealth $1/\alpha$ has false-alarm probability at most $\alpha$~\citep{Ville39}.

A core issue of this approach is choosing a suitable betting strategy. In fact, a strategy that works well for one alternative could work poorly for another, and we do not know the alternative in advance. A natural solution is to spread our initial wealth over fixed bets according to a prior $F$:
\[
\Wealth_n^F:=\int \! \Wealth_n(\beta)\,\mathrm{d}F(\beta)\,.
\]
In the online learning literature, this is known as an $F$-weighted portfolio for the case of two stocks~\citep{Cover91,CoverO96,OrabonaJ21}. Using a prior $F$ gives a valid test without choosing a single fixed bet, but leaves us with another choice: which prior should we use?

We evaluate a test by its expected rejection time under i.i.d.\ observations from an alternative law $Q$ with $\E_Q[X]\neq\mu_0$. Write $\tau_\alpha$ for the rejection time of a test with false-alarm probability at most $\alpha$. Recent optimality results by \citet{Waudby-SmithSJ25} and \citet{AgrawalRR26} provide partial answers to the choice of prior. In the bounded-mean setting, both give sufficient conditions for tests to attain
\begin{equation}
\label{eq:fixed-alternative-time}
\E_Q[\tau_\alpha]\sim\frac{\ln(1/\alpha)}{\lstar_Q}
\quad\text{as }\alpha\to0, \quad Q\in\cQ\text{ fixed},
\end{equation}
where $\beta_Q^\star$ is the fixed bet that maximizes the expected log wealth gained per observation, and $\lstar_Q:=\E_Q[\ln(1+\beta_Q^\star(X-\mu_0))]$ is its growth rate.
\citet{AgrawalRR26} also study alternatives approaching the null: when the limiting null law has positive variance, they prove matching-order expected-time bounds at fixed error probability.

However, these results still leave the choice of the prior unresolved. In fact, the fixed-alternative limit in \eqref{eq:fixed-alternative-time} can hold for priors that behave quite differently near the null, while the matching-order near-null bound in \citet{AgrawalRR26} still leaves the sharp prior-dependent constant unspecified. The prior's mass around the zero bet matters: Robbins-type priors give Law-of-the-Iterated-Logarithm (LIL) boundaries, whereas smooth priors typically incur a larger logarithmic cost~\citep{Robbins70, Orabona26b, OrabonaJ21}. So, the core question of this paper is \emph{which priors attain the optimal leading constant uniformly when the error probability and alternative vary together}.

Under technical conditions, we prove a sharp formula for $\E_Q[\tau_\alpha^F]$, the expected rejection time of a portfolio test with prior $F$ with false-alarm probability at most $\alpha$, uniform over all bounded alternatives as $\alpha\to0$ (Theorem~\ref{thm:pointwise}). The prior enters through the amount of mass it places around the zero bet, separately on the positive and negative sides. The result allows both alternatives with vanishing variance and oracle bets at the endpoints of the admissible interval.

For suitable Robbins-type priors, the expected rejection time is asymptotic, uniformly over all alternatives as $\alpha\to0$, to
\begin{equation}
\label{eq:benchmark}
B_\alpha(Q):=\frac{\ln\frac{1}{\alpha}+\LL(|\beta_Q^\star|)}{\lstar_Q},
\end{equation}
where $\LL(x):=\ln\max\{1,\ln(1/x)\}$ for $x>0$. The first term in the numerator of \eqref{eq:benchmark} captures the cost of requesting a small error probability. The second captures the cost of adapting to the unknown alternative. Our result identifies the oracle bet $|\beta_Q^\star|$ as the \emph{correct} quantity governing this adaptation cost. This is different from prior work that used $\lstar_Q$~\citep{Farrell64, RobbinsS74, AgrawalRR26}, but not in contradiction because near a null law with positive variance we have $\LL(|\beta_Q^\star|)=\LL(\lstar_Q)+O(1)$.

Under our regularity conditions, we give a necessary and sufficient condition on the prior's density near zero for attaining the best possible worst-case ratio to $B_\alpha(Q)$ (Corollary~\ref{thm:main}). This optimal ratio is 1. Theorem~\ref{thm:joint-minimax} establishes the full benchmark as the sharp minimax cost over ranges of small oracle bets, among all valid sequential tests. In contrast, Beta priors, including the uniform and Krichevsky–Trofimov priors~\citep{KrichevskyT81}, can require an arbitrarily large multiple of $B_\alpha(Q)$ when the alternative approaches the null as $\alpha\to0$.

A consequence of our result is that the exact $F$-weighted portfolio with the combined prior described in \citet[Remark~2]{OrabonaJ21} attains $B_\alpha(Q)$ uniformly. We also show how priors can complement one another: one can supply enough mass near small positive bets and another near small negative bets. Their combination can attain $B_\alpha(Q)$ uniformly even when each component has an unbounded worst-case ratio (Corollary~\ref{cor:combine}).

Then, Section~\ref{sec:confidence-sequences} combines our guarantees on portfolio testing to obtain sharp expected times to reach a prescribed confidence-sequence width using portfolio-based approaches.

Finally, Section~\ref{sec:computation} gives a concrete prior choice and an FFT-based algorithm for exact continuous-mixture evaluation given prior moments, with nearly linear computational cost in the number of observations.

\section{Related work}
\label{sec:related}

The connection between betting and probability goes back to \citet{Ville39}, whose martingale inequality supplies the validity guarantee used here. The oracle betting strategy when the payoff distribution is known is Kelly’s criterion, which maximizes the expected log wealth~\citep{Kelly56,Breiman61}. When we do not know this probability distribution, universal portfolio algorithms show how to compete with the best constant allocation in hindsight~\citep{Cover91,CoverO96}. The same betting viewpoint later became a tool for parameter-free online learning~\citep{OrabonaP16,CutkoskyO18}. These developments connect probabilistic validity, stochastic growth, and adversarial regret, but each supplies a different guarantee.

\citet{ShaferV05} develop a game-theoretic foundation for probability: forecasts specify betting prices, and nonnegative capital processes test their agreement with subsequent outcomes. Their bounded forecasting protocol (Sections~3.2–3.3) is exactly the betting setup we use here.
Modern accounts of testing by betting include \citet{Shafer21, RamdasW25}. The log-growth viewpoint underlies growth-optimal e-values~\citep{GrunwaldHK24,LarssonRR25} and e-processes~\citep{WangWZ25,Waudby-SmithSJ25}. For bounded means, \citet{Waudby-SmithR21} characterize test martingales under the constant-conditional-mean null as predictable betting products, and \citet{Clerico25} identifies these products as the admissible e-processes. This resolves which wealth processes are undominated, but it leaves open which betting strategies have the best expected rejection times.

\paragraph{From betting guarantees to concentration and prior choice.}
The construction of priors for sequential LIL boundaries was developed by \citet{DarlingR67, Robbins70}. The game-theoretic development in \citet[Section~5.2, especially equation~(5.15)]{ShaferV05} similarly combines nonnegative capital processes across betting scales to obtain a LIL bound. They credit the earlier proofs of \citet{DawidV99} and \citet{Vovk88}, and trace the underlying idea back to \citet{Ville39}. These pathwise guarantees motivate adaptation across scales, but do not characterize the expected rejection time associated with a prior. An early connection to online betting appears in the technical report of \citet{Orabona26b}. His proof method works by localizing a log-squared prior in a relative interval of bets, and then turning the resulting wealth guarantee into an empirical Bernstein LIL bound for bounded martingale differences. We will use his localization argument in our expected-time proof.

\citet[Section~7.2]{JunO19} show how noisy coin-betting potentials give time-uniform concentration through a nonnegative supermartingale. Their log-squared prior gives LIL-rate bounds for sub-Gaussian martingale differences, and their online-learning reduction extends the construction to Banach spaces. Building on this connection, \citet{OrabonaJ21} derive bounded-mean confidence sequences from universal-portfolio regret and obtain the sharp LIL leading constant with a prior adapted from \citet{Robbins70}. The mixture proposed in \citet[Remark~2]{OrabonaJ21} already combines logarithmic worst-case regret, LIL behavior, and optimal i.i.d.\ alternative growth. This idea is developed in \citet{AgrawalR26}, who give explicit regret bounds and a sharp asymptotic regret bound for their full-support Robbins prior on paths with bounded mixture wealth and diverging intrinsic variance.
\citet{RyuB22} subsequently develop efficient approximations to universal-portfolio wealth.

The relation between prior concentration and boundary growth has itself been studied systematically. For mixtures of constant-proportion product strategies, \citet{SatoMT18} give Robbins-type priors yielding the sharp self-normalized LIL and establish a correspondence between classes of prior densities near zero and upper-class boundary functions. Their results apply to bounded forecasting and concern pathwise boundaries. Instead, our criterion classifies the joint expected rejection-time cost of full-support priors. For exponential mixtures, \citet{HowardRMS21} show that a density positive and continuous at zero produces a boundary of order $\sqrt{v\ln v}$ in intrinsic time $v$, explaining a related penalty for smooth priors.

Finite-time LIL concentration and anti-concentration bounds already capture the transition between confidence and adaptation costs~\citep{Balsubramani15}. \citet{KaufmannK21} combine confidence and iterated-logarithm time corrections in mixture bounds for exponential families, including Bernoulli observations. Related developments include LIL-based bandit bounds~\citep{JamiesonMNB14} and almost-sure iterated-logarithm regret~\citep{AgrawalR25LIL}. The finer logarithmic terms have classical antecedents in upper-class integral tests~\citep{Erdos42, Feller43}. Priors are also connected to coding: describing the logarithmic betting scale incurs an iterated-logarithm cost~\citep{Rissanen83}, and the prior governs asymptotics under prior averaging more generally~\citep{ClarkeB90}.

\paragraph{Expected-time optimality.}
The small-error and near-null regimes also have a history of their own. \citet{Farrell64} proves an iterated-logarithm lower bound near the null in regular parametric models, and \citet{RobbinsS74} study power-one tests based on prior-weighted likelihood ratios. \citet{FellourisT12} study how to choose a mixture prior by expected sample size. They optimize a higher-order minimax criterion based on KL-weighted expected sample size for exponential-family alternatives in a compact interval bounded away from the null. Our criterion instead distinguishes priors through their leading cost as the alternative approaches the null. \citet[Lemma~4.1 and Remark~D.4]{ChenLQ17} give a lower bound combining confidence and the entropy of the unknown scale, up to constants, for sign testing, including a Bernoulli extension. Theorem~\ref{thm:joint-minimax} identifies the sharp sum in our bounded model, uniformly as both the error level and betting scales vanish.

Nonasymptotic expected-time bounds also combine the two costs. \citet{ShinRR21} analyze generalized-likelihood-ratio tests with iterated-logarithm boundaries, including bounded observations through their cumulant bounds; their information measure is determined by those bounds rather than by the full alternative law. For bounded adapted observations, \citet{DurandW25} derive power and expected rejection-time guarantees from regret, allowing conditional means to change and vanish over time. That setting differs from a sequence of i.i.d.\ alternative laws approaching a fixed null, which we consider.

For general composite hypotheses, \citet{AgrawalRR26} establish the small-error information lower bound and a nonparametric Farrell-type lower bound near the null, together with sufficient conditions for matching upper bounds. Their bounded-mean specialization uses the same KL distance to the null as our oracle rate, and their Proposition~5.1 combines tests attaining the two guarantees in their respective limits.
\citet{Waudby-SmithSJ25} obtain the sharp fixed-alternative expected-time limit from a deterministic sublinear portfolio-regret bound and a moment condition. Section~\ref{sec:expected-time-comparison} compares these expected-time guarantees with our uniform formula, including their assumptions and orders of limits. Related information bounds arise in best-arm identification~\citep{GarivierK16}; the bounded-support divergence studied by \citet{HondaT10} for cumulative bandit regret coincides with our oracle log wealth.

Expected-time criteria also arise for confidence sequences, which can be obtained by inverting betting tests~\citep{Hendriks18, OrabonaJ21, Waudby-SmithR21}. \citet{ShekharR23near} define effective width through the expected time to reach a prescribed width, give an information lower bound, and bound the width of betting confidence sequences with logarithmic regret when the observation law is supported strictly inside $[0,1]$. Effective width concerns two endpoints rather than a specified null mean. In Section~\ref{sec:confidence-sequences}, we obtain a sharp formula for the expected time to reach a prescribed width, improving over their result.

\section{Setting}
\label{sec:setting}

We consider bounded scalar payoffs. Without loss of generality, we will assume them to be in $[0,1]$.
Hence, let $X_1,X_2,\ldots\in[0,1]$ be i.i.d. and fix $\mu_0\in(0,1)$.
The bets that guarantee nonnegative wealth form the interval $[-1/(1-\mu_0),1/\mu_0]$, so $0\leq1+\beta (X_i-\mu_0)\leq \max\{1/\mu_0,1/(1-\mu_0)\}$. The null and alternative classes are
\[
\cP_0:=\{P:\E_P[X]=\mu_0\},\qquad
\cQ:=\{Q:\E_Q[X]\neq\mu_0\},
\]
where all distributions are supported on $[0,1]$.

We aim to produce a \emph{level-$\alpha$ family}, that is, a family $\tau=(\tau_\alpha)_{\alpha\in(0,1/2]}$ of stopping times with respect to the natural filtration such that $\sup_{P\in\cP_0} \ \Pr_{P}\{\tau_\alpha<\infty\}\leq \alpha$ for all $\alpha$. Here, stopping means rejecting the null hypothesis.

\paragraph{Portfolio tests.}
We follow the bounded forecasting protocol of \citet[Sections~3.2–3.3]{ShaferV05} to define a wealth process of a fixed bet $\beta$ as in \eqref{eq:fixed-bet-wealth}.
Given a probability prior $F$ on $[-1/(1-\mu_0),1/\mu_0]$, the $F$-weighted portfolio has portfolio wealth
\begin{equation}
\label{eq:portfolio-wealth}
\Wealth_n^F
:=\int_{-1/(1-\mu_0)}^{1/\mu_0} \Wealth_n(\beta)\,\mathrm{d}F(\beta)\,.
\end{equation}
The portfolio wealth in \eqref{eq:portfolio-wealth} is generated by the online\footnote{In the statistics literature, people use ``predictable'' instead of ``online'', but we prefer ``online'' because universal portfolio is an online algorithm.} bets
\[
\beta_n:=\frac{\int_{-1/(1-\mu_0)}^{1/\mu_0}\beta \Wealth_{n-1}(\beta)\,\mathrm{d}F(\beta)}{\Wealth_{n-1}^F},
\]
with $\beta_n:=0$ when $\Wealth_{n-1}^F=0$. To see the portfolio interpretation introduced by \citet{OrabonaJ21}, set $p:=\mu_0+\beta\mu_0(1-\mu_0)\in[0,1]$. Then,
\begin{equation}
\label{eq:two-stock-returns}
1+\beta (X_i-\mu_0)=p\frac{X_i}{\mu_0}+(1-p)\frac{1-X_i}{1-\mu_0}\,.
\end{equation}
Equation~\eqref{eq:two-stock-returns} identifies $\Wealth^F$ as a portfolio selection problem on these two stocks. Under every $P\in\cP_0$, it is a nonnegative martingale starting from 1. Ville's inequality gives
\begin{equation}
\label{eq:portfolio-test}
\tau_\alpha^F:=\inf\{n\geq1:\Wealth_n^F\geq1/\alpha\},\qquad
\sup_{P\in\cP_0} \ \Pr_{P}\{\tau_\alpha^F<\infty\}\leq\alpha\,.
\end{equation}

\paragraph{The log-optimal bettor.} For $Q\in\cQ$, define
\[
\ell_Q(\beta):=\E_Q[\ln(1+\beta (X-\mu_0))] \quad \text{ and } \quad
\lstar_Q:=\max_{\beta\in[-1/(1-\mu_0),1/\mu_0]} \ \ell_Q(\beta),
\]
and let $\beta^\star_Q$ be the maximizer. Upper semicontinuity on the compact interval $[-1/(1-\mu_0),1/\mu_0]$ gives existence, strict concavity gives uniqueness, and $\ell_Q'(0)=\E_Q[X-\mu_0]\neq0$ gives $\lstar_Q>0$: a sufficiently small bet with the same sign as $\E_Q[X-\mu_0]$ improves on $\ell_Q(0)=0$. Also,
\begin{equation}
\label{eq:oracle-growth-bound}
0<\lstar_Q\leq\min\{|\beta^\star_Q|,\ln \max\{1/\mu_0,1/(1-\mu_0)\}\},
\end{equation}
since $\ln(1+\beta (X-\mu_0))\leq\beta (X-\mu_0)\leq|\beta|$ and $1+\beta (X-\mu_0)\leq \max\{1/\mu_0,1/(1-\mu_0)\}$. Lemma~\ref{lemma:proj} identifies $\lstar_Q$ with the infimum Kullback--Leibler (KL) divergence from $Q$ to the null, the same quantity studied by \citet{HondaT10}. We allow $\beta^\star_Q$ to take endpoint values, as well as alternatives converging to the point mass at $\mu_0$.

\section{Main results}
\label{sec:main}

Here, we introduce our main result: The prior's behavior around the zero bet determines the leading cost of adaptation. All the proofs of the results of this section are deferred to Appendix~\ref{sec:proofs}.

\subsection{Uniform expected rejection times}

First, we define a logarithmic decay index for the prior around the zero bet. This will quantify the prior's contribution to the expected rejection time.
\begin{definition}[Logarithmic decay index]
\label{def:prior}
Let $F$ have a density $f$ on $[-1/(1-\mu_0),1/\mu_0]$ that is bounded below by a positive constant on each compact subset of $(-1/(1-\mu_0),1/\mu_0)\setminus\{0\}$.  Moreover, letting $h(\beta):=|\beta|f(\beta)$, assume that the limits
\begin{equation}
\label{eq:decay-indices}
a_\pm(F):=\lim_{\beta\to0^\pm} \ \frac{\ln(1/h(\beta))}{\ln\ln(1/|\beta|)}
\end{equation}
exist in $[-\infty,+\infty]$. We call $a(F):=\max\{a_+(F),a_-(F)\}$ the \emph{logarithmic decay index} of the prior $F$.
\end{definition}

On each side of zero, $h$ is the prior density in the coordinate $\ln(1/|\beta|)$. Integrability forces both limits in \eqref{eq:decay-indices} to be at least 1.

Our main result identifies the leading expected time at each alternative, uniformly as the error probability tends to zero.
\begin{theorem}
\label{thm:pointwise}
Fix $\mu_0\in(0,1)$ and a prior $F$ satisfying Definition~\ref{def:prior} with $a(F)<\infty$. For each $Q\in\cQ$, define
\begin{equation}
\label{eq:pointwise-benchmark}
a_Q(F):=
\begin{cases}
a_+(F),&\beta_Q^\star>0,\\
a_-(F),&\beta_Q^\star<0,
\end{cases}
\qquad
\widetilde B_\alpha^F(Q):=
\frac{\ln(1/\alpha)+a_Q(F)\LL(|\beta_Q^\star|)}{\lstar_Q}\,.
\end{equation}
Then the normalized rejection time converges uniformly in $L^1$:
\begin{equation}
\label{eq:pointwise-l1}
\lim_{\alpha\downarrow0}\ \sup_{Q\in\cQ}\ \E_Q\left[\left|\frac{\tau_\alpha^F}{\widetilde B_\alpha^F(Q)}-1\right|\right]=0~.
\end{equation}
In particular, for every $\epsilon>0$ there is $\alpha_0>0$ such that
$(1-\epsilon)\widetilde B_\alpha^F(Q)\leq\E_Q[\tau_\alpha^F]\leq
(1+\epsilon)\widetilde B_\alpha^F(Q)$ for all $Q\in\cQ$ and $0<\alpha\leq\alpha_0$.
\end{theorem}
\begin{proof}[Proof sketch]
For the upper bound, we localize the prior to bets $t\beta_Q^\star$ with $1-2\delta\leq t\leq1-\delta$, for fixed
$\delta\in(0,1/2)$, following \citet{Orabona26b,OrabonaJ21}.
Concavity preserves at least $(1-2\delta)\lstar_Q$ of the oracle growth. Jensen's inequality turns the reciprocal prior mass of this interval into an additive logarithmic penalty.
The decay index bounds this penalty by $(a_Q(F)+\zeta)\LL(|\beta_Q^\star|)$ plus a constant independent
of $Q,\alpha$, for any fixed $\zeta>0$.
A bounded-stopping argument then converts positive drift into an expected-time bound, including at endpoint oracle bets (Lemmas~\ref{lemma:window}--\ref{lemma:mass}).

For the lower tail, small oracle bets require accounting for the prior penalty. In this case, a change of measure controls early rejection. Afterwards, concentration controls the empirical first and second moments. These estimates remain uniform when variance vanishes because, as $|\beta_Q^\star|\to0$, $\lstar_Q=\tfrac12(\beta_Q^\star)^2 \E_Q[(X-\mu_0)^2](1+o(1))$.
On the resulting event, only a small interval of bets can produce wealth above 1. We bound the excess wealth
$\Wealth_n^F-1$: a convex-chord inequality supplies a factor vanishing linearly at zero, which cancels the density's $1/|\beta|$ singularity. The resulting bound prevents rejection before almost $\widetilde B_\alpha^F(Q)$ observations, with probability tending to 1 uniformly (Lemma~\ref{lemma:pointwise-lower}). For oracle bets away from zero, Lemma~\ref{lemma:information-early} gives the analogous probability bound at the information scale. Together with the expected-time upper bound, these early-rejection bounds control the positive and negative parts of the normalized deviation, proving uniform $L^1$ convergence.
\end{proof}

Note that the alternative may vary arbitrarily with $\alpha$, that is, no variance lower bound or support margin is assumed. If $a(F)=1$, both directional indices equal 1. Thus, $\widetilde B_\alpha^F(Q)=B_\alpha(Q)$, and the theorem gives $\mathbb E_Q[\tau_\alpha^F]/B_\alpha(Q)\to1$ uniformly over $Q$.

\subsection{Comparison with earlier expected-time bounds}
\label{sec:expected-time-comparison}

The confidence cost is unavoidable: \citet[Proposition~3.2]{Waudby-SmithSJ25} and \citet[Theorem~3.1]{AgrawalRR26} show that for every $Q\in\cQ$ and $\alpha\in(0,1)$, any stopping time $\tau$ with $\sup_{P\in\cP_0}\Pr_P\{\tau<\infty\}\leq\alpha$ satisfies the information lower bound
\begin{equation}
\label{eq:information-lower}
\E_Q[\tau]\geq\frac{\ln(1/\alpha)}{\lstar_Q}~.
\end{equation}
For a fixed alternative, Theorem~\ref{thm:pointwise} attains this lower bound asymptotically because the adaptation term is constant. Its uniform formula also covers alternatives approaching the null, where that term can dominate.

Finer distributional limits are available at a fixed alternative. Under their Assumption~4.1, \citet[Theorem~4.4]{DeepA26} prove that the stopping time $\tau_\alpha$ of their test based on the empirical KL divergence to the null satisfies, for each fixed $Q\in\cQ$,
\[
\frac{\tau_\alpha-\ln(1/\alpha)/\lstar_Q}{\sqrt{\ln(1/\alpha)}}
\xrightarrow{\mathrm d}
\mathcal N\!\left(0,\frac{\Var_Q[\ln(1+\beta_Q^\star(X-\mu_0))]}{(\lstar_Q)^3}\right)
\qquad(\alpha\downarrow0)~.
\]
Their proof first establishes a central limit theorem for the empirical KL statistic (Theorem~4.2), then transfers it to the stopping time. This describes fluctuations around the fixed-alternative benchmark on the scale $\sqrt{\ln(1/\alpha)}$. Our Theorem~\ref{thm:pointwise} gives a complementary guarantee for prior mixtures: convergence in $L^1$ relative to a prior-dependent benchmark, uniformly as the alternative varies with $\alpha$. It identifies the adaptation cost near the null, but does not give a Gaussian fluctuation limit.

Theorem~5.3 of \citet{Waudby-SmithSJ25} gives a uniform expected-time bound when $\lstar_Q\geq\underline\ell>0$ throughout the alternative class. Our Theorem~\ref{thm:pointwise} removes this growth-rate restriction in our bounded model and identifies the sharp prior-dependent adaptation cost.

For each fixed $\alpha\in(0,1/2)$, Corollary~4.8 of \citet{AgrawalRR26} constructs a test $\tau_\alpha$ such that, along every sequence $Q_k$ of bounded alternatives with means below $\mu_0$ converging weakly to $P\in\cP_0$ with $\Var_P(X)>0$,
\[
\limsup_{k\to\infty} \ \frac{\lstar_{Q_k}\E_{Q_k}[\tau_\alpha]}{\LL(\lstar_{Q_k})}<\infty\,.
\]
Their nonasymptotic bound in equation~(26) already contains both confidence and iterated-logarithm gap costs, but with generally nonsharp constants.

To compare their growth-rate scale with the oracle-bet scale in Theorem~\ref{thm:pointwise}, we first relate the two under nondegeneracy. The proof is in Appendix~\ref{sec:proof-near-null-scales}.
\begin{proposition}
\label{prop:near-null-scales}
Let $Q_k\in\cQ$ satisfy $\E_{Q_k}[X]\to\mu_0$ and $\inf_k\E_{Q_k}[(X-\mu_0)^2]>0$. Then, $|\beta_{Q_k}^\star|\to0$ and
\[
\LL(\lstar_{Q_k})=\LL(|\beta_{Q_k}^\star|)+\ln2+o(1)\,.
\]
\end{proposition}

Our upper bound also gives an explicit leading coefficient in their fixed-error near-null regime. Let $F$ satisfy Definition~\ref{def:prior} with $a(F)=1$, fix $\alpha\in(0,1/2)$, and let $Q_k\in\cQ$ converge weakly to $P\in\cP_0$ with $\Var_P(X)>0$. Then,
\[
\limsup_{k\to\infty} \ \frac{\lstar_{Q_k}\E_{Q_k}[\tau_\alpha^F]}{\LL(\lstar_{Q_k})}\leq1\,.
\]
Indeed, Proposition~\ref{prop:near-null-scales} gives $\LL(|\beta_{Q_k}^\star|)/\LL(\lstar_{Q_k})\to1$. Divide \eqref{eq:pointwise-upper} by $\LL(\lstar_{Q_k})/\lstar_{Q_k}$, send $k\to\infty$, and then send $\delta,\zeta\to0$. This gives an upper bound at fixed $\alpha$.
When $\alpha_k\to0$, Theorem~\ref{thm:pointwise} gives a matching formula for every prior satisfying Definition~\ref{def:prior} with finite indices. By Proposition~\ref{prop:near-null-scales}, we can replace $\LL(|\beta_{Q_k}^\star|)$ by $\LL(\lstar_{Q_k})$ along these nondegenerate sequences.

When variance vanishes, the distinction between the two scales can be a difference in order. For example, fix the error level $\alpha\in(0,1/2)$ and consider the deterministic alternatives
\[
 Q_\Delta:=\delta_{\mu_0+\Delta},\qquad 0<\Delta<1-\mu_0,
\]
where $\delta_x$ denotes the point mass at $x$. These alternatives converge to the null law $\delta_{\mu_0}$ as $\Delta\downarrow0$. So, we have
\[
\beta_{Q_\Delta}^\star=\frac1{\mu_0},\qquad
\lstar_{Q_\Delta}=\ln(1+\Delta/\mu_0)\sim\frac{\Delta}{\mu_0}
\qquad(\Delta\downarrow0)~.
\]
Thus, $\LL(|\beta_{Q_\Delta}^\star|)=0$, and for any fixed prior $F$ satisfying Definition~\ref{def:prior} with finite indices, the upper bound \eqref{eq:pointwise-upper} and the information lower bound \eqref{eq:information-lower} show that the expected time has order $1/\lstar_{Q_\Delta}$. Consequently, as $\Delta\downarrow0$,
\[
\E_{Q_\Delta}[\tau_\alpha^F]=\Theta\!\left(\frac1\Delta\right),\qquad
\frac{\LL(\lstar_{Q_\Delta})}{\lstar_{Q_\Delta}}
=\Theta\!\left(\frac{\ln\ln(1/\Delta)}{\Delta}\right),
\]
where the implicit constants may depend on $\alpha,F,\mu_0$, but not on $\Delta$. The growth-rate scale therefore introduces an extra factor of order $\ln\ln(1/\Delta)$. That said, this example lies outside the positive-limiting-variance assumption of \citet[Corollary~4.8]{AgrawalRR26}.

\subsection{Worst-case efficiency}

To summarize performance, we compare expected rejection time with the benchmark in \eqref{eq:benchmark}.
\begin{definition}[Efficiency]
\label{def:eff}
Define the \emph{efficiency} of a level-$\alpha$ family $\tau$ as
\begin{equation}
\label{eq:efficiency}
\eff(\tau):=\lim_{\alpha_0\downarrow0}\ \sup\left\{\frac{\E_Q[\tau_\alpha]}{B_\alpha(Q)} : 0<\alpha\leq\alpha_0,\ Q\in\cQ\right\}\,.
\end{equation}
\end{definition}
Note that the limit in \eqref{eq:efficiency} exists because the supremum decreases as $\alpha_0\downarrow0$.

To isolate the cost of adapting to an unknown betting scale, we compare tests over a range of small oracle bets. For $R\geq1$, the interval $[e^{-2R},e^{-R}]$ has length $R$ in the coordinate $\ln(1/|\beta|)$, so it contains order $R$ bets separated by a factor of 4 as $R$ grows. Note that, throughout this interval, $\LL(|\beta|)$ differs from $\ln R$ by at most $\ln2$.

The next theorem shows that this adaptation cost is unavoidable even if the test knows $R$.
\begin{theorem}
\label{thm:joint-minimax}
For $R\geq1$ and $0<\alpha\leq1/2$, define
\begin{equation}
\label{eq:oracle-range}
\cQ_R:=\{Q\in\cQ:e^{-2R}\leq|\beta_Q^\star|\leq e^{-R}\},\qquad
\mathcal T_\alpha:=\left\{\tau:\sup_{P\in\cP_0}\Pr_P\{\tau<\infty\}\leq\alpha\right\},
\end{equation}
where $\tau$ ranges over stopping times with values in $\{1,2,\ldots\}\cup\{\infty\}$, optionally using an independent random seed with the same law under every hypothesis.
Then, for each fixed $\mu_0\in(0,1)$, we have
\begin{equation}
\label{eq:joint-minimax}
\lim_{\substack{\alpha\downarrow0\\R\to\infty}}
\inf_{\tau\in\mathcal T_\alpha}\ \sup_{Q\in\cQ_R}
\frac{\E_Q[\tau]}{B_\alpha(Q)}=1~.
\end{equation}
That is, for every $\epsilon>0$, there are $\alpha_0\in(0,1/2]$ and $R_0\geq1$ such that the infimum--supremum in \eqref{eq:joint-minimax} lies in $[1-\epsilon,1+\epsilon]$ for every $0<\alpha\leq\alpha_0$ and every $R\geq R_0$. The test in the infimum may depend on both $\alpha$ and $R$.
\end{theorem}

\begin{proof}[Proof sketch]
We pack the range with $m$ Bernoulli alternatives $Q_1,\ldots,Q_m$ whose oracle bets are separated by a factor of 4, with $m$ of order $R$. Let $P$ be their common Bernoulli null, and let $Q_j^\tau$ and $P^\tau$ be the laws of the stopped history $Z$, including any independent random seed and a symbol for nonrejection. Choose the generating index $J$ uniformly and define the mixture law by
\[
\Pr\{J=j\}:=1/m\quad(1\leq j\leq m),\qquad
\overline Q^\tau:=\frac1m\sum_{j=1}^m Q_j^\tau~.
\]
If the expected times are finite, the stopped likelihood identity and the KL chain rule give
\[
\frac1m\sum_{j=1}^m\lstar_{Q_j}\E_{Q_j}[\tau]
=\KL(\overline Q^\tau\|P^\tau)+\operatorname{I}(J;Z)~,
\]
where $\operatorname{I}(J;Z)$ is the mutual information between the alternative index and the stopped history. The first term is at least $\ln(1/\alpha)$: rejection occurs almost surely under each $Q_j$, but with probability at most $\alpha$ under $P$. Moreover, validity makes early rejection unlikely, and by the time of rejection the empirical mean identifies the separated scale with error tending to zero as $\alpha\to0$, uniformly in $m$ and the scales. Fano's inequality~\citep{CoverT06} therefore bounds the second term below by $(1-o(1))\ln m-\ln2$. Lemma~\ref{lemma:scale-information} gives the finite bound. Since $\ln m=\ln R+O(1)$, this proves the lower bound; an infinite expected time makes it immediate. Theorem~\ref{thm:pointwise} supplies the upper bound using any fixed prior with index 1.
\end{proof}

Theorem~\ref{thm:joint-minimax} implies $\eff(\tau)\geq1$ for every level-$\alpha$ family. Thus, efficiency 1 means that, for every $\epsilon>0$, there exists $\alpha_0>0$ such that
\[
\E_Q[\tau_\alpha]\leq(1+\epsilon)B_\alpha(Q)
\qquad\text{for all }Q\in\mathcal Q
\text{ and } 0<\alpha\leq\alpha_0\,.
\]

For finite indices, the following characterization follows from Theorem~\ref{thm:pointwise}. For infinite indices, the lower bound in Lemma~\ref{lemma:pointwise-lower} gives the same conclusion.
\begin{corollary}
\label{thm:main}
For any prior $F$ satisfying Definition~\ref{def:prior} and the full bounded alternative class $\cQ$, we have
\[
\eff\left(\tau^F\right)=a(F)\,.
\]
\end{corollary}
\begin{proof}[Proof sketch]
For finite $a(F)$, Theorem~\ref{thm:pointwise} gives the upper bound because $\widetilde B_\alpha^F(Q)/B_\alpha(Q)\leq a(F)$.
To attain this ratio, choose the side with index $a(F)$, taking the positive side without loss of generality. For
sufficiently small $\Delta>0$, define
\[
Q_\Delta:=\mathrm{Ber}(\mu_0+\Delta),\qquad
\alpha_\Delta:=1/\LL(\beta_{Q_\Delta}^\star)~.
\]
Then, $\beta_{Q_\Delta}^\star\to0$, $\alpha_\Delta\to0$, and $\ln(1/\alpha_\Delta)=o(\LL(\beta_{Q_\Delta}^\star))$. Thus, the adaptation term dominates, and uniformity in Theorem~\ref{thm:pointwise} gives
$\E_{Q_\Delta}[\tau_{\alpha_\Delta}^F]/B_{\alpha_\Delta}(Q_\Delta)\to a(F)$.
For an infinite index, Lemma~\ref{lemma:pointwise-lower} gives arbitrarily large limiting lower bounds by allowing any finite exponent below that index.
\end{proof}

Together with Theorem~\ref{thm:joint-minimax}, Corollary~\ref{thm:main} says that a portfolio test with a prior that satisfies the conditions of Definition~\ref{def:prior} is optimal if and only if $a(F)=1$, equivalently, when the log-scale prior density decays like $(\ln\frac{1}{|\beta|})^{-1-o(1)}$. Integrability forbids $(\ln\frac{1}{|\beta|})^{-1}$ itself, so optimal priors sit at the border of integrability.

It is worth remarking that one can construct priors with efficiency 1 for which a one-sided decay-index limit does not exist, so the assumptions of Theorem~\ref{thm:pointwise} would not hold.

Let’s now see some examples.
\begin{corollary}
\label{cor:examples}
Let $L_1(\beta):=\ln\frac{e\max\{1/\mu_0,1/(1-\mu_0)\}}{|\beta|}\geq1$ and $L_2(\beta):=\ln(eL_1(\beta))\geq 1$ for $\beta\in[-1/(1-\mu_0),1/\mu_0]\setminus\{0\}$.
\begin{enumerate}[label=(\roman*),leftmargin=*,itemsep=1pt,topsep=2pt]
\item If $p=\mu_0+\beta\mu_0(1-\mu_0)$ has a $\mathrm{Beta}(\nu_1,\nu_2)$ prior, $\nu_1,\nu_2>0$, then $\eff(\tau^F)=\infty$. This includes the Beta$(1/2,1/2)$ prior $\nu_1=\nu_2=\frac12$ and the uniform prior~\citep{KrichevskyT81, CoverO96}.
\item If $f(\beta)\propto 1/(|\beta|L_1(\beta)^a)$ with $a>1$, then $\eff(\tau^F)=a$.
\item If $f_\eta(\beta):= C_\eta/(|\beta|L_1(\beta)L_2(\beta)^{1+\eta})$ with $\eta>0$ and $C_\eta$ a normalizing constant, then $\eff(\tau^{F_\eta})=1$.
\end{enumerate}
\end{corollary}
Note that the prior $f_\eta$ is a proper probability density: the change of variable $w:=L_2(\beta)$ gives
\[
\int_{0<|\beta|\leq \max\{1/\mu_0,1/(1-\mu_0)\}} \frac{\mathrm{d}\beta}{|\beta|L_1(\beta)L_2(\beta)^{1+\eta}}=\frac{2}{\eta}\,.
\]
It is the analogue for betting of the priors used by \citet{Robbins70} to obtain LIL boundaries.

\subsection{Comparison with growth and LIL guarantees}

If $\eff(\tau)=1$, \eqref{eq:information-lower} gives the fixed-alternative limit \eqref{eq:fixed-alternative-time}. The converse fails: $F$-weighted portfolios with Beta priors satisfy that limit but have infinite efficiency by Corollary~\ref{cor:examples}(i). Thus, efficiency 1 is strictly stronger than fixed-alternative expected-time optimality in this bounded model.

Sublinear portfolio regret gives oracle growth at each fixed alternative by \citet[Theorem~2.3]{Waudby-SmithSJ25}. Their Theorem~5.2 also gives distribution-uniform comparisons of log wealth. However, these guarantees do not determine the leading expected time along changing alternatives: the Beta(1/2,1/2) prior has logarithmic regret but infinite efficiency by Corollary~\ref{cor:examples}(i).

Moreover, a LIL-order bound alone also does not determine the sharp expected-time constant. For example, the log-squared prior in \citet{Orabona26b} gives a betting-based empirical Bernstein LIL bound. Its logarithmic decay index is 2, so a full-support extension satisfying Definition~\ref{def:prior} and preserving its behavior near zero has efficiency 2.

\subsection{One-sided tests and combinations of priors}

We can also consider nonnegative bets to test a one-sided mean null. In this case, only the prior's behavior on the positive side matters.
\begin{corollary}
\label{cor:one-sided}
For the null $\E_P[X]\leq\mu_0$ and alternatives $\E_Q[X]>\mu_0$, let $F$ be a probability prior on $[0,1/\mu_0]$ with a density bounded below by a positive constant on each compact subset of $(0,1/\mu_0)$. Suppose its positive-side logarithmic decay index $a_+(F)$ exists. Define the oracle, benchmark, and efficiency using these null and alternative classes and nonnegative bets. Then, $\eff(\tau^F)=a_+(F)$, and every valid test family has efficiency at least 1. If $a_+(F)<\infty$, the uniform $L^1$ limit \eqref{eq:pointwise-l1} also holds over this one-sided alternative class, with $a_Q(F)=a_+(F)$.
\end{corollary}

An $F$-weighted portfolio can also use a convex combination of priors that are effective on different sides of the null.
\begin{corollary}
\label{cor:combine}
Let $F_1,\dots,F_k$ be distributions satisfying Definition~\ref{def:prior}, and let $w_1,\dots,w_k>0$ be fixed weights summing to one. For $F:=\sum_{j=1}^k w_jF_j$,
\[
\eff(\tau^F)=\max\left\{\min_{1\leq j\leq k}a_+(F_j),\min_{1\leq j\leq k}a_-(F_j)\right\}\,.
\]
In particular, $F$ has efficiency 1 if and only if each side has at least one component with logarithmic decay index 1 on that side. If both sidewise minima are finite, Theorem~\ref{thm:pointwise} also applies, with these minima as the directional coefficients.
\end{corollary}

For example, let $f_1$ be proportional to $f_\eta$ on positive bets and to a positive constant on negative bets, and reverse these choices for $f_2$. Their pairs of one-sided logarithmic decay indices are $(1,\infty)$ and $(\infty,1)$, so both component tests have infinite efficiency, while every $F$-weighted portfolio whose prior is a fixed convex combination of these two priors with positive weights has efficiency 1. That is, different components supply the efficient behavior on the two sides.

\subsection{Losses in oracle growth}

We now consider processes that lose oracle growth. A strict loss at even one alternative prevents a process from attaining the optimal efficiency of 1.

\begin{proposition}
\label{prop:form}
Let $\Wealth$ be a finite-valued nonnegative process adapted to the natural filtration, let $\tau_\alpha:=\inf\{n:\Wealth_n\geq1/\alpha\}$, and assume $(\tau_\alpha)$ is a level-$\alpha$ family. If for some $Q\in\cQ$ and $g<\lstar_Q$ we have $\limsup_{n} \ \frac1n\ln \Wealth_n\leq g$ $Q$-a.s., then $\eff(\tau)\geq\lstar_Q/g$ (interpreted as $+\infty$ if $g\leq 0$). In particular, the prior-weighted Hoeffding process~\citep{Hoeffding63,HowardRMS21}
\begin{equation}
\label{eq:hoeffding-process}
\Wealth^{\mathrm{H}}_n
:=\int\!\exp\Big(\beta\sum_{i=1}^n (X_i-\mu_0)-\frac{n\beta^2}{8}\Big)\,\mathrm{d}F(\beta)
\end{equation}
satisfies $\eff(\tau^{\mathrm{H}})=\infty$ over $\cQ$ for any probability prior $F$ on $\R$.
\end{proposition}

For online betting strategies constrained to a compact set $\mathcal{K}\subset(-1/(1-\mu_0),1/\mu_0)$, a strict gap follows whenever $\beta^\star_Q\notin\mathcal{K}$. Indeed, strict concavity and compactness give $g_\mathcal{K}:=\max_{\beta\in \mathcal{K}}\ell_Q(\beta)<\lstar_Q$. The log increments are uniformly bounded, so the strong law for bounded martingale differences gives $\limsup_n \ n^{-1}\ln \Wealth_n\leq g_\mathcal{K}$ almost surely. Proposition~\ref{prop:form} then applies. This includes the betting strategy in \citet{CutkoskyO18} based on the Online Newton Step~\citep{HazanAK07} constrained to $[-\frac12,\frac12]$ for payoffs $X-\mu_0$.

\section{Confidence sequences by inverting the tests}
\label{sec:confidence-sequences}

In this section, we show how to recover guarantees on confidence sequences by inverting tests of candidate means~\citep{Hendriks18,OrabonaJ21,Waudby-SmithR21}. Our construction allows us to convert the uniform $L^1$ rejection-time guarantee into the expected time to reach a prescribed width. All proofs in this section are in Appendix~\ref{app:confidence-sequences}.

Let $\cP$ denote all laws on $[0,1]$, and write
\[
\mu_Q:=\E_Q[X],\qquad \sigma_Q^2:=\Var_Q(X)\qquad(Q\in\cP)~.
\]
For each candidate mean $m\in(0,1)$, choose a probability prior $F_m$ with density $f_m$ on $[-1/(1-m),1/m]$. 
Define the mixture wealth by
\begin{equation}
\label{cs:eq:wealth}
\Wealth_n(m):=\int_{-1/(1-m)}^{1/m}\prod_{i=1}^n(1+\beta(X_i-m))\,\mathrm dF_m(\beta),\quad \text{ and }\quad\Wealth_0(m):=1~.
\end{equation}
For $0<\alpha<1/2$, define the rejection time and its running inversion by
\[
\tau_\alpha(m):=\inf\{n\geq1:\Wealth_n(m)\geq1/\alpha\}~,
\]
\begin{equation}
\label{cs:eq:cs}
C_n^\alpha:=\{m\in(0,1):\tau_\alpha(m)>n\}~.
\end{equation}
Include 0 exactly when $X_1=\cdots=X_n=0$, and include 1 exactly when $X_1=\cdots=X_n=1$.

We measure precision by the time to reach a prescribed diameter, and by its inverse, the effective width of \citet[Definition~5.2]{ShekharR23near}:
\begin{equation}
\label{cs:eq:effective-width}
T_w^\alpha:=\inf\{n\geq1:\operatorname{diam}(C_n^\alpha)\leq w\},\qquad
w_{\mathrm e}(n,Q,\alpha):=\inf\{w>0:\E_Q[T_w^\alpha]\leq n\}~,
\end{equation}
where $\operatorname{diam}(\varnothing):=0$.

Coverage of the obtained sets is immediate, as stated in the following proposition.
\begin{proposition}
\label{cs:prop:validity}
For every $Q\in\cP$,
\[
\Pr_Q\{\mu_Q\in C_n^\alpha\text{ for every }n\geq1\}\geq1-\alpha~.
\]
Coverage also holds for adapted observations with a constant conditional mean. The sets are nested; if $m\mapsto\Wealth_n(m)$ is convex for every $n$ and every observation sequence, they are intervals, possibly empty.
\end{proposition}

\paragraph{Reusing a testing prior for confidence sequences.}
Let's now consider the case that we have a fixed prior and we want to use to generate the priors $F_m$. The admissible interval changes with the candidate mean $m$ in \eqref{cs:eq:wealth}. So, we can rescale the prior accordingly.
Fix a null mean $\mu_0\in(0,1)$ and a prior $F$ satisfying Definition~\ref{def:prior} with $a(F)<\infty$.  To reuse $F$, rescale each bet by the ratio of the corresponding interval endpoints:
\begin{equation}
\label{cs:eq:prior-rescaling}
\beta_m(\beta):=
\begin{cases}
\mu_0\beta/m,&\beta\geq0,\\
(1-\mu_0)\beta/(1-m),&\beta<0.
\end{cases}
\end{equation}
Then, one takes $F_m$ to be the distribution of $\beta_m(\beta)$ when $\beta\sim F$; this gives convex wealth. The next theorem will allow candidate-dependent priors uniformly comparable to this rescaled prior on compact intervals of means.

Write $\beta_Q^\star(m)$ and $\lstar_Q(m)$ for the oracle bet and growth rate with null mean $m$. Use the corresponding directional index from Definition~\ref{def:prior} in the testing benchmark:
\begin{equation}
\label{cs:eq:benchmark}
a_Q(F,m)
:=\begin{cases}
a_+(F),&m<\mu_Q,\\
a_-(F),&m>\mu_Q,
\end{cases}
\qquad
\widetilde B_\gamma(Q,m)
:=\frac{\ln(1/\gamma)+a_Q(F,m)\LL(|\beta_Q^\star(m)|)}{\lstar_Q(m)}~.
\end{equation}
Here $0<\gamma<1$ and $m\ne\mu_Q$. Rescaling and bounded density ratios preserve the two indices, so this is the benchmark of Theorem~\ref{thm:pointwise} for each candidate prior satisfying the comparison below.

\begin{theorem}
\label{cs:thm:width}
Consider i.i.d.\ observations and a prior $F$ satisfying Definition~\ref{def:prior} with $a(F)<\infty$. Suppose $m\mapsto\Wealth_n(m)$ is convex for every observation sequence and $n$. Fix $\kappa\in(0,1/4)$, and assume that constants $0<c_\kappa\leq C_\kappa<\infty$ satisfy
\begin{equation}
\label{cs:eq:prior-comparison}
c_\kappa f(\beta)
\leq f_m(\beta_m(\beta))\beta_m'(\beta)
\leq C_\kappa f(\beta)
\qquad(m\in[\kappa,1-\kappa])
\end{equation}
for almost every $\beta\in(-1/(1-\mu_0),1/\mu_0)$, where $\beta_m$ is defined in \eqref{cs:eq:prior-rescaling}.
Define
\[
\cQ_\kappa
:=\{Q\in\cP:\mu_Q\in[2\kappa,1-2\kappa]\}~.
\]
For each $Q\in\cQ_\kappa$, $0<w<\kappa$, and $0<\alpha<1/2$, let $r^\star\in(0,w)$ be the unique solution of
\begin{equation}
\label{cs:eq:balance}
\widetilde B_\alpha(Q,\mu_Q-r^\star)
=\widetilde B_\alpha(Q,\mu_Q+w-r^\star)~.
\end{equation}
Then
\begin{equation}
\label{cs:eq:uniform-width}
\lim_{\alpha\downarrow0}\ \sup_{\substack{Q\in\cQ_\kappa\\0<w<\kappa}}
\E_Q\left[\left|\frac{T_w^\alpha}{\widetilde B_\alpha(Q,\mu_Q-r^\star)}-1\right|\right]=0~.
\end{equation}
In particular, $\E_Q[T_w^\alpha]/\widetilde B_\alpha(Q,\mu_Q-r^\star)\to1$ uniformly over these laws and widths.
\end{theorem}
For the rescaled prior itself, the assumptions hold with $c_\kappa=C_\kappa=1$.

\begin{remark}
\label{cs:rem:oj}
The exact CO96/R70 mixture of \citet[Remark~2]{OrabonaJ21}, given in \eqref{eq:remark2_test} with $\mu_0=m$, satisfies the theorem with $a(F)=1$. It combines the Beta(1/2,1/2) prior on portfolio allocations for CO96 and the R70 betting prior on $[-1,1]$, with half the initial capital in each. Appendix~\ref{app:cs-oj} verifies the prior comparison and the interval property. Also, any fixed positive allocation to both components gives the same asymptotic limits.
\end{remark}

We now present simplifying corollaries of our result.
\begin{corollary}
\label{cs:cor:local}
Assume the hypotheses of Theorem~\ref{cs:thm:width} for every fixed $\kappa\in(0,1/4)$, with the same reference prior $F$. Fix $Q\in\cP$ with $\sigma_Q^2>0$. Jointly as $\alpha\downarrow0$ and $w\downarrow0$,
\begin{equation}
\label{cs:eq:local-time}
\E_Q[T_w^\alpha]
\leq(1+o(1))\frac{8\sigma_Q^2}{w^2}\left[\ln\frac{1}{\alpha}+a(F)\LL(w)\right]~.
\end{equation}
If $\alpha_n\downarrow0$ and $\ln(1/\alpha_n)=o(n)$, then
\begin{equation}
\label{cs:eq:local-effective}
w_{\mathrm e}(n,Q,\alpha_n)
\leq(1+o(1))\sqrt{\frac{8\sigma_Q^2}{n}\left[\ln\frac{1}{\alpha_n}+a(F)\ln\ln n\right]}~.
\end{equation}
If the two indices in Definition~\ref{def:prior} coincide, both bounds are asymptotic equalities. In particular, this holds when $a(F)=1$.
For each fixed $\alpha\in(0,1/2)$, we also have the upper bounds
\begin{equation}
\label{cs:eq:fixed-alpha}
\limsup_{w\downarrow0} \ \frac{w^2\E_Q[T_w^\alpha]}{8a(F)\sigma_Q^2\LL(w)}\leq1,
\qquad
\limsup_{n\to\infty} \ \frac{w_{\mathrm e}(n,Q,\alpha)}{\sqrt{8a(F)\sigma_Q^2\ln\ln n/n}}\leq1~.
\end{equation}
\end{corollary}

For deterministic observations, the iterated-logarithm penalty vanishes.
\begin{corollary}
\label{cs:cor:deterministic}
Assume the hypotheses of Theorem~\ref{cs:thm:width} for every fixed $\kappa\in(0,1/4)$, with the same reference prior $F$.
Let $Q=\delta_\mu$ with $\mu\in(0,1)$. For $0<w<\min\{\mu,1-\mu\}$, the balanced endpoints are $\mu(1-w)$ and $\mu+(1-\mu)w$, and
\[
\widetilde B_\alpha(\delta_\mu,\mu(1-w))=\frac{\ln(1/\alpha)}{-\ln(1-w)}~.
\]
Consequently, jointly as $\alpha\downarrow0$ and $w\downarrow0$,
\[
\E_{\delta_\mu}[T_w^\alpha]\sim\frac{\ln(1/\alpha)}{w}~.
\]
If $\alpha_n\downarrow0$ and $\ln(1/\alpha_n)=o(n)$, then $w_{\mathrm e}(n,\delta_\mu,\alpha_n)\sim\ln(1/\alpha_n)/n$.
\end{corollary}

\subsection{Comparison with effective-width bounds}
\citet[Theorem~5.6]{ShekharR23near} bound the same effective width for a betting confidence sequence with logarithmic regret, assuming that the observation law is supported strictly inside $[0,1]$. For a fixed law with positive variance, expanding their information radii gives
\[
w_{\mathrm e}(n,Q,\alpha)\leq(4+o(1))\sigma_Q\sqrt{\frac{\ln(1/\alpha)+2\ln n+C_Q}{n}}~,
\]
when the numerator is $o(n)$, where $C_Q$ depends only on the law. Our Corollary~\ref{cs:cor:local} replaces the logarithmic sampling penalty by an iterated logarithm. It also improves the leading constant in this upper bound by a factor $\sqrt2$ when $\alpha\to0$ and $\ln(1/\alpha)\gg\ln n$. Finally, Theorem~\ref{cs:thm:width} allows endpoint support and vanishing variance.

\section{An optimal prior and fast portfolio evaluation}
\label{sec:computation}

We now choose a concrete prior for the bounded-mean test and describe how to evaluate its portfolio wealth. Consider the mixture of \citet[Remark~2]{OrabonaJ21}:
\[
F_w:=wF_{\mathrm{CO96}}+(1-w)F_{\mathrm{R70}},\qquad w:=1/2,
\]
where $F_{\mathrm{CO96}}$ is the prior induced by the Beta(1/2,1/2) distribution on portfolio allocations~\citep{CoverO96}, and $F_{\mathrm{R70}}$ is the fixed betting prior on $[-1,1]$ used by \citet{OrabonaJ21}, adapted from \citet{Robbins70}. Both are defined in Proposition~\ref{prop:oj-exact-mixture} below. We compute the exact mixture wealth and reject at its first threshold crossing:
\begin{equation}
\label{eq:remark2_test}
\Wealth_n^{F_w}=\tfrac12\Wealth_n^{F_{\mathrm{CO96}}}+\tfrac12\Wealth_n^{F_{\mathrm{R70}}},\qquad
\tau_\alpha^{F_w}=\inf\{n\geq1:\Wealth_n^{F_w}\geq1/\alpha\}\,.
\end{equation}
Ville's inequality gives $\sup_{P\in\cP_0}\Pr_P\{\tau_\alpha^{F_w}<\infty\}\leq\alpha$.

We now show that $\eff(\tau^{F_w})=1$.
\begin{proposition}
\label{prop:oj-exact-mixture}
Fix $\mu_0\in(0,1)$ and a weight $w\in(0,1)$, independent of the observations and the error probability. Let $F_{\mathrm{CO96}}$ be the prior induced by $p\sim\mathrm{Beta}(1/2,1/2)$ under $p=\mu_0+\beta\mu_0(1-\mu_0)$. Its density on $(-1/(1-\mu_0),1/\mu_0)$ satisfies
\begin{equation}
\label{eq:oj-kt-density}
f_{\mathrm{CO96}}(\beta):=\frac{\mu_0(1-\mu_0)}{\pi\sqrt{p(1-p)}}
\geq\frac{2\mu_0(1-\mu_0)}{\pi}\,.
\end{equation}
Define the $f_{\mathrm{R70}}$ density by
\begin{equation}
\label{eq:oj-robbins-density}
f_{\mathrm{R70}}(\beta):=\frac{\ln\ln c}{2|\beta|\ln(c/|\beta|)[\ln\ln(c/|\beta|)]^2}
\quad (0<|\beta|\leq1),\qquad c:=6.6e,
\end{equation}
and set $f_{\mathrm{R70}}(\beta):=0$ for $|\beta|>1$; its value at zero is immaterial. Let $F_{\mathrm{R70}}$ denote its distribution and define
\[
F_w:=wF_{\mathrm{CO96}}+(1-w)F_{\mathrm{R70}},\qquad
f_w:=wf_{\mathrm{CO96}}+(1-w)f_{\mathrm{R70}}\,.
\]
Then, the test based on the wealth $\Wealth_n^{F_w}=w\Wealth_n^{F_{\mathrm{CO96}}}+(1-w)\Wealth_n^{F_{\mathrm{R70}}}$ has $\eff(\tau^{F_w})=1$ over the full bounded alternative class $\cQ$, and satisfies Theorem~\ref{thm:pointwise} with $a_Q(F_w)=1$.
\end{proposition}
\begin{proof}
The substitution $u=\ln(c/|\beta|)$ shows that $f_{\mathrm{R70}}$ integrates to 1. The density $f_w$ has a positive lower bound on every compact subset of the interior away from zero, by \eqref{eq:oj-kt-density}. Since $f_{\mathrm{CO96}}$ is bounded near zero and $f_{\mathrm{R70}}(\beta)\to\infty$, we have $f_w(\beta)\sim(1-w)f_{\mathrm{R70}}(\beta)$ as $\beta\to0^\pm$. Equation~\eqref{eq:oj-robbins-density} therefore gives $a_+(F_w)=a_-(F_w)=1$. Theorem~\ref{thm:pointwise} and Corollary~\ref{thm:main} apply, including the jumps at $\beta=\pm1$.
\end{proof}

Using \eqref{eq:pointwise-upper} we can be even more precise: for every $\delta\in(0,1/2)$ and $\epsilon>0$, we have
\[
\E_Q[\tau_\alpha^{F_w}]\leq
\frac{\ln(1/\alpha)+(1+\epsilon)\LL(|\beta_Q^\star|)+C_{\mu_0,\delta,\epsilon}}{(1-2\delta)\lstar_Q}
\qquad(Q\in\cQ),
\]
where the constant is independent of $Q$ and $\alpha$. Theorem~\ref{thm:pointwise} also gives $\E_Q[\tau_\alpha^{F_w}]/B_\alpha(Q)\to1$ uniformly over $Q$ as $\alpha\to0$.

\subsection{Comparison with regret-based tests}
\label{sec:regret-comparison}

The fixed positive CO96 weight also retains logarithmic portfolio regret: $\Wealth_n^{F_w}\geq w\Wealth_n^{F_{\mathrm{CO96}}}$ adds at most $\ln(1/w)$ to the CO96 regret bound. Thus, the mixture inherits the known growth guarantees of \citet[Theorems~2.3 and~5.2]{Waudby-SmithSJ25} while attaining efficiency 1.

The computable PRECiSE bounds in \citet[Algorithms~1 and~3]{OrabonaJ21} use regret-based wealth lower bounds, whose inversion can give different confidence sequences. Using a regret upper bound to lower bound wealth requires a separate expected-time analysis. Theorem~\ref{thm:oj-efficient-regret} in Appendix~\ref{app:oj-mixture} derives such a bound for the exact prior mixture in \citet[Remark~2]{OrabonaJ21}, while preserving efficiency 1 and the uniform expected-time formula. In contrast, Proposition~\ref{prop:ar-regret-counterexample} in Appendix~\ref{app:ar-counterexample} gives a sequence along which the test based on the explicit regret bound by \citet{AgrawalR26} takes twice the benchmark and twice the expected time of its exact mixture asymptotically. Thus, this proxy has efficiency at least 2, while the exact mixture has efficiency 1.

\subsection{Numerical evaluation by a regional grid and FFT}

Evaluating the wealth in \eqref{eq:remark2_test} can be computationally demanding. \citet[Remark~2.2]{Waudby-SmithSJ25} suggested using an exact implementation of $F$-weighted portfolios for small sample sizes, and the regret-based approach of \citet{OrabonaJ21} for large sample sizes. Here, we show how to evaluate the continuous mixture efficiently.

We combine a fixed spatial partition with FFT integration inside each region. Every region retains the continuous conditional prior, represented by its Chebyshev moments through the required polynomial degree. Given these moments exactly, the algorithm evaluates the original mixture exactly in exact arithmetic. We use a one-dimensional specialization of the grid operations and \texttt{Prefix} recursion of \citet[Appendices~E.3--E.4]{LiuPO26}. Their algorithm handles a restarting three-stock portfolio; the fixed-start two-stock specialization and regional combination below are our derivation. We provide this implementation as OPT-FFT (Optimal Portfolio Testing with Fast Fourier Transforms) at \url{https://github.com/AdaMLLab/OPT-FFT}.

To describe the computation, use the portfolio coordinate
\[
p:=\mu_0+\beta\mu_0(1-\mu_0)\in[0,1],
\]
and let $\overline F$ be its distribution under $F$. By \eqref{eq:two-stock-returns}, a single observation contributes the linear polynomial
\[
r_i(p):=p\frac{X_i}{\mu_0}+(1-p)\frac{1-X_i}{1-\mu_0}~.
\]
The fixed-allocation wealth after $n$ observations is the degree-$n$ polynomial $\prod_{i=1}^n r_i(p)$. Its integral under $\overline F$ is $\Wealth_n^F$, and the online allocation to the first stock is
\[
p_n:=\frac{\int_0^1 p\prod_{i=1}^{n-1}r_i(p)\,\mathrm{d}\overline F(p)}
{\int_0^1\prod_{i=1}^{n-1}r_i(p)\,\mathrm{d}\overline F(p)}~.
\]
If the denominator is zero, wealth remains zero and we may choose any allocation.

Before observing data, partition the allocation interval into $K$ regions. On each region, map the coordinate affinely to $z\in[-1,1]$ and use the conditional prior's moments to construct an integration rule on Chebyshev points. The rule covers every polynomial degree needed up to the chosen horizon $T$. The \texttt{Prefix} recursion evaluates the integrals of each prefix product and that product multiplied by $z$. Splitting the observations into two halves, we absorb the first half's product into the integration functional for the second half. Polynomial multiplication is pointwise on a sufficiently large grid; interpolation and its transpose use cosine transforms computed by FFT. The ratio of the two integrals gives the posterior mean of $z$, which we map back to the allocation interval.

Finally, combine the regional allocations using their preceding wealths and the original prior interval masses. In exact arithmetic, each regional wealth is the continuous conditional mixture, so their weighted sum is exactly $\Wealth_n^F$. Appendix~\ref{app:prior-moments} gives the integration rule, pseudocode in Algorithms~\ref{alg:regional-fft}--\ref{alg:regional-prefix}, and a proof.

\paragraph{Why the grid improves numerical stability.}
Note that the grid is not necessary for the exactness of the algorithm in exact arithmetic.
However, after a long run favoring one direction, fixed-allocation wealths can differ by many orders of magnitude. Then, cosine interpolation and its transpose combine signed terms: roundoff that is small relative to the largest entries can overwhelm small contributions that become important after a reversal~\citep[Appendix~E.6]{LiuPO26}. This affects both CO96 and R70 prior. A single common rescaling prevents overflow but does not restore lost relative accuracy.

For this reason we use the decomposition in a grid that permits independent rescaling of products and functionals in different parts of the bet interval. Note that we do not claim any bound of the error due to the floating point truncations, but in our experiments this method substantially improved the accuracy.
We store each region's accumulated wealth logarithmically and combine regions by positive weights, so a region can regain posterior weight after a reversal. Within a narrow region, we compute the posterior mean of the local coordinate $z$ before mapping it back. This avoids computing an allocation from nearly proportional numerator and denominator functionals; an error in the local mean is also multiplied by the region's half-width. We place more boundaries near the zero bet and the support endpoints, where small betting scales and concentrated likelihoods need resolution.

\paragraph{Calculations of the moments.}
The R70 singularity at zero requires a dedicated moment calculation. We isolate it in an interval centered at zero, use its analytic prior mass, and pair the contributions at $\beta$ and $-\beta$ in the change of variables described in Appendix~\ref{app:prior-moments}. The continuous prior retains all its mass near zero. The FFT integration weights can still be signed, and products can still concentrate inside a region, so these measures do not give a uniform roundoff guarantee.

\paragraph{Cost and numerical choices.}
Given the conditional moments, $K$ regions cost $O(KT\ln^2(T+1))$ arithmetic operations and $O(KT)$ memory per candidate mean. We choose $K=O(\ln^2(T+1))$, giving $O(T\ln^4(T+1))=\widetilde O(T)$ operations and $O(T\ln^2(T+1))$ memory. Numerical moment preprocessing also has cost $\widetilde O(T)$ for fixed accuracy parameters.

In floating point, we project each computed regional allocation onto its interval, combine regions using their preceding wealths, and accumulate the actual returns of the resulting predictable allocations. This construction retains betting validity under the null, whereas substituting arbitrary approximate mixture values need not. The exact-arithmetic algorithm computes the continuous mixture covered by our expected-time results; their extension to the rounded strategy requires a separate analysis.

\section{Conclusion}
\label{sec:conclusion}

We have presented a novel characterization of how a prior's density near zero determines expected rejection time in bounded i.i.d.\ testing. Under Definition~\ref{def:prior}, the expected-time formula is uniform over all alternatives as $\alpha\to0$ for finite indices, and the efficiency identity $\eff(\tau^F)=a(F)$ also covers infinite indices. The adaptation cost depends on the oracle bet through $\LL(|\beta_Q^\star|)$, including when variance vanishes. Thus, priors with index 1 are optimal under this criterion, certifying existing Robbins-type constructions. These priors also attain the joint minimax constant over ranges of small oracle bets among all valid sequential tests. We also present a rule for combining priors showing that different components can supply the required mass on the two sides of zero: two priors with infinite efficiency can together attain efficiency 1.

The uniform $L^1$ limit also gives a general route from rejection times to confidence-sequence widths through test inversion and balanced endpoint benchmarks. This yields sharp expected times to reach a prescribed width, including at laws with vanishing variance.

Our characterization also certifies the expected-time optimality of the prior proposed in \citet[Remark~2]{OrabonaJ21}. For this prior, we combine a spatial grid with FFT recursion inside each region, retaining betting validity. Given exact prior moments, this computes the continuous mixture exactly in exact arithmetic, using $\widetilde O(T)$ operations per candidate mean.

It remains open how to distinguish between different priors with index 1. For example, higher-order expected-time bounds could guide the choice of their parameters and hedging weights.

\section*{Acknowledgments}
The author acknowledges the use of GPT 6 and Claude 5.5 for assistance with proofs and code. The author has checked all details for correctness and takes full responsibility for the content of this paper.

\bibliographystyle{abbrvnat}
\begingroup
\RaggedRight
\setlength{\bibsep}{6pt plus 1pt minus 1pt}
\bibliography{learning}

\begin{thebibliography}{53}
\providecommand{\natexlab}[1]{#1}
\providecommand{\url}[1]{\texttt{#1}}
\expandafter\ifx\csname urlstyle\endcsname\relax
  \providecommand{\doi}[1]{doi: #1}\else
  \providecommand{\doi}{doi: \begingroup \urlstyle{rm}\Url}\fi

\bibitem[Agrawal and Ramdas(2025)]{AgrawalR25LIL}
S.~Agrawal and A.~Ramdas.
\newblock Eventually {LIL} regret: Almost sure $\ln\ln t$ regret for a
  sub-{G}aussian mixture on unbounded data.
\newblock \emph{arXiv preprint arXiv:2512.12325}, 2025.
\newblock URL \url{https://arxiv.org/pdf/2512.12325v1}.

\bibitem[Agrawal and Ramdas(2026)]{AgrawalR26}
S.~Agrawal and A.~Ramdas.
\newblock Cover meets {Robbins} while betting on bounded data: $\ln n$ regret
  and almost sure $\ln\ln n$ regret.
\newblock \emph{arXiv preprint arXiv:2604.20172v3}, 2026.
\newblock URL \url{https://arxiv.org/pdf/2604.20172v3}.

\bibitem[Agrawal et~al.(2026)Agrawal, Ram, and Ramdas]{AgrawalRR26}
S.~Agrawal, A.~Ram, and A.~Ramdas.
\newblock On stopping times of power-one sequential tests: Tight lower and
  upper bounds.
\newblock \emph{arXiv preprint arXiv:2504.19952v2}, 2026.
\newblock URL \url{https://arxiv.org/pdf/2504.19952v2}.

\bibitem[Balsubramani(2015)]{Balsubramani15}
A.~Balsubramani.
\newblock Sharp finite-time iterated-logarithm martingale concentration.
\newblock \emph{arXiv preprint arXiv:1405.2639v4}, 2015.
\newblock URL \url{https://arxiv.org/pdf/1405.2639v4}.

\bibitem[Breiman(1961)]{Breiman61}
L.~Breiman.
\newblock Optimal gambling systems for favorable games.
\newblock In \emph{Proceedings of the Fourth Berkeley Symposium on Mathematical
  Statistics and Probability}, pages 65--78, 1961.
\newblock URL
  \url{https://digicoll.lib.berkeley.edu/record/112884/files/math_s4_v1_article-05.pdf}.

\bibitem[Chen et~al.(2017)Chen, Li, and Qiao]{ChenLQ17}
L.~Chen, J.~Li, and M.~Qiao.
\newblock Towards instance optimal bounds for best arm identification.
\newblock In \emph{Proceedings of the 2017 Conference on Learning Theory},
  volume~65 of \emph{Proceedings of Machine Learning Research}, pages 535--592.
  PMLR, 2017.
\newblock URL \url{https://proceedings.mlr.press/v65/chen17b/chen17b.pdf}.

\bibitem[Clarke and Barron(1990)]{ClarkeB90}
B.~S. Clarke and A.~R. Barron.
\newblock Information-theoretic asymptotics of {Bayes} methods.
\newblock \emph{IEEE Transactions on Information Theory}, 36\penalty0
  (3):\penalty0 453--471, 1990.
\newblock URL
  \url{https://bertrandsclarke.weebly.com/uploads/2/7/3/6/2736392/clarkebarron1990.pdf}.

\bibitem[Clerico(2025)]{Clerico25}
E.~Clerico.
\newblock On the optimality of coin-betting for mean estimation.
\newblock \emph{International Journal of Approximate Reasoning}, page 109550,
  2025.
\newblock URL \url{https://arxiv.org/pdf/2412.02640v4}.

\bibitem[Cover(1991)]{Cover91}
T.~M. Cover.
\newblock Universal portfolios.
\newblock \emph{Mathematical Finance}, pages 1--29, 1991.
\newblock URL
  \url{https://isl.stanford.edu/~cover/papers/universal_portfolios.pdf}.

\bibitem[Cover and Ordentlich(1996)]{CoverO96}
T.~M. Cover and E.~Ordentlich.
\newblock Universal portfolios with side information.
\newblock \emph{IEEE Transactions on Information Theory}, 42\penalty0
  (2):\penalty0 348--363, 1996.
\newblock URL
  \url{https://isl.stanford.edu/~cover/papers/cover_ordentlich_96.pdf}.

\bibitem[Cover and Thomas(2006)]{CoverT06}
T.~M. Cover and J.~A. Thomas.
\newblock \emph{Elements of information theory}.
\newblock Wiley-Interscience, 2006.

\bibitem[Cutkosky and Orabona(2018)]{CutkoskyO18}
A.~Cutkosky and F.~Orabona.
\newblock Black-box reductions for parameter-free online learning in {Banach}
  spaces.
\newblock In \emph{Proc. of the Conference on Learning Theory (COLT)}, 2018.
\newblock URL \url{https://arxiv.org/pdf/1802.06293}.

\bibitem[Darling and Robbins(1967)]{DarlingR67}
D.~A. Darling and H.~Robbins.
\newblock Confidence sequences for mean, variance, and median.
\newblock \emph{Proceedings of the National Academy of Sciences}, 58\penalty0
  (1):\penalty0 66--68, 1967.
\newblock URL \url{https://www.pnas.org/doi/pdf/10.1073/pnas.58.1.66}.

\bibitem[Dawid and Vovk(1999)]{DawidV99}
A.~P. Dawid and V.~G. Vovk.
\newblock Prequential probability: Principles and properties.
\newblock \emph{Bernoulli}, 5\penalty0 (1):\penalty0 125--162, 1999.
\newblock \doi{10.2307/3318616}.
\newblock URL
  \url{https://projecteuclid.org/download/pdf_1/euclid.bj/1173707098}.

\bibitem[Deep and Agrawal(2026)]{DeepA26}
V.~Deep and S.~Agrawal.
\newblock Beyond first-order asymptotics in sequential mean testing.
\newblock In \emph{Proceedings of the 43rd International Conference on Machine
  Learning}, volume 306 of \emph{Proceedings of Machine Learning Research},
  pages 23527--23561. PMLR, 2026.
\newblock URL
  \url{https://raw.githubusercontent.com/mlresearch/v306/main/assets/deep26a/deep26a.pdf}.

\bibitem[Durand and Wintenberger(2025)]{DurandW25}
A.~Durand and O.~Wintenberger.
\newblock Power comparison of sequential testing by betting procedures.
\newblock \emph{arXiv preprint arXiv:2504.00593v2}, 2025.
\newblock URL \url{https://arxiv.org/pdf/2504.00593v2}.

\bibitem[Erd\H{o}s(1942)]{Erdos42}
P.~Erd\H{o}s.
\newblock On the law of the iterated logarithm.
\newblock \emph{Annals of Mathematics}, 43\penalty0 (3):\penalty0 419--436,
  1942.
\newblock URL \url{https://www.renyi.hu/~p_erdos/1942-04.pdf}.

\bibitem[Farrell(1964)]{Farrell64}
R.~H. Farrell.
\newblock Asymptotic behavior of expected sample size in certain one sided
  tests.
\newblock \emph{The Annals of Mathematical Statistics}, 35\penalty0
  (1):\penalty0 36--72, 1964.
\newblock URL
  \url{https://projecteuclid.org/journals/annals-of-mathematical-statistics/volume-35/issue-1/Asymptotic-Behavior-of-Expected-Sample-Size-in-Certain-One-Sided/10.1214/aoms/1177703731.pdf}.

\bibitem[Feller(1943)]{Feller43}
W.~Feller.
\newblock The general form of the so-called law of the iterated logarithm.
\newblock \emph{Transactions of the American Mathematical Society},
  54:\penalty0 373--402, 1943.
\newblock URL
  \url{https://www.ams.org/tran/1943-054-03/S0002-9947-1943-0009263-7/S0002-9947-1943-0009263-7.pdf}.

\bibitem[Fellouris and Tartakovsky(2012)]{FellourisT12}
G.~Fellouris and A.~G. Tartakovsky.
\newblock Nearly minimax one-sided mixture-based sequential tests.
\newblock \emph{Sequential Analysis}, 31\penalty0 (3):\penalty0 297--325, 2012.
\newblock \doi{10.1080/07474946.2012.694346}.
\newblock URL \url{https://arxiv.org/pdf/1110.0902v2}.

\bibitem[Garivier and Kaufmann(2016)]{GarivierK16}
A.~Garivier and E.~Kaufmann.
\newblock Optimal best arm identification with fixed confidence.
\newblock In \emph{Proc. of the Conference on Learning Theory (COLT)}, pages
  998--1027, 2016.
\newblock URL \url{https://proceedings.mlr.press/v49/garivier16a.pdf}.

\bibitem[Gr{\"u}nwald et~al.(2024)Gr{\"u}nwald, de~Heide, and
  Koolen]{GrunwaldHK24}
P.~Gr{\"u}nwald, R.~de~Heide, and W.~M. Koolen.
\newblock Safe testing.
\newblock \emph{Journal of the Royal Statistical Society Series B}, 86\penalty0
  (5):\penalty0 1091--1128, 2024.
\newblock URL \url{https://arxiv.org/pdf/1906.07801}.

\bibitem[Hazan et~al.(2007)Hazan, Agarwal, and Kale]{HazanAK07}
E.~Hazan, A.~Agarwal, and S.~Kale.
\newblock Logarithmic regret algorithms for online convex optimization.
\newblock \emph{Machine Learning}, 69\penalty0 (2-3):\penalty0 169--192, 2007.
\newblock URL
  \url{https://link.springer.com/content/pdf/10.1007/s10994-007-5016-8.pdf}.

\bibitem[Hendriks(2018)]{Hendriks18}
H.~Hendriks.
\newblock Test martingales for bounded random variables.
\newblock \emph{arXiv preprint arXiv:1801.09418}, 2018.
\newblock URL \url{https://arxiv.org/abs/1801.09418}.

\bibitem[Hoeffding(1963)]{Hoeffding63}
W.~Hoeffding.
\newblock Probability inequalities for sums of bounded random variables.
\newblock \emph{Journal of the American Statistical Association}, 58\penalty0
  (301):\penalty0 13--30, 1963.
\newblock URL
  \url{https://www.cs.rpi.edu/academics/courses/spring06/random/hoefding.pdf}.

\bibitem[Honda and Takemura(2010)]{HondaT10}
J.~Honda and A.~Takemura.
\newblock An asymptotically optimal bandit algorithm for bounded support
  models.
\newblock In \emph{Proc. of the Conference on Learning Theory (COLT)}, pages
  67--79, 2010.
\newblock URL \url{https://www.learningtheory.org/colt2010/papers/29honda.pdf}.

\bibitem[Howard et~al.(2021)Howard, Ramdas, McAuliffe, and Sekhon]{HowardRMS21}
S.~R. Howard, A.~Ramdas, J.~McAuliffe, and J.~Sekhon.
\newblock Time-uniform, nonparametric, nonasymptotic confidence sequences.
\newblock \emph{The Annals of Statistics}, 49\penalty0 (2):\penalty0
  1055--1080, 2021.
\newblock URL \url{https://arxiv.org/pdf/1810.08240}.

\bibitem[Jamieson et~al.(2014)Jamieson, Malloy, Nowak, and
  Bubeck]{JamiesonMNB14}
K.~Jamieson, M.~Malloy, R.~Nowak, and S.~Bubeck.
\newblock lil' {UCB}: An optimal exploration algorithm for multi-armed bandits.
\newblock In \emph{Proc. of the Conference on Learning Theory (COLT)}, pages
  423--439, 2014.
\newblock URL \url{https://proceedings.mlr.press/v35/jamieson14.pdf}.

\bibitem[Jun and Orabona(2019)]{JunO19}
K.-S. Jun and F.~Orabona.
\newblock Parameter-free online convex optimization with sub-exponential noise.
\newblock In \emph{Proc. of the Conference on Learning Theory (COLT)}, 2019.
\newblock URL \url{https://proceedings.mlr.press/v99/jun19a/jun19a.pdf}.

\bibitem[Kaufmann and Koolen(2021)]{KaufmannK21}
E.~Kaufmann and W.~M. Koolen.
\newblock Mixture martingales revisited with applications to sequential tests
  and confidence intervals.
\newblock \emph{Journal of Machine Learning Research}, 22\penalty0
  (246):\penalty0 1--44, 2021.
\newblock URL \url{https://ir.cwi.nl/pub/31315/31315.pdf}.

\bibitem[Kelly(1956)]{Kelly56}
J.~L. Kelly, jr.
\newblock A new interpretation of information rate.
\newblock \emph{IRE Transactions on Information Theory}, 2\penalty0
  (3):\penalty0 185--189, 1956.
\newblock URL
  \url{https://ieeexplore.ieee.org/stamp/stamp.jsp?arnumber=1056803}.

\bibitem[Krichevsky and Trofimov(1981)]{KrichevskyT81}
R.~Krichevsky and V.~Trofimov.
\newblock The performance of universal encoding.
\newblock \emph{IEEE Trans. on Information Theory}, 27\penalty0 (2):\penalty0
  199--207, 1981.
\newblock URL
  \url{https://ieeexplore.ieee.org/stamp/stamp.jsp?arnumber=1056331}.

\bibitem[Larsson et~al.(2025)Larsson, Ramdas, and Ruf]{LarssonRR25}
M.~Larsson, A.~Ramdas, and J.~Ruf.
\newblock The numeraire e-variable and reverse information projection.
\newblock \emph{The Annals of Statistics}, 53\penalty0 (3):\penalty0
  1015--1043, 2025.
\newblock URL \url{https://arxiv.org/pdf/2402.18810}.

\bibitem[Liu et~al.(2026)Liu, Prinster, and Orabona]{LiuPO26}
T.~Liu, D.~Prinster, and F.~Orabona.
\newblock Monitoring for shifts at unknown times: Parameter-free conditional
  conformal martingales.
\newblock \emph{arXiv preprint arXiv:2610.17343}, 2026.
\newblock URL \url{https://arxiv.org/pdf/2610.17343v1}.

\bibitem[Orabona(2018)]{Orabona26b}
F.~Orabona.
\newblock From betting to empirical {Bernstein} {LIL}.
\newblock \emph{arXiv preprint arXiv:2605.22124}, 2018.
\newblock URL \url{https://arxiv.org/pdf/2605.22124v1}.
\newblock Technical report dated January~22, 2018; posted May~21, 2026 as a
  verbatim copy of the 2017--2018 report.

\bibitem[Orabona and Jun(2024)]{OrabonaJ21}
F.~Orabona and K.-S. Jun.
\newblock Tight concentrations and confidence sequences from the regret of
  universal portfolio.
\newblock \emph{{IEEE} Trans. Inf. Theory}, 70\penalty0 (1):\penalty0 436--455,
  2024.
\newblock URL \url{https://arxiv.org/pdf/2110.14099}.

\bibitem[Orabona and P\'al(2016)]{OrabonaP16}
F.~Orabona and D.~P\'al.
\newblock Coin betting and parameter-free online learning.
\newblock In D.~D. Lee, M.~Sugiyama, U.~V. Luxburg, I.~Guyon, and R.~Garnett,
  editors, \emph{Advances in Neural Information Processing Systems 29}, pages
  577--585. Curran Associates, Inc., 2016.
\newblock URL \url{https://arxiv.org/pdf/1602.04128.pdf}.

\bibitem[Ramdas and Wang(2025)]{RamdasW25}
A.~Ramdas and R.~Wang.
\newblock Hypothesis testing with e-values.
\newblock \emph{Foundations and Trends in Statistics}, 1\penalty0
  (1-2):\penalty0 1--390, 2025.
\newblock URL \url{https://arxiv.org/pdf/2410.23614}.

\bibitem[Rissanen(1983)]{Rissanen83}
J.~Rissanen.
\newblock A universal prior for integers and estimation by minimum description
  length.
\newblock \emph{The Annals of Statistics}, 11\penalty0 (2):\penalty0 416--431,
  1983.
\newblock URL
  \url{https://www.ccs.neu.edu/home/jaa/CS188.01W/Papers/Rissanen83.pdf}.

\bibitem[Robbins(1970)]{Robbins70}
H.~Robbins.
\newblock Statistical methods related to the law of the iterated logarithm.
\newblock \emph{The Annals of Mathematical Statistics}, 41\penalty0
  (5):\penalty0 1397--1409, 1970.
\newblock URL
  \url{https://projecteuclid.org/journals/annals-of-mathematical-statistics/volume-41/issue-5/Statistical-Methods-Related-to-the-Law-of-the-Iterated-Logarithm/10.1214/aoms/1177696786.pdf}.

\bibitem[Robbins and Siegmund(1974)]{RobbinsS74}
H.~Robbins and D.~Siegmund.
\newblock The expected sample size of some tests of power one.
\newblock \emph{The Annals of Statistics}, 2\penalty0 (3):\penalty0 415--436,
  1974.
\newblock URL
  \url{https://projecteuclid.org/journals/annals-of-statistics/volume-2/issue-3/The-Expected-Sample-Size-of-Some-Tests-of-Power-One/10.1214/aos/1176342704.pdf}.

\bibitem[Ryu and Bhatt(2024)]{RyuB22}
J.~J. Ryu and A.~Bhatt.
\newblock On confidence sequences for bounded random processes via universal
  gambling strategies.
\newblock \emph{IEEE Transactions on Information Theory}, 2024.
\newblock URL \url{https://arxiv.org/pdf/2207.12382}.

\bibitem[Sato et~al.(2018)Sato, Miyabe, and Takemura]{SatoMT18}
R.~Sato, K.~Miyabe, and A.~Takemura.
\newblock Relation between the rate of convergence of strong law of large
  numbers and the rate of concentration of {Bayesian} prior in game-theoretic
  probability.
\newblock \emph{Stochastic Processes and their Applications}, 128\penalty0
  (5):\penalty0 1466--1484, 2018.
\newblock \doi{10.1016/j.spa.2017.07.014}.
\newblock URL \url{https://arxiv.org/pdf/1604.07911v1}.

\bibitem[Shafer(2021)]{Shafer21}
G.~Shafer.
\newblock Testing by betting: A strategy for statistical and scientific
  communication.
\newblock \emph{Journal of the Royal Statistical Society Series A},
  184\penalty0 (2):\penalty0 407--431, 2021.
\newblock URL
  \url{https://glennshafer.com/assets/downloads/articles/article104_jrss_shafer-testingbybetting-with-discussion-response.pdf}.

\bibitem[Shafer and Vovk(2001)]{ShaferV05}
G.~Shafer and V.~Vovk.
\newblock \emph{Probability and finance: it's only a game!}
\newblock John Wiley \& Sons, 2001.

\bibitem[Shekhar and Ramdas(2023)]{ShekharR23near}
S.~Shekhar and A.~Ramdas.
\newblock On the near-optimality of betting confidence sets for bounded means.
\newblock \emph{arXiv preprint arXiv:2310.01547v2}, 2023.
\newblock URL \url{https://arxiv.org/pdf/2310.01547v2}.

\bibitem[Shin et~al.(2021)Shin, Ramdas, and Rinaldo]{ShinRR21}
J.~Shin, A.~Ramdas, and A.~Rinaldo.
\newblock Nonparametric iterated-logarithm extensions of the sequential
  generalized likelihood ratio test.
\newblock \emph{IEEE Journal on Selected Areas in Information Theory},
  2\penalty0 (2):\penalty0 691--704, 2021.
\newblock URL \url{https://arxiv.org/pdf/2010.08082v4}.

\bibitem[Ville(1939)]{Ville39}
J.~Ville.
\newblock \emph{Étude critique de la notion de collectif}.
\newblock Gauthier-Villars, Paris, 1939.
\newblock URL \url{https://www.numdam.org/item/THESE_1939__218__1_0.pdf}.

\bibitem[Vovk and Wang(2021)]{VovkW21}
V.~Vovk and R.~Wang.
\newblock E-values: Calibration, combination and applications.
\newblock \emph{The Annals of Statistics}, 49\penalty0 (3):\penalty0
  1736--1754, 2021.

\bibitem[Vovk(1988)]{Vovk88}
V.~G. Vovk.
\newblock {Kolmogorov--Stout} law of the iterated logarithm.
\newblock \emph{Mathematical Notes}, 44\penalty0 (1):\penalty0 502--507, 1988.
\newblock \doi{10.1007/BF01158114}.
\newblock URL
  \url{https://link.springer.com/content/pdf/10.1007/BF01158114.pdf}.

\bibitem[Wang et~al.(2025)Wang, Wang, and Ziegel]{WangWZ25}
Q.~Wang, R.~Wang, and J.~Ziegel.
\newblock E-backtesting.
\newblock \emph{Management Science}, 2025.
\newblock URL \url{https://arxiv.org/pdf/2209.00991}.

\bibitem[Waudby-Smith and Ramdas(2024)]{Waudby-SmithR21}
I.~Waudby-Smith and A.~Ramdas.
\newblock Estimating means of bounded random variables by betting.
\newblock \emph{Journal of the Royal Statistical Society Series B: Statistical
  Methodology}, 86\penalty0 (1):\penalty0 1--27, 2024.
\newblock URL \url{https://arxiv.org/pdf/2010.09686}.

\bibitem[Waudby-Smith et~al.(2026)Waudby-Smith, Sandoval, and
  Jordan]{Waudby-SmithSJ25}
I.~Waudby-Smith, R.~J. Sandoval, and M.~I. Jordan.
\newblock Instance-log-optimality of portfolio-based e-processes and their
  sequential hypothesis tests.
\newblock \emph{arXiv preprint arXiv:2504.02818v2}, 2026.
\newblock URL \url{https://arxiv.org/pdf/2504.02818v2}.

\end{thebibliography}
\endgroup

\appendix

\section{Proofs of Section~\ref{sec:main}}
\label{sec:proofs}

Throughout this section, constants may depend on $\mu_0$ and $F$, and on explicitly fixed approximation parameters, but never on $Q$ or $\alpha$. Recall that $0\leq1+\beta (X-\mu_0)\leq \max\{1/\mu_0,1/(1-\mu_0)\}$ for every admissible bet. The oracle bet $\beta^\star_Q$ maximizes $\ell_Q(\beta)=\E_Q[\ln(1+\beta (X-\mu_0))]$, with value $\lstar_Q>0$.

First of all, we introduce the lemmas that we will use in our proofs.

The first lemma adapts the KL-projection construction in \citet[Lemma~A.4]{Waudby-SmithSJ25} to the bounded-mean null and bounds all positive-order moments of the oracle log multiplier uniformly over alternatives. The case $s=1$ supplies the integrability needed below.
\begin{lemma}
\label{lemma:proj}
For every $Q\in\cQ$, there exists $P^\star\in\cP_0$ such that $Q\ll P^\star$ and
\begin{equation}
\label{eq:oracle-projection}
\frac{\mathrm{d}Q}{\mathrm{d}P^\star}(X)=1+\beta^\star_Q(X-\mu_0),\qquad
\KL(Q\|P^\star)=\lstar_Q=\inf_{P\in\cP_0}\KL(Q\|P)\,.
\end{equation}
Moreover, for every $s>0$,
\begin{equation}
\label{eq:oracle-log-moments}
\sup_{Q\in\cQ} \ \E_Q\left[|\ln(1+\beta^\star_Q(X-\mu_0))|^s\right]
\leq (\ln \max\{1/\mu_0,1/(1-\mu_0)\})^s+(s/e)^s
<\infty\,.
\end{equation}
\end{lemma}
\begin{proof}
Let $r:=\E_Q[(1+\beta^\star_Q(X-\mu_0))^{-1}]$. For $0<t<1$, the function $t\mapsto\ell_Q(t\beta^\star_Q)$ is concave and nondecreasing, since its maximum is at $t=1$. On any closed interval of $t$ values inside $(0,1)$, the denominator $1+t\beta^\star_Q(X-\mu_0)\geq1-t$ stays bounded away from zero and the numerator $\beta^\star_Q(X-\mu_0)$ is bounded. We can therefore differentiate under the expectation to obtain $\E_Q[\beta^\star_Q(X-\mu_0)/(1+t\beta^\star_Q(X-\mu_0))]\geq0$. So, we obtain $\E_Q[(1+t\beta^\star_Q(X-\mu_0))^{-1}]=1-t\E_Q[\beta^\star_Q(X-\mu_0)/(1+t\beta^\star_Q(X-\mu_0))]\leq1$. Fatou's lemma gives $r\leq1$ and $\Pr_{Q}\{\beta^\star_Q(X-\mu_0)=-1\}=0$.

For $u\geq0$, maximizing $u^se^{-u}$ at $u=s$ gives $u^s\leq(s/e)^se^u$. Taking $u:=-\ln(1+\beta^\star_Q(X-\mu_0))$ on $\{1+\beta^\star_Q(X-\mu_0)<1\}$ and using $r\leq1$ bounds the negative-part $s$th moment by $(s/e)^s$; the positive-part moment is at most $(\ln \max\{1/\mu_0,1/(1-\mu_0)\})^s$. This proves \eqref{eq:oracle-log-moments}, in particular integrability.

We now construct the null law by dividing $Q$ by the oracle multiplier; if this leaves total mass below one, we add the missing mass where that multiplier is zero.

If $\beta^\star_Q\in(-1/(1-\mu_0),1/\mu_0)$, differentiation at the maximizer gives $\E_Q[(X-\mu_0)/(1+\beta^\star_Q(X-\mu_0))]=0$ and $r=1$. Define $P^\star$ by $\mathrm{d}P^\star:=(1+\beta^\star_Q(X-\mu_0))^{-1}\mathrm{d}Q$ in this case. Otherwise, let $x_b:=\mu_0-1/\beta^\star_Q\in\{0,1\}$ and define
\begin{equation}
\label{eq:projection-law}
\mathrm{d}P^\star(x):=\frac{\mathrm{d}Q(x)}{1+\beta^\star_Q(x-\mu_0)}+(1-r)\,\mathrm{d}\delta_{x_b}(x)\,.
\end{equation}
Since $Q$ is the law of $X$,
\[
\int_{[0,1]}\frac{\mathrm{d}Q(x)}{1+\beta^\star_Q(x-\mu_0)}=\E_Q[(1+\beta^\star_Q(X-\mu_0))^{-1}]=r\,.
\]
Here, $\delta_{x_b}$ is the unit point mass at $x_b$. Since $r\leq1$, both terms in \eqref{eq:projection-law} are nonnegative, and the preceding identity shows that their masses sum to $r+(1-r)=1$. In the interior case, $r=1$ already. Thus, $P^\star$ is a probability measure. The identity $\beta^\star_Q(X-\mu_0)/(1+\beta^\star_Q(X-\mu_0))=1-(1+\beta^\star_Q(X-\mu_0))^{-1}$ gives $\E_Q[(X-\mu_0)/(1+\beta^\star_Q(X-\mu_0))]=(1-r)/\beta^\star_Q$. In the endpoint case, the atom supplies the missing mass and corrects the mean:
\[
\E_{P^\star}[X-\mu_0]=\frac{1-r}{\beta^\star_Q}+(1-r)(x_b-\mu_0)=0\,.
\]
The interior case follows from $r=1$ and $\E_Q[(X-\mu_0)/(1+\beta^\star_Q(X-\mu_0))]=0$. At an endpoint, $1+\beta^\star_Q(x_b-\mu_0)=0$. If $\Pr_Q\{X=x_b\}>0$, the oracle log multiplier would therefore equal $-\infty$ with positive probability. Since its positive part is bounded, this would give $\ell_Q(\beta^\star_Q)=-\infty$, contradicting optimality because the zero bet gives $\ell_Q(0)=0$. Thus, $\Pr_Q\{X=x_b\}=0$, so the likelihood ratio is zero at the additional atom; elsewhere it follows directly from the density in \eqref{eq:projection-law}. Thus, $\KL(Q\|P^\star)=\E_Q[\ln(1+\beta^\star_Q(X-\mu_0))]=\lstar_Q$, by the definition of oracle log growth. For any $P\in\cP_0$ with $\KL(Q\|P)<\infty$, Jensen's inequality gives
\begin{equation}
\label{eq:projection-jensen}
\lstar_Q-\KL(Q\|P)=\E_Q\!\left[\ln\frac{1+\beta^\star_Q(X-\mu_0)}{(\mathrm{d}Q/\mathrm{d}P)(X)}\right]\leq\ln\E_P[1+\beta^\star_Q(X-\mu_0)]=0\,.
\end{equation}
In \eqref{eq:projection-jensen}, canceling the density $\mathrm{d}Q/\mathrm{d}P$ in the expectation of the ratio leaves the integral of $1+\beta^\star_Q(X-\mu_0)$ only over $\{\mathrm{d}Q/\mathrm{d}P>0\}$. Extending that integral to the whole space can only increase it, since the multiplier is nonnegative. The last equality in \eqref{eq:projection-jensen} uses $\E_P[X-\mu_0]=0$, the defining property of the null class $\cP_0$. The case $\KL(Q\|P)=\infty$ is immediate.
\end{proof}

\begin{lemma}[Information bound for early rejection]
\label{lemma:information-early}
Fix $\mu_0\in(0,1)$. There is $C_{\mu_0}<\infty$ such that, for every $Q\in\cQ$, $0<\alpha<1$, $\epsilon\in(0,1)$, and stopping time $\tau$ satisfying $\sup_{P\in\cP_0}\Pr_P\{\tau<\infty\}\leq\alpha$,
\begin{equation}
\label{eq:information-early}
\Pr_Q\left\{\tau\leq\frac{(1-\epsilon)\ln(1/\alpha)}{\lstar_Q}\right\}
\leq\alpha^{\epsilon/2}+\frac{4C_{\mu_0}}{\epsilon^2\ln(1/\alpha)}~.
\end{equation}
\end{lemma}
\begin{proof}
Let $P^\star$ be the null projection from Lemma~\ref{lemma:proj}. Its likelihood ratio is the oracle multiplier, bounded above by $\max\{1/\mu_0,1/(1-\mu_0)\}$. The ratio
\[
\frac{r(\ln r)^2}{r\ln r-r+1}
\]
extends continuously at $r=0$ and $r=1$, with values 0 and 2, respectively. It is therefore bounded by some $C_{\mu_0}<\infty$ on $[0,\max\{1/\mu_0,1/(1-\mu_0)\}]$. Integrating the resulting inequality under $P^\star$, and using $\E_{P^\star}[\mathrm dQ/\mathrm dP^\star]=1$, gives
\begin{equation}
\label{eq:oracle-log-variance}
\E_Q\left[\ln^2(1+\beta_Q^\star(X-\mu_0))\right]\leq C_{\mu_0}\KL(Q\|P^\star)=C_{\mu_0}\lstar_Q~.
\end{equation}
This also holds at endpoint oracle bets, including the additional null atom where the likelihood ratio is zero.

Define the log likelihood increments, their product likelihood, and the deterministic time horizon by
\[
Z_i:=\ln(1+\beta_Q^\star(X_i-\mu_0)),\qquad
L_n:=\exp\left(\sum_{i=1}^n Z_i\right),\qquad
N:=\left\lfloor\frac{(1-\epsilon)\ln(1/\alpha)}{\lstar_Q}\right\rfloor~.
\]
If $N=0$, the event is empty. Otherwise, change of measure on each event $\{\tau=k\}$, $k\leq N$, gives
\[
\Pr_Q\{\tau\leq N,\ \ln L_\tau<(1-\epsilon/2)\ln(1/\alpha)\}
\leq\alpha^{-(1-\epsilon/2)}\Pr_{P^\star}\{\tau\leq N\}\leq\alpha^{\epsilon/2}~.
\]
On the remaining early-rejection event, some $k\leq N$ satisfies $\sum_{i=1}^k(Z_i-\lstar_Q)\geq\epsilon\ln(1/\alpha)/2$. The increments are independent under $Q$ and have mean $\lstar_Q$. Kolmogorov's maximal inequality and \eqref{eq:oracle-log-variance} give
\[
\Pr_Q\left\{\max_{k\leq N}\sum_{i=1}^k(Z_i-\lstar_Q)\geq\frac{\epsilon\ln(1/\alpha)}2\right\}
\leq\frac{4N C_{\mu_0}\lstar_Q}{\epsilon^2\ln^2(1/\alpha)}
\leq\frac{4C_{\mu_0}}{\epsilon^2\ln(1/\alpha)}~.
\]
Adding the two probabilities proves the claim. Only the projection's null validity was used, so the same bound applies to the one-sided setting of Corollary~\ref{cor:one-sided}.
\end{proof}

The next lemma contains the main upper-bound argument. We restrict the prior to an interval around the optimal bets. Concavity preserves most of the optimal growth, and Jensen's inequality converts the interval's prior mass into an additive penalty. This localization appears in \citet[Theorem~1]{Orabona26b} and in \citet[equations~(16)–(18)]{OrabonaJ21}.
\begin{lemma}
\label{lemma:window}
For $\delta\in(0,1/2)$ and $Q\in\cQ$, let
\begin{equation}
\label{eq:localization-interval}
\scI_\delta(\beta^\star_Q):=\{t\beta^\star_Q:1-2\delta\leq t\leq1-\delta\}\,.
\end{equation}
Recall from \eqref{eq:portfolio-test} that $\tau^F_\alpha$ is the first time the portfolio wealth reaches $1/\alpha$. If $F(\scI_\delta(\beta^\star_Q))>0$, then for every $\alpha\in(0,1)$,
\begin{equation}
\label{eq:localization-time}
\E_Q[\tau^F_\alpha]\leq
\frac{\ln(1/\alpha)+\ln(1/F(\scI_\delta(\beta^\star_Q)))+\ln \max\{1/\mu_0,1/(1-\mu_0)\}}{(1-2\delta)\lstar_Q}\,.
\end{equation}
\end{lemma}
\begin{proof}
Let $H_Q$ be $F$ conditioned on the interval $\scI_\delta(\beta^\star_Q)$ in \eqref{eq:localization-interval}, and define the i.i.d. increments
\[
Z_i:=\int\ln(1+\beta (X_i-\mu_0))\,\mathrm{d}H_Q(\beta)\,.
\]
For any $\beta=t\beta^\star_Q$ in $\scI_\delta(\beta^\star_Q)$, we have $1+\beta (X_i-\mu_0)=(1-t)+t(1+\beta^\star_Q(X_i-\mu_0))\geq\delta$. Thus, $1+\beta(X_i-\mu_0)\in[\delta,\max\{1/\mu_0,1/(1-\mu_0)\}]$, and $Z_i\in[\ln\delta,\ln \max\{1/\mu_0,1/(1-\mu_0)\}]$. Concavity gives $\ell_Q(t\beta^\star_Q)\geq t\ell_Q(\beta^\star_Q)+(1-t)\ell_Q(0)=t\lstar_Q$. Hence, $\E_Q[Z_i]\geq(1-2\delta)\lstar_Q>0$. Jensen's inequality gives, pathwise,
\[
\Wealth^F_n\geq F(\scI_\delta(\beta^\star_Q))\int\exp\left(\sum_{i=1}^n\ln(1+\beta (X_i-\mu_0))\right)\,\mathrm{d}H_Q(\beta)
\geq F(\scI_\delta(\beta^\star_Q))\exp\left(\sum_{i=1}^n Z_i\right)\,.
\]
Set $A_Q:=\ln(1/\alpha)+\ln(1/F(\scI_\delta(\beta^\star_Q)))>0$ and $T:=\inf\{n\geq1:\sum_{i=1}^nZ_i\geq A_Q\}$. Taking logarithms in the preceding inequality gives $\tau^F_\alpha\leq T$. For each integer $N$, the sum stopped at $T\wedge N$ is below $A_Q$ if no crossing has occurred; otherwise it is at most $A_Q+\ln \max\{1/\mu_0,1/(1-\mu_0)\}$, since each increment is at most $\ln \max\{1/\mu_0,1/(1-\mu_0)\}$. Since $T\wedge N=\sum_{i=1}^N\indevent\{T\geq i\}$ and $\{T\geq i\}$ depends only on the first $i-1$ observations, independence gives
\[
(1-2\delta)\lstar_Q\E_Q[T\wedge N]\leq\sum_{i=1}^N\E_Q[\indevent\{T\geq i\}Z_i]\leq A_Q+\ln \max\{1/\mu_0,1/(1-\mu_0)\}\,.
\]
This bound is uniform in $N$, so monotone convergence proves $\E_Q[T]<\infty$ and gives \eqref{eq:localization-time}.
\end{proof}

Finally, we prove a lemma that allows us to introduce the logarithmic decay index of the prior.
\begin{lemma}
\label{lemma:mass}
If $a(F)<\infty$, then for every $\delta\in(0,1/2)$ and $\epsilon>0$ there exists $C_{\delta,\epsilon}<\infty$ such that, uniformly over $Q\in\cQ$,
\[
\ln\frac{1}{F(\scI_\delta(\beta^\star_Q))}
\leq(a_Q(F)+\epsilon)\LL(|\beta_Q^\star|)+C_{\delta,\epsilon}\,.
\]
\end{lemma}
\begin{proof}
Recall the limits defining $a_\pm(F)$ in \eqref{eq:decay-indices}. For some $r\in(0,1/e)$, on the side of $\beta_Q^\star$ we have
\[
f(\beta)\geq\frac{1}{|\beta|(\ln(1/|\beta|))^{a_Q(F)+\epsilon}}
\quad (0<|\beta|\leq r)\,.
\]
If $|\beta^\star_Q|\leq r$, then $\scI_\delta(\beta^\star_Q)$ has length $\delta|\beta^\star_Q|$, and every $\beta\in \scI_\delta(\beta^\star_Q)$ satisfies $(1-2\delta)|\beta^\star_Q|\leq|\beta|\leq(1-\delta)|\beta^\star_Q|$. Integrating the preceding lower bound on $f$ over this interval gives
\begin{equation}
\label{eq:interval-mass-lower}
F(\scI_\delta(\beta^\star_Q))\geq
\frac{\delta}{1-\delta}\left(\ln\frac{1}{(1-2\delta)|\beta^\star_Q|}\right)^{-(a_Q(F)+\epsilon)}\,.
\end{equation}
Also,
\[
\ln\frac{1}{(1-2\delta)|\beta^\star_Q|}\leq
\left(1+\ln\frac1{1-2\delta}\right)\max\left\{1,\ln\frac1{|\beta_Q^\star|}\right\}\,.
\]
Taking logarithms in \eqref{eq:interval-mass-lower} and using the preceding inequality proves the stated bound in this case. If $|\beta^\star_Q|>r$, then $\scI_\delta(\beta^\star_Q)$ lies in the compact set $[-(1-\delta)/(1-\mu_0),(1-\delta)/\mu_0]\cap\{|\beta|\geq(1-2\delta)r\}\subset(-1/(1-\mu_0),1/\mu_0)\setminus\{0\}$. On this set $f$ has a positive lower bound independent of $Q$, and $\scI_\delta(\beta^\star_Q)$ has length $\delta|\beta^\star_Q|\geq\delta r$, giving a uniform positive lower bound on $F(\scI_\delta(\beta^\star_Q))$. Enlarging $C_{\delta,\epsilon}$ completes the proof.
\end{proof}

\subsection{A uniform pointwise lower bound}
\label{sec:proof-lower}

The following finite lower bound also applies to an infinite directional index. It also supplies the matching lower bound for the efficiency characterization.

\begin{lemma}[Pointwise lower bound at small oracle bets]
\label{lemma:pointwise-lower}
Let $F$ satisfy Definition~\ref{def:prior}. Fix a finite $a'\in(0,a_+(F))$ and $\epsilon\in(0,1/2)$. There are constants $b_0,A_0,C,\kappa>0$, independent of $Q,\alpha$, such that if $0<\beta_Q^\star\leq b_0$ and $\ln(1/\alpha)\geq A_0$, then
\begin{equation}
\label{eq:finite-pointwise-lower}
\E_Q[\tau_\alpha^F]\geq
\left\lfloor
\frac{(1-\epsilon)\{\ln(1/\alpha)+a'\LL(\beta_Q^\star)\}}{\lstar_Q}
\right\rfloor(1-C\alpha^\kappa)\,.
\end{equation}
The analogous statement holds for negative oracle bets, using $a_-(F)$ and $|\beta_Q^\star|$.
\end{lemma}
\begin{proof}
We split time at a deterministic cutoff $m$ that will be defined in the following. Change of measure controls rejection before $m$; afterwards, concentration and a bound on the prior mass keep the portfolio below the rejection threshold until the horizon in \eqref{eq:finite-pointwise-lower}.

We will choose $b_0,A_0,C,\kappa$ to depend only on $\mu_0,F,a',\epsilon$.

Now, we compute some preliminary inequalities for the quantities of interest.
Fix $Q\in\cQ$ with $0<\beta_Q^\star\leq b_0$. Since $\ell_Q$ is concave and $\ell_Q'(0)=\E_Q[X-\mu_0]$, a positive oracle bet implies
\begin{align*}
\Delta&:=\E_Q[X-\mu_0]>0,\\
v&:=\E_Q[(X-\mu_0)^2]>0\,.
\end{align*}
We take $b_0\leq1/4$ and reduce it further below. Then, the oracle is interior, and its first-order condition gives
\[
\E_Q\left[\frac{X-\mu_0}{1+\beta_Q^\star(X-\mu_0)}\right]=0 \quad \text{ and } \quad
\Delta=\beta_Q^\star\,\E_Q\left[\frac{(X-\mu_0)^2}{1+\beta_Q^\star(X-\mu_0)}\right]\,.
\]
Since $|X-\mu_0|\leq1$,
\begin{equation}
\label{eq:small-bet-relations}
\frac{\beta_Q^\star v}{1+\beta_Q^\star}\leq\Delta\leq\frac{\beta_Q^\star v}{1-\beta_Q^\star}
\quad \text{ and } \quad
\frac{(\beta_Q^\star)^2v}{2(1+\beta_Q^\star)^2}\leq\lstar_Q\leq\frac{(\beta_Q^\star)^2v}{2(1-\beta_Q^\star)^2}.
\end{equation}
The first pair follows by bounding the denominator in the first-order identity between $1-\beta_Q^\star$ and $1+\beta_Q^\star$. For the second pair, subtract the zero quantity $\E_Q[\beta_Q^\star(X-\mu_0)/(1+\beta_Q^\star(X-\mu_0))]$ from $\lstar_Q$. For $|t|\leq\beta_Q^\star$, using $\ln(1+t)-\frac{t}{1+t}=t^2\int_0^1\frac{s}{(1+st)^2}\,\mathrm ds$, we have
\[
\frac{t^2}{2(1+\beta_Q^\star)^2}
\leq\ln(1+t)-\frac{t}{1+t}
\leq\frac{t^2}{2(1-\beta_Q^\star)^2}\,.
\]
Substitute $t=\beta_Q^\star(X-\mu_0)$ and take expectations to obtain the claimed bounds.
These inequalities imply, uniformly as $\beta_Q^\star\to0$ even if $v\to0$,
\begin{equation}
\label{eq:small-bet-asymptotics}
\Delta=\beta_Q^\star v(1+O(\beta_Q^\star)),\qquad
\lstar_Q=\tfrac12(\beta_Q^\star)^2v(1+O(\beta_Q^\star)),\qquad
\frac{\Delta^2}{2v\lstar_Q}=1+O(\beta_Q^\star)\,.
\end{equation}
For all $0<\beta_Q^\star\leq1/4$, using also $v\leq1$, convenient weaker bounds are
\begin{equation}
\label{eq:small-bet-constants}
\beta_Q^\star\Delta\leq5\lstar_Q,\qquad
\frac{\Delta^2}{v\lstar_Q}\geq\frac12,\qquad
\frac{v}{\lstar_Q}\geq\frac1{(\beta_Q^\star)^2},
\quad \text{ and } \quad \lstar_Q\leq1\,.
\end{equation}

\paragraph{Early rejection.}
The null projection from Lemma~\ref{lemma:proj} is
$\mathrm{d}P(x):=\mathrm{d}Q(x)/(1+\beta_Q^\star(x-\mu_0))$. To check that $P\in\cP_0$, the first-order condition gives
\[
\E_Q\left[\frac1{1+\beta_Q^\star(X-\mu_0)}\right]
=1-\beta_Q^\star\E_Q\left[\frac{X-\mu_0}{1+\beta_Q^\star(X-\mu_0)}\right]=1,\qquad
\E_P[X-\mu_0]=0\,.
\]
The multiplier is strictly positive, so the likelihood ratio of the first $n$ observations under $Q$ to those under $P$ is
$L_n:=\prod_{i=1}^n(1+\beta_Q^\star(X_i-\mu_0))$.
Independence and the definition of $P$ give
\begin{align*}
\E_P[L_n^2]
=\left(\E_P[(1+\beta_Q^\star(X-\mu_0))^2]\right)^n
=\left(\E_Q[1+\beta_Q^\star(X-\mu_0)]\right)^n
=(1+\beta_Q^\star\Delta)^n\,.
\end{align*}
Set
\[
m:=\left\lceil \frac{\ln\frac{1}{\alpha}}{20\lstar_Q}\right\rceil,
\]
which is deterministic for fixed $Q,\alpha$. Since $\{\tau_\alpha^F\leq m\}$ is measurable at time $m$, change of measure, Cauchy--Schwarz, and null validity in \eqref{eq:portfolio-test} give
\begin{align}
\Pr_Q\{\tau_\alpha^F\leq m\}
&=\E_P[L_m\indevent\{\tau_\alpha^F\leq m\}]
\leq\sqrt{\E_P[L_m^2]\Pr_P\{\tau_\alpha^F\leq m\}}
\leq\sqrt{\alpha(1+\beta_Q^\star\Delta)^m}
\leq \sqrt{\alpha \exp(m \beta_Q^\star\Delta)} \nonumber\\
&\leq \sqrt{\alpha}\exp(\ln(1/\alpha)/8+5/2)
= e^{5/2}\alpha^{3/8}, \label{eq:pointwise-early}
\end{align}
where the second-to-last inequality uses $1+t\leq e^t$, while for the last inequality, \eqref{eq:small-bet-constants} and the definition of $m$ give $m\beta_Q^\star\Delta\leq5m\lstar_Q\leq\ln(1/\alpha)/4+5$.

\paragraph{Uniform concentration after $m$.}
Define
\[
S_n:=\sum_{i=1}^n(X_i-\mu_0)
\quad \text{ and } \quad
V_n:=\sum_{i=1}^n(X_i-\mu_0)^2\,.
\]
Choose $\theta,\rho\in(0,1/4]$, independently of $Q,\alpha$, small enough that
\begin{equation}
\label{eq:pointwise-approximation-parameters}
\frac{(1+\theta)^2(1+\rho)^2}{1-\theta}\leq1+\frac\epsilon4\,.
\end{equation}
The parameter $\rho$ will determine the range of the quadratic approximation below. Define
\[
\mathcal G:=\bigcap_{n\geq m}
\{|S_n-n\Delta|\leq\theta n\Delta,\
V_n\geq(1-\theta)nv\}\,.
\]
On $\mathcal G$, we have $(1-\theta)n\Delta\leq S_n\leq(1+\theta)n\Delta$ and $V_n\geq(1-\theta)nv$ simultaneously for every $n\geq m$.

If $|Z|\leq1$, $\E_Q[Z]=0$, and $\E_Q[Z^2]\leq v$, then, for $|\lambda|\leq1$, we have
\[
\E_Q[e^{\lambda Z}]\leq1+\lambda\E_Q[Z]+\lambda^2\E_Q[Z^2]
\leq1+\lambda^2v\leq e^{\lambda^2v},
\]
where the first inequality uses $e^t\leq1+t+t^2$ on $[-1,1]$.
Both $X-\E_Q[X]=(X-\mu_0)-\Delta$ and $(X-\mu_0)^2-v$ satisfy these conditions, because their absolute values are at most 1, and
\[
\Var_Q(X)=v-\Delta^2\leq v,\qquad
\Var_Q((X-\mu_0)^2)\leq\E_Q[(X-\mu_0)^4]\leq v\,.
\]
By \eqref{eq:small-bet-relations}, $\theta\Delta/(2v)\leq\theta\beta_Q^\star/(2(1-\beta_Q^\star))\leq1$, and also $\theta/2\leq1$. Independence and the preceding exponential-moment bound therefore make
\[
M_n^\pm:=\exp\left(\pm\frac{\theta\Delta}{2v}(S_n-n\Delta)
-\frac{n\theta^2\Delta^2}{4v}\right),\qquad
\widetilde M_n:=\exp\left(-\frac\theta2(V_n-nv)-\frac{n\theta^2v}{4}\right)
\]
nonnegative $Q$-supermartingales starting from 1. If $\pm(S_n-n\Delta)>\theta n\Delta$ at any $n\geq m$, then $M_n^\pm>\exp(n\theta^2\Delta^2/(4v))\geq\exp(m\theta^2\Delta^2/(4v))$. Similarly, $V_n<(1-\theta)nv$ at any $n\geq m$ implies $\widetilde M_n>\exp(m\theta^2v/4)$.

Since each supermartingale starts from 1, Ville's inequality bounds the probability of crossing each fixed threshold by its reciprocal. Applying this separately to the three events gives
\begin{align*}
\Pr_Q\{\exists n\geq m:S_n-n\Delta>\theta n\Delta\}
&\leq\Pr_Q\left\{\sup_{n\geq0}M_n^+\geq\exp\left(\frac{m\theta^2\Delta^2}{4v}\right)\right\}
\leq\exp\left(-\frac{m\theta^2\Delta^2}{4v}\right),\\
\Pr_Q\{\exists n\geq m:n\Delta-S_n>\theta n\Delta\}
&\leq\Pr_Q\left\{\sup_{n\geq0}M_n^-\geq\exp\left(\frac{m\theta^2\Delta^2}{4v}\right)\right\}
\leq\exp\left(-\frac{m\theta^2\Delta^2}{4v}\right),\\
\Pr_Q\{\exists n\geq m:V_n<(1-\theta)nv\}
&\leq\Pr_Q\left\{\sup_{n\geq0}\widetilde M_n\geq\exp\left(\frac{m\theta^2v}{4}\right)\right\}
\leq\exp\left(-\frac{m\theta^2v}{4}\right)~.
\end{align*}
The event $\mathcal G^c$, the complement of $\mathcal G$, is the union of these three events: at some $n\geq m$, either $S_n$ is too large, $S_n$ is too small, or $V_n$ is too small. Hence, the union bound gives
\begin{align}
\label{eq:pointwise-good-event}
\Pr_Q\{\mathcal G^c\}
&\leq\Pr_Q\{\exists n\geq m:S_n-n\Delta>\theta n\Delta\}
+\Pr_Q\{\exists n\geq m:n\Delta-S_n>\theta n\Delta\}\nonumber\\
&\quad+\Pr_Q\{\exists n\geq m:V_n<(1-\theta)nv\}\nonumber\\
&\leq2\exp\left(-\frac{m\theta^2\Delta^2}{4v}\right)
+\exp\left(-\frac{m\theta^2v}{4}\right)
\leq2\alpha^{\theta^2/160}+\alpha^{\theta^2/(80(\beta_Q^\star)^2)}
\leq3\alpha^{\theta^2/160},
\end{align}
where the second-to-last inequality uses $m\geq\ln(1/\alpha)/(20\lstar_Q)$ and \eqref{eq:small-bet-constants}, while the last inequality follows from $\beta_Q^\star\leq1/4$ and $0<\alpha<1$.

\paragraph{A deterministic bound on the portfolio.}
Set
\begin{equation}
\label{eq:pointwise-horizon}
H:=\ln(1/\alpha)+a'\LL(\beta_Q^\star),\quad \text{ and } \quad
N:=\left\lfloor\frac{(1-\epsilon)H}{\lstar_Q}\right\rfloor\,.
\end{equation}
The floor and ceiling give
\[
N-m\geq\frac{(19/20-\epsilon)\ln(1/\alpha)}{\lstar_Q}-2>0
\]
when $\ln(1/\alpha)\geq5$, since $H\geq\ln(1/\alpha)$, $\epsilon<1/2$, and $\lstar_Q\leq1$.

Fix $m<n\leq N$, work on $\mathcal G$, and define the concave function
\[
g_n(\beta):=\ln\Wealth_n(\beta),
\]
where $g_n(0)=0$. On $[-\rho,\rho]$, the second derivative of $\ln(1+t)$ is at most $-1/(1+\rho)^2$. Taylor's theorem at zero therefore gives
$\ln(1+t)\leq t-t^2/(2(1+\rho)^2)$.
Apply this inequality with $t=\beta(X_i-\mu_0)$ and sum over $i$. Since $|X_i-\mu_0|\leq1$ and we are on $\mathcal{G}$, for $0\leq\beta\leq\rho$ we obtain
\[
g_n(\beta)
\leq\beta S_n-\frac{\beta^2V_n}{2(1+\rho)^2}
\leq \beta S_n-\frac{\beta^2 (1-\theta) n v}{2(1+\rho)^2}
= nv'\left(x_n\beta-\frac{\beta^2}{2}\right),
\]
where
\[
v':=\frac{(1-\theta)v}{(1+\rho)^2}
\quad \text{ and } \quad x_n:=\frac{S_n}{nv'}\,.
\]
Here, $x_n$ is the unconstrained maximizer of this quadratic upper bound on log-wealth.
Completing the square gives
\[
nv'\left(x_n\beta-\frac{\beta^2}{2}\right)
=K_n-\frac{nv'}2(\beta-x_n)^2\leq K_n,
\quad \text{ where} \quad
K_n:=\frac{S_n^2}{2nv'}=\frac{nv'x_n^2}{2}\,.
\]
On $\mathcal G$, the bounds on $S_n$ and \eqref{eq:small-bet-relations} imply
\[
0<x_n\leq
\frac{(1+\theta)(1+\rho)^2}{(1-\theta)(1-\beta_Q^\star)}\beta_Q^\star
\leq C_0\beta_Q^\star,
\quad \text{ where } \quad
C_0:=\frac{4(1+\theta)(1+\rho)^2}{3(1-\theta)}\,.
\]
Reducing $b_0$ so that $2C_0b_0\leq\rho$ ensures $2x_n\leq\rho$. Moreover,
\[
K_n\leq n\lstar_Q\,
\frac{(1+\theta)^2(1+\rho)^2}{1-\theta}
\frac{\Delta^2}{2v\lstar_Q}\,.
\]
By \eqref{eq:pointwise-approximation-parameters}, the first ratio is at most $1+\epsilon/4$; by \eqref{eq:small-bet-asymptotics}, the second ratio tends uniformly to 1 as $\beta_Q^\star\to0$. We can therefore reduce $b_0$, independently of $Q,\alpha$, so that
\[
\frac{\Delta^2}{2v\lstar_Q}\leq\frac{1+\epsilon/2}{1+\epsilon/4}\,.
\]
The product of the two ratios is then at most $1+\epsilon/2$, giving
\[
K_n\leq(1+\epsilon/2)n\lstar_Q
\leq(1+\epsilon/2)(1-\epsilon)H
=(1-\epsilon/2-\epsilon^2/2)H
\leq(1-\epsilon/2)H\,.
\]
We have established, simultaneously for all $m<n\leq N$ on $\mathcal G$,
\begin{equation}
\label{eq:pointwise-quadratic}
0<x_n\leq C_0\beta_Q^\star,\qquad
2x_n\leq\rho,\qquad
K_n\leq(1-\epsilon/2)H\,.
\end{equation}

We now bound the wealth above its initial value of 1. First, consider $\beta$ outside $(0,2x_n)$. For $\beta<0$, the inequality $\ln(1+t)\leq t$ gives $g_n(\beta)\leq\beta S_n\leq0$. For $2x_n\leq\beta\leq\rho$, the quadratic bound gives
$g_n(\beta)\leq nv'\beta(x_n-\beta/2)\leq0$.
For $\beta\geq\rho$, concavity gives
\[
g_n(\rho)\geq\frac{\rho}{\beta}g_n(\beta)
+\left(1-\frac{\rho}{\beta}\right)g_n(0)
=\frac{\rho}{\beta}g_n(\beta)\,.
\]
So, using the quadratic bound at $\rho$ together with $2x_n\leq \rho$, gives $g_n(\beta)\leq(\beta/\rho)g_n(\rho)\leq0$. Endpoint zeros are interpreted as log wealth $-\infty$. Since $F$ has total mass 1 and $e^{g_n(\beta)}\leq1$ outside $(0,2x_n)$, we have
\[
\Wealth_n^F
=1+\int\!\big(e^{g_n(\beta)}-1\big)\,\mathrm dF(\beta)
\leq1+\int_0^{2x_n}\!\big(e^{g_n(\beta)}-1\big)f(\beta)\,\mathrm d\beta\,.
\]

For $0\leq\beta\leq2x_n$, the quadratic bound and $K_n=nv'x_n^2/2$ give
\[
g_n(\beta)\leq K_n\left(\frac{2\beta}{x_n}-\frac{\beta^2}{x_n^2}\right),
\qquad
0\leq\frac{2\beta}{x_n}-\frac{\beta^2}{x_n^2}
=1-\left(1-\frac{\beta}{x_n}\right)^2\leq1\,.
\]
By convexity, the function $t\mapsto e^{tK_n}$ lies below the chord joining its values at 0 and 1. Thus,
\[
e^{g_n(\beta)}-1
\leq\exp\left(K_n\left(\frac{2\beta}{x_n}-\frac{\beta^2}{x_n^2}\right)\right)-1
\leq(e^{K_n}-1)\left(\frac{2\beta}{x_n}-\frac{\beta^2}{x_n^2}\right)\,.
\]

Since $a'<a_+(F)$, the limit in \eqref{eq:decay-indices} lets us choose
\[
\beta_1\in(0,1/e)
\quad\text{such that}\quad
\beta f(\beta)
\leq\left[\ln\frac{1}{\beta}\right]^{-a'}
\qquad(0<\beta\leq\beta_1)\,.
\]
Reduce $b_0$ so that $2C_0b_0\leq\beta_1$ and $b_0\leq(2C_0)^{-2}$. Then
\[
\ln\frac{1}{2x_n}
\geq\ln\frac{1}{\beta_Q^\star}-\ln(2C_0)
\geq\frac{1}{2}\ln\frac{1}{\beta_Q^\star}\,.
\]
For every $0<\beta\leq2x_n$, the density envelope therefore gives
\begin{align*}
\beta f(\beta)
\leq\left[\ln\frac{1}{\beta}\right]^{-a'}
\leq\left[\ln \frac{1}{2x_n}\right]^{-a'}
\leq2^{a'}\left[\ln\frac{1}{\beta_Q^\star}\right]^{-a'}
=2^{a'}e^{-a'\LL(\beta_Q^\star)},
\end{align*}
where the last equality uses $\beta_Q^\star\leq1/4<1/e$, so that $\LL(\beta_Q^\star)=\ln\ln(1/\beta_Q^\star)$.

Combining the preceding bounds gives
\begin{align*}
\Wealth_n^F
&\leq1+(e^{K_n}-1)\int_0^{2x_n}
\left(\frac{2\beta}{x_n}-\frac{\beta^2}{x_n^2}\right)f(\beta)\,\mathrm d\beta
\leq1+2^{a'}e^{-a'\LL(\beta_Q^\star)}(e^{K_n}-1)
\int_0^{2x_n}\left(\frac2{x_n}-\frac{\beta}{x_n^2}\right)\,\mathrm d\beta\\
&=1+2^{a'+1}e^{-a'\LL(\beta_Q^\star)}(e^{K_n}-1),
\end{align*}
where the integral is
\[
\int_0^{2x_n}\left(\frac2{x_n}-\frac{\beta}{x_n^2}\right)\,\mathrm d\beta
=\left[\frac{2\beta}{x_n}-\frac{\beta^2}{2x_n^2}\right]_0^{2x_n}
=4-2=2\,.
\]
The factor multiplying $f(\beta)$ vanishes linearly at zero, canceling the $1/\beta$ in the density bound. Thus, the integral is finite even when $a'<1$.

By $e^{K_n}-1\leq e^{K_n}$, the definition of $H$, and \eqref{eq:pointwise-quadratic},
\[
\alpha e^{-a'\LL(\beta_Q^\star)}(e^{K_n}-1)
\leq\alpha e^{K_n-a'\LL(\beta_Q^\star)}
=e^{K_n-H}\leq e^{-\epsilon H/2}\,.
\]
Multiplying the wealth bound by $\alpha$, we conclude that
\begin{equation}
\label{eq:pointwise-wealth-envelope}
\alpha\Wealth_n^F
\leq\alpha+2^{a'+1}e^{-\epsilon H/2}
\qquad(m<n\leq N)\,.
\end{equation}
Choose $A_0$ large enough that
\[
A_0\geq5
\quad\text{and}\quad
e^{-A_0}+2^{a'+1}e^{-\epsilon A_0/2}<1\,.
\]
Such a choice exists because both exponential terms tend to zero as $A_0\to\infty$. Whenever $H\geq\ln(1/\alpha)\geq A_0$, \eqref{eq:pointwise-wealth-envelope} gives
\[
\alpha\Wealth_n^F\leq e^{-A_0}+2^{a'+1}e^{-\epsilon A_0/2}<1\,.
\]
Thus, on $\mathcal G$, no rejection occurs at any $m<n\leq N$.

\paragraph{Conclusion.}
The preceding bound gives
$\{\tau_\alpha^F\leq N\}\subseteq\{\tau_\alpha^F\leq m\}\cup\mathcal G^c$.
Combining \eqref{eq:pointwise-early} and \eqref{eq:pointwise-good-event},
\begin{equation}
\label{eq:pointwise-lower-tail}
\Pr_Q\{\tau_\alpha^F\leq N\}
\leq e^{5/2}\alpha^{3/8}+3\alpha^{\theta^2/160}
\leq C\alpha^\kappa,
\qquad C:=e^{5/2}+3,\quad\kappa:=\theta^2/160\,.
\end{equation}
Finally,
\[
\E_Q[\tau_\alpha^F]\geq N\Pr_Q\{\tau_\alpha^F>N\}
\geq N(1-C\alpha^\kappa),
\]
which proves \eqref{eq:finite-pointwise-lower}. For a negative oracle bet, replace $X,\mu_0,\beta$ by $1-X,1-\mu_0,-\beta$ to obtain the analogous result with $a_+(F)$ replaced by $a_-(F)$.
\end{proof}

\subsection{Proof of Theorem~\ref{thm:pointwise}}
\label{sec:proof-pointwise}

We use the following elementary criterion to combine an expected-time upper bound with control of early rejection. Its comparison form will also apply to confidence-sequence widths and to the regret-based test in Appendix~\ref{app:oj-mixture}.

\begin{lemma}[A uniform $L^1$ criterion]
\label{lemma:l1-criterion}
Let $Z_{\alpha,Q}$ be nonnegative integrable random variables under $Q$, indexed by $0<\alpha<1$ and a family of laws $\mathcal Q$, such that
\[
\limsup_{\alpha\downarrow0}\ \sup_{Q\in\mathcal Q}\E_Q[Z_{\alpha,Q}]\leq1~.
\]
Then $\sup_{Q\in\mathcal Q}\E_Q[|Z_{\alpha,Q}-1|]\to0$ as $\alpha\downarrow0$ if, for every $\eta\in(0,1)$,
\[
\sup_{Q\in\mathcal Q}\Pr_Q\{Z_{\alpha,Q}<1-\eta\}\longrightarrow0~.
\]
The same conclusion holds if there are random variables $0\leq U_{\alpha,Q}\leq Z_{\alpha,Q}$ almost surely such that
\[
\sup_{Q\in\mathcal Q}\E_Q[|U_{\alpha,Q}-1|]\longrightarrow0\qquad(\alpha\downarrow0)~.
\]
\end{lemma}
\begin{proof}
For every $x\geq0$, $|x-1|=x-1+2(1-x)_+$. Also,
\[
(1-x)_+\leq\eta+\indevent\{x<1-\eta\}~,
\]
because the left side is at most $\eta$ when $x\geq1-\eta$ and at most 1 otherwise. Taking expectations gives
\[
\E_Q[|Z_{\alpha,Q}-1|]
\leq\E_Q[Z_{\alpha,Q}]-1+2\eta+2\Pr_Q\{Z_{\alpha,Q}<1-\eta\}~.
\]
Taking the supremum over $Q$ and then the upper limit as $\alpha\downarrow0$ bounds the left side by $2\eta$: the expectation minus 1 has upper limit at most zero, and the probability tends uniformly to zero. Sending $\eta\downarrow0$ proves the first conclusion.

For the comparison form, $(1-Z_{\alpha,Q})_+\leq(1-U_{\alpha,Q})_+\leq|U_{\alpha,Q}-1|$. The same identity therefore gives
\[
\E_Q[|Z_{\alpha,Q}-1|]\leq\E_Q[Z_{\alpha,Q}]-1+2\E_Q[|U_{\alpha,Q}-1|]~.
\]
After taking the supremum and upper limit, the right side is at most zero. Nonnegativity proves the claimed convergence.
\end{proof}

\begin{proof}[Proof of Theorem~\ref{thm:pointwise}]
Fix $\delta\in(0,1/2)$ and $\zeta>0$. Enlarge the finite constant $C_{\delta,\zeta}$ from Lemma~\ref{lemma:mass}, if necessary, so that it is nonnegative, and define
\[
C'_{\delta,\zeta}:=C_{\delta,\zeta}
+\ln\max\{1/\mu_0,1/(1-\mu_0)\}<\infty\,.
\]
This constant is independent of $Q,\alpha$. So, combining Lemma~\ref{lemma:window} with Lemma~\ref{lemma:mass} gives
\begin{equation}
\label{eq:pointwise-upper}
\E_Q[\tau_\alpha^F]\leq
\frac{\ln(1/\alpha)+(a_Q(F)+\zeta)\LL(|\beta_Q^\star|)+C'_{\delta,\zeta}}{(1-2\delta)\lstar_Q},
\end{equation}
uniformly in $Q,\alpha$. Dividing by $\widetilde B_\alpha^F(Q)$ gives
\begin{align*}
\frac{\E_Q[\tau_\alpha^F]}{\widetilde B_\alpha^F(Q)}
&\leq\frac{1}{1-2\delta}
\left(1+\frac{\zeta\LL(|\beta_Q^\star|)+C'_{\delta,\zeta}}
{\ln(1/\alpha)+a_Q(F)\LL(|\beta_Q^\star|)}\right)
\leq\frac{1}{1-2\delta}
\left(1+\zeta+\frac{C'_{\delta,\zeta}}{\ln(1/\alpha)}\right),
\end{align*}
where in the second inequality, $a_Q(F)\geq1$ makes the denominator at least $\LL(|\beta_Q^\star|)$, and it is also at least $\ln(1/\alpha)>0$. Since $C'_{\delta,\zeta}$ is finite for fixed $\delta,\zeta$, its contribution tends to zero uniformly over $Q$ as $\alpha\downarrow0$. Hence,
\[
\limsup_{\alpha\downarrow0}\ \sup_{Q\in\cQ}
\frac{\E_Q[\tau_\alpha^F]}{\widetilde B_\alpha^F(Q)}
\leq\frac{1+\zeta}{1-2\delta}\,.
\]
Now send $\delta,\zeta\downarrow0$ to obtain
\[
\limsup_{\alpha\downarrow0}\ \sup_{Q\in\cQ}
\frac{\E_Q[\tau_\alpha^F]}{\widetilde B_\alpha^F(Q)}\leq1~.
\]
The constant $C'_{\delta,\zeta}$ need not remain bounded in this last step, because we have already taken the limit in $\alpha$ with $\delta,\zeta$ fixed.

To control the lower tail, fix $\eta\in(0,1)$ and choose $\epsilon\in(0,1/2)$ so that $(1-\epsilon)^2>1-\eta$. Define
\[
a'_s:=(1-\epsilon)a_s(F)\in(0,a_s(F)),\qquad s\in\{+,-\}~.
\]
Apply Lemma~\ref{lemma:pointwise-lower} separately on the two sides. Choose $b_0$ as the smaller betting cutoff, $A_0$ as the larger confidence cutoff, $C$ as the larger constant, and $\kappa$ as the smaller exponent. Since $0<\alpha<1$, decreasing the exponent increases $\alpha^\kappa$. Thus, for $\ln(1/\alpha)\geq A_0$, these choices give a common probability bound in \eqref{eq:pointwise-lower-tail}, with constants independent of $Q,\alpha$.

For $0<|\beta_Q^\star|\leq b_0$, let $s$ be the side containing $\beta_Q^\star$. The choice $a'_s=(1-\epsilon)a_Q(F)$ gives
\begin{align*}
\ln(1/\alpha)+a'_s\LL(|\beta_Q^\star|)
&=\ln(1/\alpha)+(1-\epsilon)a_Q(F)\LL(|\beta_Q^\star|)\\
&\geq(1-\epsilon)\{\ln(1/\alpha)+a_Q(F)\LL(|\beta_Q^\star|)\}~.
\end{align*}
Consequently, the lemma's time horizon satisfies
\[
N:=\left\lfloor\frac{(1-\epsilon)\{\ln(1/\alpha)+a'_s\LL(|\beta_Q^\star|)\}}{\lstar_Q}\right\rfloor
\geq\left\lfloor(1-\epsilon)^2\widetilde B_\alpha^F(Q)\right\rfloor~.
\]
One factor $1-\epsilon$ comes from the horizon and the other from the preceding comparison of the two numerators. The floor loses at most one observation, whose normalized cost is bounded by \eqref{eq:oracle-growth-bound}:
\[
\frac1{\widetilde B_\alpha^F(Q)}\leq
\frac{\ln\max\{1/\mu_0,1/(1-\mu_0)\}}{\ln(1/\alpha)}~.
\]
It follows that
\[
\frac{N}{\widetilde B_\alpha^F(Q)}\geq(1-\epsilon)^2-
\frac{\ln\max\{1/\mu_0,1/(1-\mu_0)\}}{\ln(1/\alpha)}~.
\]
Since $(1-\epsilon)^2>1-\eta$, this ratio is at least $1-\eta$ whenever
\[
\ln(1/\alpha)\geq
\frac{\ln\max\{1/\mu_0,1/(1-\mu_0)\}}{(1-\epsilon)^2-(1-\eta)}~.
\]
This cutoff is independent of $Q$. Equation~\eqref{eq:pointwise-lower-tail} therefore gives
\[
\Pr_Q\{\tau_\alpha^F\leq(1-\eta)\widetilde B_\alpha^F(Q)\}\leq C\alpha^\kappa~.
\]

For $|\beta_Q^\star|>b_0$, the adaptation term is bounded by $a(F)\LL(b_0)$. Hence, whenever
\[
\ln(1/\alpha)\geq\frac{2(1-\eta)a(F)\LL(b_0)}\eta~,
\]
we have $(1-\eta)a(F)\LL(b_0)\leq(\eta/2)\ln(1/\alpha)$ and therefore
\[
(1-\eta)\widetilde B_\alpha^F(Q)\leq\frac{(1-\eta/2)\ln(1/\alpha)}{\lstar_Q}
\]
uniformly on this range. Apply Lemma~\ref{lemma:information-early} with $\epsilon=\eta/2$ to obtain
\[
\Pr_Q\{\tau_\alpha^F\leq(1-\eta)\widetilde B_\alpha^F(Q)\}
\leq\alpha^{\eta/4}+\frac{16C_{\mu_0}}{\eta^2\ln(1/\alpha)}~.
\]
For all $\alpha$ satisfying the three cutoffs above, the supremum over $Q\in\cQ$ is at most the maximum of this bound and $C\alpha^\kappa$. Both tend to zero as $\alpha\downarrow0$, proving
\begin{equation}
\label{eq:pointwise-early-uniform}
\lim_{\alpha\downarrow0}\ \sup_{Q\in\cQ}\Pr_Q\{\tau_\alpha^F\leq(1-\eta)\widetilde B_\alpha^F(Q)\}=0~.
\end{equation}

Apply Lemma~\ref{lemma:l1-criterion} to the normalized rejection times
\[
Z_{\alpha,Q}:=\frac{\tau_\alpha^F}{\widetilde B_\alpha^F(Q)}\qquad(Q\in\cQ)~.
\]
They are nonnegative and integrable by \eqref{eq:pointwise-upper}, and the preceding upper-bound argument gives the required uniform upper limit of their expectations. Equation~\eqref{eq:pointwise-early-uniform} supplies the lower-tail condition for every $\eta\in(0,1)$, since the event with a strict inequality is contained in the displayed event. The lemma proves \eqref{eq:pointwise-l1}.
\end{proof}

\subsection{Proof of Corollary~\ref{thm:main}}
\label{sec:proof-upper}

\begin{proof}
For the upper bound, suppose $a(F)<\infty$. Theorem~\ref{thm:pointwise} gives $\E_Q[\tau_\alpha^F]/\widetilde B_\alpha^F(Q)\to1$ uniformly over $Q$ as $\alpha\to0$. Since $1\leq a_Q(F)\leq a(F)$,
\[
\frac{\widetilde B_\alpha^F(Q)}{B_\alpha(Q)}
=\frac{\ln(1/\alpha)+a_Q(F)\LL(|\beta_Q^\star|)}
{\ln(1/\alpha)+\LL(|\beta_Q^\star|)}\leq a(F)\,.
\]
Hence, $\eff(\tau^F)\leq a(F)$. The same upper bound is automatic for $a(F)=\infty$.

For the lower bound, choose a side with index $a(F)$, taking the positive side without loss of generality. For $0<\Delta<1-\mu_0$, let
\[
Q_\Delta:=\mathrm{Ber}(\mu_0+\Delta)
\quad\text{and}\quad
\sigma_0^2:=\mu_0(1-\mu_0)\,.
\]
Then
\begin{equation}
\label{eq:bernoulli-asymptotics}
\begin{aligned}
\beta_{Q_\Delta}^\star&=\Delta/\sigma_0^2,\\
\lstar_{Q_\Delta}
&=\inf_{P\in\cP_0}\KL(Q_\Delta\|P)
=\KLBern(\mu_0+\Delta;\mu_0)
=\frac{\Delta^2}{2\sigma_0^2}(1+O(\Delta)),
\end{aligned}
\end{equation}
where $\KLBern(p;q)$ denotes the Kullback--Leibler divergence between Bernoulli distributions with parameters $p$ and $q$:
\[
\KLBern(p;q):=p\ln\frac{p}{q}+(1-p)\ln\frac{1-p}{1-q},
\qquad p\in[0,1],\quad q\in(0,1),
\]
with zero-weight terms interpreted as zero.
The first equality for $\lstar_{Q_\Delta}$ follows from Lemma~\ref{lemma:proj}. To evaluate the infimum, note that every $P\in\cP_0$ satisfies
\[
\Pr_P\{X=1\}\leq\E_P[X]=\mu_0,
\qquad
\Pr_P\{X=0\}\leq\E_P[1-X]=1-\mu_0\,.
\]
Since $Q_\Delta$ is supported on $\{0,1\}$,
\begin{align*}
\KL(Q_\Delta\|P)
&=(\mu_0+\Delta)\ln\frac{\mu_0+\Delta}{\Pr_P\{X=1\}}
+(1-\mu_0-\Delta)\ln\frac{1-\mu_0-\Delta}{\Pr_P\{X=0\}}\\
&\geq\KLBern(\mu_0+\Delta;\mu_0)\,.
\end{align*}
Equality holds for $P=\mathrm{Ber}(\mu_0)$, which belongs to $\cP_0$. Thus, the Bernoulli divergence is the attained infimum over the full null class.

For sufficiently small $\Delta$, set
\[
\alpha_\Delta:=\frac{1}{\LL(\beta_{Q_\Delta}^\star)}\,.
\]
Here, $\alpha_\Delta$ is the prescribed upper bound on the false-alarm probability under every null law, and we evaluate the test at this level under $Q_\Delta$. As $\Delta\downarrow0$ with $\mu_0$ fixed, $Q_\Delta$ approaches the null law $\mathrm{Ber}(\mu_0)$ while the permitted false-alarm probability $\alpha_\Delta$ also tends to zero. Our choice makes $\alpha_\Delta$ decrease slowly enough that the confidence cost is negligible compared with the adaptation cost:
\[
\ln(1/\alpha_\Delta)=o\!\left(\LL(\beta_{Q_\Delta}^\star)\right)
\qquad(\Delta\downarrow0)\,.
\]

Fix any finite $a'\in(0,a(F))$ and $\epsilon\in(0,1/2)$. The constants in Lemma~\ref{lemma:pointwise-lower} are uniform over $Q,\alpha$, so we can apply it along the changing alternatives and levels $(Q_\Delta,\alpha_\Delta)$: for all sufficiently small $\Delta$, we have $0<\beta_{Q_\Delta}^\star\leq b_0$ and $\ln(1/\alpha_\Delta)\geq A_0$. Enlarge $A_0$ so that $0\leq1-C\alpha_\Delta^\kappa\leq1$. Dividing \eqref{eq:finite-pointwise-lower} by $B_{\alpha_\Delta}(Q_\Delta)$ and using $\lfloor t\rfloor\geq t-1$ gives
\[
\frac{\E_{Q_\Delta}[\tau_{\alpha_\Delta}^F]}{B_{\alpha_\Delta}(Q_\Delta)}
\geq(1-\epsilon)(1-C\alpha_\Delta^\kappa)
\frac{\ln(1/\alpha_\Delta)+a'\LL(\beta_{Q_\Delta}^\star)}
{\ln(1/\alpha_\Delta)+\LL(\beta_{Q_\Delta}^\star)}
-\frac{1}{B_{\alpha_\Delta}(Q_\Delta)}~.
\]
The floor loses at most one observation, and \eqref{eq:oracle-growth-bound} gives
\[
\frac{1}{B_{\alpha_\Delta}(Q_\Delta)}
\leq\frac{\ln\max\{1/\mu_0,1/(1-\mu_0)\}}{\ln(1/\alpha_\Delta)}\longrightarrow0~.
\]
In the lower bound, dividing the numerator and denominator of the fraction containing $a'$ by $\LL(\beta_{Q_\Delta}^\star)$ shows that this fraction tends to $a'$, while $1-C\alpha_\Delta^\kappa\to1$. Hence,
\[
\liminf_{\Delta\downarrow0}\frac{\E_{Q_\Delta}[\tau_{\alpha_\Delta}^F]}{B_{\alpha_\Delta}(Q_\Delta)}\geq(1-\epsilon)a'~.
\]
For every $\alpha_0>0$, all sufficiently small $\Delta$ satisfy $\alpha_\Delta\leq\alpha_0$, so the supremum in Definition~\ref{def:eff} is at least $(1-\epsilon)a'$. Sending $\alpha_0\downarrow0$ and then $\epsilon\downarrow0$ gives $\eff(\tau^F)\geq a'$. Finally, let $a'\uparrow a(F)$ if $a(F)<\infty$, or $a'\to\infty$ otherwise. This proves the lower bound in both cases.
\end{proof}

\subsection{Proofs of Corollaries~\ref{cor:examples},~\ref{cor:one-sided}, ~\ref{cor:combine}, and Proposition~\ref{prop:form}}

\begin{proof}[Proof of Corollary~\ref{cor:examples}]
All three priors satisfy the local positivity part of Definition~\ref{def:prior}, so we only need to compute $a(F)$ and apply Corollary~\ref{thm:main}. Recall that $L_1(\beta)=\ln(e\max\{1/\mu_0,1/(1-\mu_0)\}/|\beta|)$ and $L_2(\beta)=\ln(eL_1(\beta))$.

For (i), the density of $\beta$ is $f(\beta)=\mu_0(1-\mu_0)\,b(\mu_0+\beta\mu_0(1-\mu_0))$, where $b$ is the Beta density, which is continuous and positive at $\mu_0\in(0,1)$. Hence, $\ln\frac{1}{h(\beta)}=\ln\frac{1}{|\beta|}+O(1)$ and $a_\pm(F)=+\infty$.

For (ii), $\ln\frac{1}{h(\beta)}=a\ln L_1(\beta)+O(1)$ and $\ln L_1(\beta)/\ln\ln\frac{1}{|\beta|}\to1$, so $a_\pm(F)=a$.

For (iii), $\ln\frac{1}{h(\beta)}=\ln L_1(\beta)+(1+\eta)\ln L_2(\beta)-\ln C_\eta$ and $\ln L_2=\ln(1+\ln L_1)=o(\ln L_1)$, so $a_\pm(F_\eta)=1$.
\end{proof}

\begin{proof}[Proof of Corollary~\ref{cor:one-sided}]
For nonnegative bets, $\E_P[1+\beta(X-\mu_0)]\leq1$ under the one-sided null, so the portfolio wealth is a nonnegative supermartingale and the test is valid. Since $\ell_Q'(0)=\E_Q[X]-\mu_0>0$, the oracle bet is positive, so $\scI_\delta(\beta^\star_Q)\subset(0,1/\mu_0)$ and the bounds in Lemmas~\ref{lemma:window} and~\ref{lemma:mass} use only $a_+(F)$. The projection in Lemma~\ref{lemma:proj} has mean $\mu_0$ and belongs to the one-sided null, giving the universal lower bound. Lemma~\ref{lemma:pointwise-lower} uses only the positive-side density bound and the null projection with mean $\mu_0$. Lemma~\ref{lemma:information-early} also applies because that projection belongs to the one-sided null. Thus, the proof of both the expected-time and $L^1$ limits applies unchanged. The positive Bernoulli alternatives in the proof of Corollary~\ref{thm:main} give the matching efficiency lower bound.
\end{proof}

\begin{proof}[Proof of Corollary~\ref{cor:combine}]
Let $f_j$ be the component densities. The combined density $f(\beta)=\sum_{j=1}^k w_jf_j(\beta)$ satisfies
\begin{equation}
\label{eq:combined-density}
\max_j\{w_j|\beta|f_j(\beta)\}
\leq |\beta|f(\beta)\leq k\max_j\{w_j|\beta|f_j(\beta)\}\,.
\end{equation}
The negative logarithm of the maximum in \eqref{eq:combined-density} is $\min_j\{\ln(1/(|\beta|f_j(\beta)))+\ln(1/w_j)\}$. Taking negative logarithms in \eqref{eq:combined-density}, we see that $\ln(1/(|\beta|f(\beta)))$ differs from this minimum by at most $\ln k$. Divide by $\ln\ln(1/|\beta|)$ and pass to each one-sided limit. Since $k$ and the weights are fixed, $a_\pm(F)=\min_j a_\pm(F_j)$, including infinite logarithmic decay indices. Finally, apply Corollary~\ref{thm:main} and Theorem~\ref{thm:pointwise} when both minima are finite.
\end{proof}

\begin{proof}[Proof of Proposition~\ref{prop:form}]
If $g\leq0$, the assumption also holds with any $g'\in(0,\lstar_Q)$ in place of $g$, and letting $g'\to0$ gives $\eff(\tau)=\infty$; so, assume $g>0$ and fix $\eta>0$. $Q$-a.s., there exists $n_0(\omega)<\infty$ such that $\ln \Wealth_n\leq n(g+\eta)$ for all $n\geq n_0$. For sufficiently small $\alpha$, the threshold $1/\alpha$ exceeds all the finitely many wealth values before $n_0$. Thus, $\tau_\alpha\geq n_0$ and
\[
\tau_\alpha\geq\frac{\ln(1/\alpha)}{g+\eta}\,.
\]
This inequality also holds when $\tau_\alpha=\infty$. So, $\liminf_{\alpha\to0}\tau_\alpha/\ln\frac1\alpha\geq\frac{1}{g+\eta}$ $Q$-a.s., and by Fatou's lemma $\liminf_{\alpha\to0}\E_Q[\tau_\alpha]/\ln\frac1\alpha\geq\frac{1}{g+\eta}$. Since $B_\alpha(Q)\sim\ln\frac1\alpha/\lstar_Q$ as $\alpha\to0$, we get $\eff(\tau)\geq\lstar_Q/(g+\eta)$, and we conclude by letting $\eta\to0$.

For the Hoeffding process in \eqref{eq:hoeffding-process}, let $S_n:=\sum_{i=1}^n(X_i-\mu_0)$. Then, $\exp(\beta S_n-n\beta^2/8)$ is a supermartingale under every $P\in\cP_0$ by Hoeffding's lemma, so $(\tau^{\mathrm H}_\alpha)$ is a level-$\alpha$ family. Moreover, completing the square gives $\beta S_n-n\beta^2/8=2S_n^2/n-(n/8)(\beta-4S_n/n)^2\leq2S_n^2/n$, so $\frac1n\ln \Wealth^{\mathrm H}_n\leq2(S_n/n)^2\to2(\E_Q[X-\mu_0])^2$ $Q$-a.s. To see that the growth gap is unbounded, take $Q_\Delta:=\delta_{\mu_0+\Delta}$ for $0<\Delta<1-\mu_0$. Then, $\beta^\star_{Q_\Delta}=1/\mu_0$ and
\[
\frac{\lstar_{Q_\Delta}}{2(\E_{Q_\Delta}[X-\mu_0])^2}
=\frac{\ln(1+\Delta/\mu_0)}{2\Delta^2}\longrightarrow\infty\,.
\]
Applying the proposition separately at each fixed $Q_\Delta$ proves infinite efficiency.
\end{proof}

\subsection{Proof of Proposition~\ref{prop:near-null-scales}}
\label{sec:proof-near-null-scales}

\begin{proof}
We first show that the oracle bet tends to zero. For each fixed $c\in(0,1)$, the bound $1+c(X-\mu_0)\leq1+c$ gives
\begin{align*}
\ell_{Q_k}'(c)
=\E_{Q_k}[X-\mu_0]-c\E_{Q_k}\left[\frac{(X-\mu_0)^2}{1+c(X-\mu_0)}\right]
\leq\E_{Q_k}[X]-\mu_0-\frac{c}{1+c}\E_{Q_k}[(X-\mu_0)^2]<0
\end{align*}
eventually, since the mean tends to $\mu_0$ and the second moment is bounded below. Concavity gives $\beta_{Q_k}^\star<c$ eventually; applying the same argument to $1-X$ at null mean $1-\mu_0$ gives $\beta_{Q_k}^\star>-c$ eventually. Since $c$ is arbitrary, $|\beta_{Q_k}^\star|\to0$.

Equation~\eqref{eq:small-bet-asymptotics}, with reflection for negative bets, now gives
\[
\lstar_{Q_k}=\tfrac12(\beta_{Q_k}^\star)^2\E_{Q_k}[(X-\mu_0)^2](1+o(1))~.
\]
The second moment is bounded above by 1 and bounded away from zero, so
\[
\ln(1/\lstar_{Q_k})=2\ln(1/|\beta_{Q_k}^\star|)+O(1)~.
\]
Taking logarithms again proves the claimed relation between the two $\LL$ terms.
\end{proof}

\subsection{Proof of Theorem~\ref{thm:joint-minimax}}
\label{sec:proof-joint-minimax}

Recall from \eqref{eq:oracle-range} that $\cQ_R$ consists of the alternatives with $e^{-2R}\leq|\beta_Q^\star|\leq e^{-R}$, and $\mathcal T_\alpha$ contains all valid level-$\alpha$ stopping times. We first prove a finite lower bound whose error is uniform in the number and location of the betting scales.

\begin{lemma}[Confidence and scale information]
\label{lemma:scale-information}
Define
\[
\sigma_0^2:=\mu_0(1-\mu_0),\qquad
\rho_\alpha:=\alpha^{3/8}+2\alpha^{\sigma_0^2/64}~.
\]
Fix $0<\alpha\leq e^{-1}$, an integer $m\geq2$, and gaps $\Delta_1, \dots, \Delta_m$ satisfying
\[
0<\Delta_m<\cdots<\Delta_1\leq\sigma_0^2/e,\qquad
\Delta_{j+1}\leq\Delta_j/4\quad(1\leq j<m)~.
\]
For the Bernoulli laws
\[
P:=\mathrm{Ber}(\mu_0),\qquad Q_j:=\mathrm{Ber}(\mu_0+\Delta_j)
\quad(1\leq j\leq m),
\]
every stopping time with $\Pr_P\{\tau<\infty\}\leq\alpha$ satisfies
\begin{equation}
\label{eq:finite-scale-information}
\frac1m\sum_{j=1}^m\lstar_{Q_j}\E_{Q_j}[\tau]
\geq\ln(1/\alpha)+(1-\rho_\alpha)\ln m-\ln2~.
\end{equation}
In particular, the maximum of these normalized expected times satisfies the same lower bound. The result includes stopping times using independent randomization.
\end{lemma}

\begin{proof}
If any of the expectations on the left is infinite, the claim is immediate. We may therefore assume that they are all finite. This gives power one at the finitely many alternatives used in the proof.

For these Bernoulli laws, recall from \eqref{eq:bernoulli-asymptotics} that
\begin{equation}
\label{eq:bernoulli-oracle}
\beta_{Q_j}^\star=\frac{\Delta_j}{\sigma_0^2},\qquad
\lstar_{Q_j}=\KL(Q_j\|P)=\KLBern(\mu_0+\Delta_j;\mu_0),
\end{equation}
where $P$ attains the infimum divergence over the full bounded null class, as shown in the proof of Corollary~\ref{thm:main}.

Let $U$ be the independent random seed, with the same law under every hypothesis, where a constant seed covers deterministic tests. Define
\[
\mathscr F_n:=\sigma(U,X_1,\ldots,X_n)\qquad(n\geq0)~.
\]

\emph{Rejection is unlikely before enough observations have arrived.}
For each alternative, choose the deterministic time
\[
n_j:=\left\lfloor\frac{\sigma_0^2\ln(1/\alpha)}{4\Delta_j^2}\right\rfloor~.
\]
The quantity inside the floor is at least $e^2/(4\sigma_0^2)>2$. Hence, $n_j\geq1$, and rounding loses at most half of that quantity:
\begin{equation}
\label{eq:cutoff-size}
n_j\Delta_j^2\geq\frac{\sigma_0^2\ln(1/\alpha)}8~.
\end{equation}
The likelihood ratio on the first $n$ observations is
\[
L_{j,n}:=\prod_{i=1}^n\frac{\mathrm dQ_j}{\mathrm dP}(X_i)
=\prod_{i=1}^n\left(1+\frac{\Delta_j}{\sigma_0^2}(X_i-\mu_0)\right)~.
\]
Independence and the definition of $\sigma_0$ give
\[
\E_P[L_{j,n}^2]=\left(1+\frac{\Delta_j^2}{\sigma_0^2}\right)^n
\leq\exp\left(\frac{n\Delta_j^2}{\sigma_0^2}\right)~.
\]
The seed contributes no likelihood factor. Since $\{\tau\leq n_j\}\in\mathscr F_{n_j}$, change of measure and Cauchy--Schwarz yield
\begin{equation}
\label{eq:early-rejection}
\Pr_{Q_j}\{\tau\leq n_j\}
\leq\sqrt{\Pr_P\{\tau\leq n_j\}\E_P[L_{j,n_j}^2]}
\leq\sqrt{\alpha\exp(\ln(1/\alpha)/4)}
=\alpha^{3/8}~.
\end{equation}

\emph{The empirical mean identifies the scale after that time.}
Define the cumulative payoff by
\[
S_n:=\sum_{i=1}^n(X_i-\mu_0)~.
\]
Under $Q_j$, Hoeffding's lemma~\citep{Hoeffding63} implies that, for each $s\in\{-1,1\}$,
\[
M_n^{s,j}:=\exp\left(s\Delta_j(S_n-n\Delta_j)-\frac{n\Delta_j^2}{8}\right),
\qquad n\geq0,
\]
is a nonnegative supermartingale starting at one. If $s(S_n-n\Delta_j)\geq n\Delta_j/4$ at some $n\geq n_j$, then $M_n^{s,j}\geq\exp(n_j\Delta_j^2/8)$. Applying Ville's inequality~\citep{Ville39} separately for each $s$ gives
\[
\Pr_{Q_j}\{\exists n\geq n_j:s(S_n-n\Delta_j)\geq n\Delta_j/4\}
\leq e^{-n_j\Delta_j^2/8}
\leq\alpha^{\sigma_0^2/64},
\]
where the final inequality uses \eqref{eq:cutoff-size}. A union bound over the two signs now yields
\begin{equation}
\label{eq:late-concentration}
\Pr_{Q_j}\{\exists n\geq n_j:|S_n/n-\Delta_j|\geq\Delta_j/4\}
\leq2\alpha^{\sigma_0^2/64}~.
\end{equation}

The intervals $[3\Delta_j/4,5\Delta_j/4]$ are disjoint by the separation assumption. Define a decoder from the stopped observations by
\[
\widehat J:=
\begin{cases}
j,&\tau<\infty\text{ and }S_\tau/\tau\in[3\Delta_j/4,5\Delta_j/4],\\
1,&\text{otherwise}.
\end{cases}
\]
Under $Q_j$, the stopping time is finite almost surely. Combining \eqref{eq:early-rejection} and \eqref{eq:late-concentration}, we obtain
\begin{equation}
\label{eq:decoder-error}
\Pr_{Q_j}\{\widehat J\neq j\}\leq\rho_\alpha
\qquad(1\leq j\leq m)~.
\end{equation}

\emph{Rejection and scale identification require additive information.}
Record the seed and the finite stopped history, using a distinguished symbol for nonrejection:
\[
Z:=\begin{cases}(U,\tau,X_1,\ldots,X_\tau),&\tau<\infty,\\
\dagger,&\tau=\infty.
\end{cases}
\]
Define the laws of the stopped history by
\[
Q_j^\tau:=\operatorname{Law}_{Q_j}(Z),\qquad P^\tau:=\operatorname{Law}_P(Z)~.
\]
Their density ratio on finite histories is $L_{j,\tau}$. To compute its expected logarithm, first stop at $\tau\wedge N$. Independence gives
\[
\E_{Q_j}[\ln L_{j,\tau\wedge N}]
=\sum_{i=1}^N\E_{Q_j}\left[\indevent\{\tau\geq i\}\ln\frac{\mathrm dQ_j}{\mathrm dP}(X_i)\right]
=\KL(Q_j\|P)\E_{Q_j}[\tau\wedge N]~.
\]
Each log likelihood increment is bounded, and $\E_{Q_j}[\tau]<\infty$. Thus, dominated convergence applies to the stopped log likelihood, while monotone convergence applies to $\tau\wedge N$. Using \eqref{eq:bernoulli-oracle}, we conclude that
\begin{equation}
\label{eq:stopped-information}
\KL(Q_j^\tau\|P^\tau)=\lstar_{Q_j}\E_{Q_j}[\tau]~.
\end{equation}

For the information calculation, draw the generating index uniformly and then draw observations under $Q_J$:
\[
\Pr\{J=j\}:=1/m\quad(1\leq j\leq m),\qquad
\overline Q^\tau:=\frac1m\sum_{j=1}^m Q_j^\tau~.
\]
The mutual information between the index and the stopped history is
\[
\operatorname{I}(J;Z):=\frac1m\sum_{j=1}^m\KL(Q_j^\tau\|\overline Q^\tau)~.
\]
Splitting each log density ratio at the mixture density gives the KL chain rule
\begin{equation}
\label{eq:information-decomposition}
\frac1m\sum_{j=1}^m\KL(Q_j^\tau\|P^\tau)
=\KL(\overline Q^\tau\|P^\tau)+\operatorname{I}(J;Z)~.
\end{equation}
Since $\Pr_{\overline Q^\tau}\{Z\neq\dagger\}=1$ and $\Pr_{P^\tau}\{Z\neq\dagger\}\leq\alpha$, data processing through the rejection event gives
\[
\KL(\overline Q^\tau\|P^\tau)\geq\ln(1/\alpha)~.
\]

Fano's inequality~\citep{CoverT06} bounds the second term. We include its short argument. Define the decoding-error indicator by
\[
D:=\indevent\{\widehat J\neq J\}~.
\]
Conditional on $Z$, this indicator determines whether $J$ equals the decoded value; on an error there are at most $m-1$ remaining choices. Therefore, with entropy and mutual information measured in nats,
\begin{align*}
\entropy(J\mid Z)
&\leq\entropy(D\mid Z)+\entropy(J\mid D,Z)
\leq\ln2+\Pr\{D=1\}\ln m
\leq\ln2+\rho_\alpha\ln m,\\
\operatorname{I}(J;Z)
&=\ln m-\entropy(J\mid Z)
\geq(1-\rho_\alpha)\ln m-\ln2~.
\end{align*}
Substituting these two bounds and \eqref{eq:stopped-information} into \eqref{eq:information-decomposition} proves \eqref{eq:finite-scale-information}.
\end{proof}

We can now prove Theorem~\ref{thm:joint-minimax}.
\begin{proof}[Proof of Theorem~\ref{thm:joint-minimax}]
\emph{Lower bound.}
For $R\geq\ln4$, choose the number of alternatives and their mean gaps as
\[
m:=\left\lfloor\frac{R}{\ln4}\right\rfloor+1,\qquad
\Delta_j:=\mu_0(1-\mu_0)e^{-R}4^{-(j-1)},\qquad
Q_j:=\mathrm{Ber}(\mu_0+\Delta_j)\quad(1\leq j\leq m)~.
\]
The gaps satisfy Lemma~\ref{lemma:scale-information}, and \eqref{eq:bernoulli-oracle} gives
\[
\beta_{Q_j}^\star=e^{-R}4^{-(j-1)}\in[e^{-2R},e^{-R}]~.
\]
Thus, all these alternatives belong to $\cQ_R$. Every $\tau\in\mathcal T_\alpha$ is valid under their common null $P=\mathrm{Ber}(\mu_0)$. Lemma~\ref{lemma:scale-information} gives the finite bound
\begin{equation}
\label{eq:finite-minimax}
\inf_{\tau\in\mathcal T_\alpha} \sup_{Q\in\cQ_R} \ \frac{\lstar_Q\E_Q[\tau]}{\ln(1/\alpha)+\ln R}
\geq\frac{\ln(1/\alpha)+(1-\rho_\alpha)\ln m-\ln2} {\ln(1/\alpha)+\ln R}~.
\end{equation}
For $R\geq\ln4$, the inequalities $R/\ln4<m\leq2R/\ln4$ imply that $\ln m-\ln R$ is bounded by an absolute constant in magnitude. Also $\rho_\alpha\to0$ as $\alpha\downarrow0$ for fixed $\mu_0$. Hence, the right-hand side of \eqref{eq:finite-minimax} tends to 1 jointly, regardless of the relative rates of $\alpha$ and $R$. In particular, for some absolute $C<\infty$ and all sufficiently small $\alpha$, it is at least
\[
1-\rho_\alpha-\frac{C}{\ln(1/\alpha)+\ln R}~.
\]

\emph{Upper bound and simultaneous attainment.}
Fix a prior $F$ satisfying Definition~\ref{def:prior} with $a(F)=1$; Corollary~\ref{cor:examples}(iii) supplies such priors. Its test belongs to $\mathcal T_\alpha$ by \eqref{eq:portfolio-test}. Since both directional indices equal one, Theorem~\ref{thm:pointwise} gives, for every $\eta>0$, a cutoff $\alpha_0>0$ such that $\E_Q[\tau_\alpha^F]\leq(1+\eta)B_\alpha(Q)$ for every $Q\in\cQ$ and every $0<\alpha\leq\alpha_0$. Thus, for all these error levels and every $R\geq1$, \eqref{eq:range-adaptation} gives
\[
\sup_{Q\in\cQ_R} \ \frac{\lstar_Q\E_Q[\tau_\alpha^F]}{\ln(1/\alpha)+\ln R}
\leq (1+\eta)\left(1+\frac{\ln2}{\ln(1/\alpha)+\ln R}\right)~.
\]
Taking the joint upper limit and then sending $\eta\downarrow0$ proves the upper half of \eqref{eq:joint-minimax}. The finite lower bound applies to this fixed mixture as well, proving simultaneous attainment.

Finally, observe that
\begin{equation}
\label{eq:range-adaptation}
\ln R\leq\LL(|\beta_Q^\star|)\leq\ln R+\ln2\qquad(Q\in\cQ_R)~.
\end{equation}
Hence, \eqref{eq:range-adaptation} bounds the two normalizing denominators between $\ln(1/\alpha)+\ln R$ and that quantity plus $\ln2$. Their ratio tends to 1 uniformly over $\cQ_R$, proving the stated result.
\end{proof}

\section{Proofs of Section~\ref{sec:confidence-sequences}}
\label{app:confidence-sequences}

\subsection{Proof of Proposition~\ref{cs:prop:validity}}

\begin{proof}[Proof of Proposition~\ref{cs:prop:validity}]
At the true mean, each fixed-bet product in \eqref{cs:eq:wealth} is a nonnegative martingale under the conditional-mean null. The prior mixture starts at one, so Ville's inequality gives
\[
\Pr_Q\left\{\sup_{n\geq1}\Wealth_n(\mu_Q)\geq1/\alpha\right\}\leq\alpha~.
\]
This proves coverage for interior means. If the conditional mean is 0 or 1, boundedness forces every observation to equal that endpoint almost surely, and our convention gives coverage as well.

If the wealth is convex in $m$, every set $\{m\in(0,1):\Wealth_s(m)<1/\alpha\}$ is an interval. Their running intersection is a nested interval, possibly empty. Adding the stipulated endpoints preserves this property: on an all-zero path, every admissible bet gives $1-\beta m\leq1/(1-m)$, so $\Wealth_s(m)\leq(1-m)^{-s}$ and the intersection contains all sufficiently small positive $m$. Thus, adding 0 leaves an interval. Reflection proves the claim at 1.
\end{proof}

\subsection{Proof of Theorem~\ref{cs:thm:width}}

First we need a lemma.
\begin{lemma}
\label{cs:lem:balance}
Fix $Q\in\cP$ with $\mu_Q\in(0,1)$ and $0<\gamma<1$. The function $m\mapsto\widetilde B_\gamma(Q,m)$ is continuous and strictly increasing on $(0,\mu_Q)$, continuous and strictly decreasing on $(\mu_Q,1)$, and divergent as $m\to\mu_Q$ from either side.
\end{lemma}
\begin{proof}
For $m<\mu_Q$, express log growth in the normalized betting coordinate:
\[
\phi_m(t):=\E_Q[\ln(1-t+tX/m)]\qquad(0\leq t\leq1)~.
\]
Its unique maximizer is the fraction of the maximal positive bet,
\[
\theta_Q^\star(m):=m\beta_Q^\star(m)>0~,
\]
and its maximum is $\lstar_Q(m)$. For $m_1<m_2$,
\begin{equation}
\label{cs:eq:objective-continuity}
0\leq\phi_{m_1}(t)-\phi_{m_2}(t)\leq\ln(m_2/m_1)
\end{equation}
whenever these values are finite; at $t=1$ they are either both finite or both $-\infty$. This bound gives continuity of the optimal value. Using the optimizer at $m_2$ also gives strict decrease of that value, since $\Pr_Q\{X>0\}>0$.

Upper semicontinuity in $t$, compactness, uniqueness of the maximizer, and \eqref{cs:eq:objective-continuity} give continuity of $\theta_Q^\star(m)$, including at endpoint optima. Its monotonicity follows from decreasing differences: for $t<1$,
\[
\frac{\partial^2}{\partial m\,\partial t}\ln(1-t+tx/m)=-\frac{x}{[m+t(x-m)]^2}\leq0~.
\]
The difference $\phi_{m_1}(t)-\phi_{m_2}(t)$ is nondecreasing in $t$, also at a finite endpoint by limits. If $\theta_Q^\star(m_1)<\theta_Q^\star(m_2)$, the two strict optimality inequalities contradict this monotonicity. Hence, $\theta_Q^\star(m)$ is nonincreasing. Since $\beta_Q^\star(m)=\theta_Q^\star(m)/m>0$, the oracle bet is nonincreasing as well.

The coefficient in \eqref{cs:eq:benchmark} is the fixed index $a_+(F)$ on this side. Its numerator is therefore nondecreasing and its positive denominator strictly decreases, proving strict increase of the benchmark. As $m\uparrow\mu_Q$, \eqref{cs:eq:objective-continuity} and Jensen's inequality give $\lstar_Q(m)\to0$, so the benchmark diverges. Reflection proves the upper-side statements with the fixed index $a_-(F)$.
\end{proof}

\begin{proof}[Proof of Theorem~\ref{cs:thm:width}]
We first obtain the required rejection-time limits from Corollary~\ref{cor:one-sided}. Define the positive- and negative-bet contributions of the rescaled reference prior by
\[
A_n^\pm(m):=\int_{\{\beta:\,\pm\beta>0\}}\prod_{i=1}^n(1+\beta_m(\beta)(X_i-m))\,\mathrm dF(\beta)~.
\]
Their initial masses are
\[
q_+:=F((0,1/\mu_0]),\qquad q_-:=F([-1/(1-\mu_0),0))~.
\]
Compact-set positivity gives $q_+,q_->0$, and $q_++q_-=1$. For positive bets, each product factor is $1-\mu_0\beta+\mu_0\beta X_i/m$, which is nonincreasing in $m$. For negative bets, it is $1+(1-\mu_0)\beta-(1-\mu_0)\beta(1-X_i)/(1-m)$, which is nondecreasing. Hence, $A_n^+$ is nonincreasing and $A_n^-$ nondecreasing. Changing variables in the prior comparison \eqref{cs:eq:prior-comparison} gives
\begin{equation}
\label{cs:eq:wealth-comparison}
c_\kappa(A_n^+(m)+A_n^-(m))\leq\Wealth_n(m)
\leq C_\kappa(A_n^+(m)+A_n^-(m))
\qquad(m\in[\kappa,1-\kappa])~.
\end{equation}
Integrating the density bounds also gives $c_\kappa\leq1\leq C_\kappa$. Define the reference crossing times in the direction of the alternative by
\[
\sigma_\gamma(Q,m):=
\begin{cases}
\inf\{n\geq1:A_n^+(m)\geq1/\gamma\},&m<\mu_Q,\\
\inf\{n\geq1:A_n^-(m)\geq1/\gamma\},&m>\mu_Q.
\end{cases}
\]

For $m<\mu_Q$, rescale the observations and normalize the positive part of the prior at the fixed null mean $\kappa$:
\[
\widehat X_i:=\frac{\kappa X_i}{m},\qquad
f_\kappa^+(\beta):=\frac{\kappa}{\mu_0 q_+}f\left(\frac{\kappa\beta}{\mu_0}\right)
\quad(0<\beta<1/\kappa)~.
\]
The density integrates to 1. Since $m\geq\kappa$, the transformed observations lie in $[0,1]$ and have mean greater than $\kappa$. Rescaling the bet preserves each factor:
\[
1+\frac{\mu_0\beta}{\kappa}(\widehat X_i-\kappa)
=1+\beta_m(\beta)(X_i-m)\qquad(\beta>0)~.
\]
Thus, $A_n^+(m)/q_+$ is the exact one-sided mixture with density $f_\kappa^+$, and $\sigma_\gamma(Q,m)$ is its rejection time at level $q_+\gamma$. Its oracle growth is $\lstar_Q(m)$ and its oracle bet is $m\beta_Q^\star(m)/\kappa$. Scaling and normalization preserve the index $a_+(F)$ and compact-set positivity, so Corollary~\ref{cor:one-sided} applies uniformly over all transformed laws.

For $m>\mu_Q$, use the reflected observations and normalized negative prior:
\[
\widehat X_i:=\frac{\kappa(1-X_i)}{1-m},\qquad
f_\kappa^-(\beta):=\frac{\kappa}{(1-\mu_0)q_-}
f\left(-\frac{\kappa\beta}{1-\mu_0}\right)
\quad(0<\beta<1/\kappa)~.
\]
The same argument applies to $A_n^-(m)/q_-$ at level $q_-\gamma$, with index $a_-(F)$, growth $\lstar_Q(m)$, and oracle bet $(1-m)|\beta_Q^\star(m)|/\kappa$.

We also need fixed changes in the error level. For $x>0$ and $1\leq v\leq1/\kappa$,
\begin{equation}
\label{cs:eq:ll-shift}
0\leq\LL(x)-\LL(vx)\leq\ln(1+\ln(1/\kappa))~.
\end{equation}
Indeed, $\ln(1/x)=\ln(1/(vx))+\ln v$ bounds $\max\{1,\ln(1/x)\}$ by $(1+\ln(1/\kappa))\max\{1,\ln(1/(vx))\}$. For each fixed $d>0$, the positive-side testing benchmark for $\sigma_{d\alpha}(Q,m)$ is
\[
\frac{\ln(1/(q_+d\alpha))+a_+(F)\LL(m|\beta_Q^\star(m)|/\kappa)}{\lstar_Q(m)}~,
\]
with $q_-,a_-(F),1-m$ in place of $q_+,a_+(F),m$ on the negative side. By \eqref{cs:eq:ll-shift}, its absolute relative difference from \eqref{cs:eq:benchmark} is at most
\begin{equation}
\label{cs:eq:confidence-shift}
\frac{\max_{s\in\{+,-\}}|\ln(q_s d)|+a(F)\ln(1+\ln(1/\kappa))}{\ln(1/\alpha)}~.
\end{equation}
For $x\geq0$ and $r>0$, the inequality $|rx-1|\leq r|x-1|+|r-1|$ therefore transfers the testing limit to our benchmark. In particular,
\begin{equation}
\label{cs:eq:l1}
\rho_\kappa(\alpha):=
\sup_{\substack{d\in\{c_\kappa,2C_\kappa\},\ Q\in\cQ_\kappa\\
m\in[\kappa,1-\kappa],\ m\ne\mu_Q}}
\E_Q\left[\left|\frac{\sigma_{d\alpha}(Q,m)}{\widetilde B_\alpha(Q,m)}-1\right|\right]
\longrightarrow0~.
\end{equation}
We henceforth take $\alpha$ sufficiently small that $2C_\kappa\alpha<1$ and $\rho_\kappa(\alpha)<\infty$; all the corresponding times are integrable.

By Lemma~\ref{cs:lem:balance}, as $r$ increases from 0 to $w$, the lower-endpoint benchmark decreases continuously from infinity and the upper-endpoint benchmark increases continuously to infinity. They meet exactly once, giving $r^\star$. Define the balanced endpoints and their common benchmark by
\[
l:=\mu_Q-r^\star,\qquad u:=\mu_Q+w-r^\star,
\]
\begin{equation}
\label{cs:eq:balanced-benchmark}
K_\alpha(Q,w):=\widetilde B_\alpha(Q,l)=\widetilde B_\alpha(Q,u)~.
\end{equation}
Both endpoints lie in $(\kappa,1-\kappa)$, since $\mu_Q\in[2\kappa,1-2\kappa]$ and $w<\kappa$.

To bound the width on every path, we use convexity and the empirical mean
\[
\overline X_n:=\frac1n\sum_{i=1}^n X_i~.
\]
At an interior empirical mean, the arithmetic--geometric mean inequality gives
\[
\prod_{i=1}^n(1+\beta(X_i-\overline X_n))
\leq\left(\frac1n\sum_{i=1}^n(1+\beta(X_i-\overline X_n))\right)^n=1
\Rightarrow \Wealth_n(\overline X_n)\leq 1~.
\]
If $A_n^+(l)\geq1/(c_\kappa\alpha)$, then $\overline X_n>l$: otherwise the same inequality bounds each positive-bet product at $l$ by 1, so $A_n^+(l)\leq q_+\leq1$. At this crossing, \eqref{cs:eq:wealth-comparison} gives $\Wealth_n(l)\geq1/\alpha$. If $\overline X_n<1$, its wealth is at most 1. If $\overline X_n=1$, all observations are one and $\Wealth_n(m)\leq m^{-n}$, which tends to 1 as $m\uparrow1$. In either case, there is a point $z\in(l,1)$ with $\Wealth_n(z)<1/\alpha$. For any $0<m<l$, convexity gives
\[
\Wealth_n(l)\leq\frac{z-l}{z-m}\Wealth_n(m)
+\frac{l-m}{z-m}\Wealth_n(z),
\]
which forces $\Wealth_n(m)>1/\alpha$. Thus, every $m\leq l$ is excluded from the running inversion, including 0 because $\overline X_n>l$. Reflection proves that crossing by $A_n^-(u)$ excludes every $m\geq u$. Consequently, on every path,
\begin{equation}
\label{cs:eq:pathwise}
T_w^\alpha\leq\max\{\sigma_{c_\kappa\alpha}(Q,l),\sigma_{c_\kappa\alpha}(Q,u)\}~.
\end{equation}
Equation~\eqref{cs:eq:l1} proves integrability. Using $\max\{x,y\}\leq1+|x-1|+|y-1|$ also gives
\[
\frac{\E_Q[T_w^\alpha]}{K_\alpha(Q,w)}\leq1+2\rho_\kappa(\alpha)~.
\]

For the lower bound, Ville's inequality applied to $\Wealth_n(\mu_Q)$ gives
\[
\mathcal E_\alpha:=\left\{\sup_{n\geq1}\Wealth_n(\mu_Q)<\frac{c_\kappa}{2C_\kappa\alpha}\right\},
\qquad \Pr_Q\{\mathcal E_\alpha^c\}\leq\frac{2C_\kappa\alpha}{c_\kappa}~.
\]
By continuity and monotonicity of the benchmarks, choose $l'<l$ and $u'>u$, still in $(\kappa,1-\kappa)$, with
\[
1-\alpha\leq\frac{\widetilde B_\alpha(Q,m)}{K_\alpha(Q,w)}\leq1
\qquad(m\in\{l',u'\})~.
\]
On $\mathcal E_\alpha$, monotonicity of the reference contributions and \eqref{cs:eq:wealth-comparison} imply, for every $n$,
\[
C_\kappa A_n^-(l')\leq C_\kappa A_n^-(\mu_Q)
\leq\frac{C_\kappa}{c_\kappa}\Wealth_n(\mu_Q)<\frac1{2\alpha}~.
\]
Before $A_n^+(l')$ reaches $1/(2C_\kappa\alpha)$, the upper bound in \eqref{cs:eq:wealth-comparison} therefore keeps $\Wealth_n(l')<1/\alpha$. Reflection gives the same conclusion at $u'$. Both endpoints remain in the running inversion until the first of these crossings, and $u'-l'>w$. Hence
\[
T_w^\alpha\geq\min\{\sigma_{2C_\kappa\alpha}(Q,l'),\sigma_{2C_\kappa\alpha}(Q,u')\}
\qquad\text{on }\mathcal E_\alpha~.
\]
Define the bounded lower comparison
\[
U_{\alpha,Q,w}:=\indevent\{\mathcal E_\alpha\}
\min\left\{1,\frac{\sigma_{2C_\kappa\alpha}(Q,l')}{K_\alpha(Q,w)},
\frac{\sigma_{2C_\kappa\alpha}(Q,u')}{K_\alpha(Q,w)}\right\}~.
\]
It satisfies $0\leq U_{\alpha,Q,w}\leq T_w^\alpha/K_\alpha(Q,w)$ on every path. For $x\geq0$ and $0\leq r\leq1$, $|rx-1|\leq|x-1|+1-r$. Apply this with each endpoint's benchmark ratio, whose deficit from 1 is at most $\alpha$. Since $1-\min\{1,x,y\}\leq|x-1|+|y-1|$, we obtain
\[
\E_Q[|U_{\alpha,Q,w}-1|]\leq
\frac{2C_\kappa\alpha}{c_\kappa}+2\rho_\kappa(\alpha)+2\alpha\longrightarrow0
\]
uniformly over $Q,w$. The expected upper bound and this lower comparison verify Lemma~\ref{lemma:l1-criterion} and prove \eqref{cs:eq:uniform-width}.

Finally, if $F_m$ is the distribution of $\beta_m(\beta)$ under $F$, then $\Wealth_n(m)=A_n^+(m)+A_n^-(m)$. In the product factors displayed at the start of the proof, the constant term and the coefficient of $1/m$ or $1/(1-m)$ are nonnegative. Expanding the positive-bet product therefore gives a nonnegative linear combination of $m^{-k}$, $0\leq k\leq n$, and the negative-bet product gives one of $(1-m)^{-k}$. These functions are convex; integration and addition preserve convexity.
\end{proof}

\subsection{Verification of Remark~\ref{cs:rem:oj}}
\label{app:cs-oj}

The integrals in \citet[equations~(5), (11), and Remark~2]{OrabonaJ21} are
\begin{align}
\label{cs:eq:oj-wealth}
\Wealth_n(m)
=\frac12\int_0^1\prod_{i=1}^n\left(p\frac{X_i}{m}+(1-p)\frac{1-X_i}{1-m}\right)\frac{\mathrm dp}{\pi\sqrt{p(1-p)}}
+\frac12\int_{-1}^1\prod_{i=1}^n(1+\beta(X_i-m))f_{\mathrm{R70}}(\beta)\,\mathrm d\beta~.
\end{align}

To check the prior comparison, use the fraction of the maximal directional bet as a coordinate. For $0<t<1$, define the corresponding density by
\[
h_m(t):=\frac{f_m(t/m)}m
=\frac{1}{2\pi\sqrt{(1-t)(m/(1-m)+t)}}
+\frac{f_{\mathrm{R70}}(t/m)}{2m}\indevent\{t\leq m\},
\]
\[
h_m(-t):=\frac{f_m(-t/(1-m))}{1-m}=h_{1-m}(t)~.
\]

As a reference prior satisfying Definition~\ref{def:prior}, take $\mu_0=1/2$ and the density
\[
f(\beta):=\frac14 f_{\mathrm{R70}}(\beta/2)
+\frac1{16\sqrt{1-|\beta|/2}}\qquad(0<|\beta|<2)~.
\]
Each term integrates to $1/2$ on $(-2,2)$; the density is bounded below on every compact subset away from zero and has both indices equal to 1. The rescaling \eqref{cs:eq:prior-rescaling} gives
\[
f_m(\beta_m(\beta))\beta_m'(\beta)=\tfrac12 h_m(\beta/2)~.
\]
Thus, it suffices to compare $h_m(t)$ with $2f(2t)$. Uniformly for $m\in[\kappa,1-\kappa]$, the CO96 term of $h_m(t)$ is bounded above and below by positive constant multiples of $(1-t)^{-1/2}$: its other denominator factor satisfies
\[
\frac{\kappa}{1-\kappa}\leq\frac m{1-m}+t\leq\frac1\kappa~.
\]
For $t\leq m$, the R70 density ratio is
\[
\frac{f_{\mathrm{R70}}(t/m)/m}{f_{\mathrm{R70}}(t)}
=\frac{\ln(c/t)}{\ln(cm/t)}
\left(\frac{\ln\ln(c/t)}{\ln\ln(cm/t)}\right)^2~.
\]
Both ratios are at least 1 and uniformly bounded. Indeed, $\ln(c/t)=\ln(cm/t)+\ln(1/m)$, $\ln(cm/t)\geq\ln c>1$, and $0\leq\ln(1/m)\leq\ln(1/\kappa)$ give
\[
\frac{\ln(c/t)}{\ln(cm/t)}\leq1+\frac{\ln(1/\kappa)}{\ln c},\qquad
\frac{\ln\ln(c/t)}{\ln\ln(cm/t)}
\leq1+\frac{\ln(1+\ln(1/\kappa)/\ln c)}{\ln\ln c}~.
\]
For $t>m\geq\kappa$, the omitted R70 term is bounded on $[\kappa,1]$ and is therefore dominated by a constant multiple of the CO96 term. These bounds and reflection prove \eqref{cs:eq:prior-comparison}, with constants depending only on $\kappa$. Fixed positive mixing weights only change these constants.

Finally, \citet[Theorem~10]{OrabonaJ21} gives convexity in $m$ for the R70 integral, whose betting prior is fixed on $[-1,1]$. For CO96, expand the product at each fixed $p$: it is a nonnegative linear combination of $m^{-k}(1-m)^{-(n-k)}$, $0\leq k\leq n$. These terms are log-convex, hence convex, because
\[
\frac{\mathrm d^2}{\mathrm dm^2}\ln\bigl(m^{-k}(1-m)^{-(n-k)}\bigr)
=\frac{k}{m^2}+\frac{n-k}{(1-m)^2}\geq0~.
\]
Integration and addition preserve convexity. Thus, Theorem~\ref{cs:thm:width} applies with $a(F)=1$, and Proposition~\ref{cs:prop:validity} gives nested confidence intervals. Their endpoint conventions also agree with the wealth limits: at 0 the CO96 integral is at most 1 on an all-zero path, and otherwise diverges by Fatou's lemma, since each fixed $0<p<1$ gives a diverging product. The R70 integral tends to 1 on an all-zero path. Reflection gives the corresponding limits at 1.

\subsection{Proofs of simplifying corollaries}

\begin{proof}[Proof of Corollary~\ref{cs:cor:local}]
The derivative argument in the proof of Proposition~\ref{prop:near-null-scales}, with the candidate mean varying, gives $\beta_Q^\star(\mu_Q\pm r)\to0$ as $r\downarrow0$. Since $\E_Q[(X-(\mu_Q\pm r))^2]=\sigma_Q^2+r^2$, the small-bet identities \eqref{eq:small-bet-asymptotics} give
\begin{equation}
\label{cs:eq:local-oracles}
\beta_Q^\star(\mu_Q\pm r)=\mp\frac r{\sigma_Q^2}(1+o(1)),\qquad
\lstar_Q(\mu_Q\pm r)=\frac{r^2}{2\sigma_Q^2}(1+o(1))~.
\end{equation}
Consequently, $\LL(|\beta_Q^\star(\mu_Q\pm r)|)=\LL(r)+o(1)$. For $r=cw$, uniformly over $c$ in each compact subinterval of $(0,1)$ and in the joint limit,
\[
\widetilde B_\alpha(Q,\mu_Q-cw)
=\frac{2\sigma_Q^2}{c^2w^2}\{\ln(1/\alpha)+a_+(F)\LL(w)\}(1+o(1))~,
\]
with $a_-(F)$ in the reflected formula for the upper endpoint. Recall that $K_\alpha(Q,w)$ in \eqref{cs:eq:balanced-benchmark} is the common cost at the balanced endpoints. For any interval of width $w$ containing $\mu_Q$, at least one endpoint cost is as large as this common cost: moving the lower endpoint to the left decreases its cost and increases the upper endpoint's cost, and conversely. In particular,
\[
K_\alpha(Q,w)\leq
\max\{\widetilde B_\alpha(Q,\mu_Q-w/2),\widetilde B_\alpha(Q,\mu_Q+w/2)\}
\leq(1+o(1))\frac{8\sigma_Q^2}{w^2}\{\ln(1/\alpha)+a(F)\LL(w)\}~.
\]
Theorem~\ref{cs:thm:width} gives $\E_Q[T_w^\alpha]\sim K_\alpha(Q,w)$, proving \eqref{cs:eq:local-time}. If $a_+(F)=a_-(F)$, the expansions at $r=w/4$ and $r=3w/4$ place $r^\star/w$ between these two fractions eventually. The expansions are uniform there, and equality of the endpoint costs forces $r^\star/w\to1/2$. Thus, \eqref{cs:eq:local-time} is an asymptotic equality in this case.

To invert the expected-time bound, define
\[
v_n:=\sqrt{\frac{8\sigma_Q^2}{n}\{\ln(1/\alpha_n)+a(F)\ln\ln n\}}~.
\]
The assumptions imply $v_n\to0$. For every fixed $c>0$,
\[
\frac{\ln(1/\alpha_n)+a(F)\LL(cv_n)}{\ln(1/\alpha_n)+a(F)\ln\ln n}\longrightarrow1~.
\]
Indeed, when $\ln(1/\alpha_n)/\ln\ln n$ is bounded, $\ln(1/v_n)\sim\tfrac12\ln n$; when that ratio diverges, both iterated-logarithm terms are negligible relative to $\ln(1/\alpha_n)$. These two cases cover all subsequences. Equation~\eqref{cs:eq:local-time} gives
\[
\limsup_n\frac{\E_Q[T_{cv_n}^{\alpha_n}]}n\leq c^{-2}~.
\]
For $c>1$, the expected time is therefore below $n$ eventually, so $w_{\mathrm e}(n,Q,\alpha_n)\leq cv_n$. Sending $c\downarrow1$ proves \eqref{cs:eq:local-effective}. If the indices coincide, the expected-time ratio tends to $c^{-2}$. Taking $c<1$ then gives the matching lower bound by monotonicity in the target width.

For fixed $\alpha\in(0,1/2)$, reduce the error probability slowly with the target width:
\[
\alpha_w:=\min\{\alpha,\exp(-\sqrt{\LL(w)})\}~.
\]
Then $\alpha_w\to0$ and $\ln(1/\alpha_w)=o(\LL(w))$. Since $C_n^\alpha\subseteq C_n^{\alpha_w}$, we have $T_w^\alpha\leq T_w^{\alpha_w}$ pathwise. Applying \eqref{cs:eq:local-time} at $(\alpha_w,w)$ gives the first upper bound in \eqref{cs:eq:fixed-alpha}. For the second, set
\[
v_n:=\sqrt{8a(F)\sigma_Q^2\ln\ln n/n}~.
\]
For any fixed $c>1$, the first bound gives $\limsup_n\E_Q[T_{cv_n}^\alpha]/n\leq c^{-2}<1$, so $w_{\mathrm e}(n,Q,\alpha)\leq cv_n$ eventually. Send $c\downarrow1$. These fixed-$\alpha$ conclusions are upper bounds.
\end{proof}

\begin{proof}[Proof of Corollary~\ref{cs:cor:deterministic}]
For $Q=\delta_\mu$, the oracle bet is $1/m$ below $\mu$ and $-1/(1-m)$ above $\mu$. Its magnitude exceeds 1, so its $\LL$ term is zero on both sides, regardless of the prior indices. The two rates are
\[
\lstar_Q(m)=\ln(\mu/m)\quad(m<\mu),\qquad
\lstar_Q(m)=\ln((1-\mu)/(1-m))\quad(m>\mu)~.
\]
At $m=\mu(1-w)$ and $m=\mu+(1-\mu)w$, both rates equal $-\ln(1-w)$. Equation~\eqref{cs:eq:balance} gives the claimed benchmark. Theorem~\ref{cs:thm:width} and $-\ln(1-w)\sim w$ give the expected time. Inverting at widths $c\ln(1/\alpha_n)/n$, with $c$ on either side of 1, gives the effective-width limit as in the preceding proof.
\end{proof}

\section{An efficient regret bound for the CO96/R70 mixture}
\label{app:oj-mixture}

Here, we derive a wealth lower bound for the exact prior mixture in \citet[Remark~2]{OrabonaJ21}. The resulting test retains the uniform $L^1$ guarantee of Theorem~\ref{thm:pointwise} and hence efficiency 1. It is a separate test from exact mixture inversion and from the PRECiSE bounds in their Algorithms~1 and~3.

The equal allocation in Remark~2 corresponds to $w=1/2$. For a regret bound, we need a density lower bound along the segment from zero to the hindsight-optimal bet. The R70 component supplies this bound near zero, and the CO96 component covers the rest of the admissible interval.

\begin{theorem}[A regret bound preserving efficiency]
\label{thm:oj-efficient-regret}
Let $\mu_0,w,F_w,f_w,f_{\mathrm{R70}}$ be as in Proposition~\ref{prop:oj-exact-mixture}. Define
\begin{equation}
\label{eq:oj-density-minorant}
\rho_w:=\min\left\{1-w,\frac{2w\mu_0(1-\mu_0)}{\pi f_{\mathrm{R70}}(1)}\right\},\qquad
\underline f_w(\beta):=\rho_w f_{\mathrm{R70}}(\min\{|\beta|,1\}),
\end{equation}
for $\beta\neq0$. For any sequence $X_1,X_2,\ldots\in[0,1]$ and $n\geq1$, choose a hindsight-optimal bet and its wealth:
\[
\widehat\beta_n\in\argmax_{\beta\in[-1/(1-\mu_0),1/\mu_0]}\Wealth_n(\beta),\qquad
\Wealth_n^\star:=\Wealth_n(\widehat\beta_n)\,.
\]
Take $\widehat\beta_n:=0$ whenever zero is a maximizer. Define the logarithmic-mean regret bound by
\begin{equation}
\label{eq:oj-lm-regret}
\Regret_n^{\mathrm{LM}}:=
\ln\frac{\ln\Wealth_n^\star}{1-1/\Wealth_n^\star}
-\ln\bigl(|\widehat\beta_n|\underline f_w(\widehat\beta_n)\bigr)
\end{equation}
when $\widehat\beta_n\neq0$, and set $\Regret_n^{\mathrm{LM}}:=+\infty$ otherwise. Combining this with the regret bound of universal portfolio for 2 stocks~\citep{CoverO96} gives
\begin{equation}
\label{eq:oj-efficient-regret}
\widehat{\Regret}_n:=\min\left\{
\ln\frac{\sqrt\pi\,\Gamma(n+1)}{w\,\Gamma(n+1/2)},\ \Regret_n^{\mathrm{LM}}
\right\},\qquad
\underline\Wealth_n:=\Wealth_n^\star\exp(-\widehat{\Regret}_n)\,.
\end{equation}

Then, for every $n\geq1$,
\[
\ln\Wealth_n^\star-\ln\Wealth_n^{F_w}
\leq\widehat{\Regret}_n\,.
\]
Consequently, $\underline\tau_\alpha:=\inf\{n\geq1:\underline\Wealth_n\geq1/\alpha\}$ defines a level-$\alpha$ family. For i.i.d. observations, its normalized rejection time converges uniformly in $L^1$:
\begin{equation}
\label{eq:oj-pointwise-proxy}
\lim_{\alpha\downarrow0}\sup_{Q\in\cQ}
\E_Q\left[\left|\frac{\underline\tau_\alpha}{B_\alpha(Q)}-1\right|\right]=0~.
\end{equation}
In particular, $\E_Q[\underline\tau_\alpha]/B_\alpha(Q)\to1$ uniformly over $Q\in\cQ$, and $\eff(\underline\tau)=1$.
\end{theorem}
\begin{proof}
On the admissible interval, $0<\underline f_w\leq f_w$ almost everywhere by \eqref{eq:oj-kt-density} and the definition of $\rho_w$. The $f_{\mathrm{R70}}$ density is decreasing in $|\beta|$ on $(0,1]$, as used in \citet[Section~V-H]{OrabonaJ21}. Hence, $\underline f_w(-\beta)=\underline f_w(\beta)$ and $\underline f_w$ is nonincreasing in $|\beta|$, while $b\underline f_w(b)$ is nondecreasing for $b>0$. Near zero, $\underline f_w$ has the same logarithmic decay index as $f_{\mathrm{R70}}$.

We first repeat the log-concavity argument in \citet[Section~V-H, equation~(16)]{OrabonaJ21}, using $\underline f_w$ along the segment from zero to $\widehat\beta_n$. If $\widehat\beta_n\neq0$, then $\Wealth_n^\star>1$, and
\begin{align*}
\Wealth_n^{F_w}
&\geq |\widehat\beta_n|\underline f_w(\widehat\beta_n)
\int_0^1\Wealth_n(t\widehat\beta_n)\,\mathrm{d}t\\
&\geq |\widehat\beta_n|\underline f_w(\widehat\beta_n)
\int_0^1(\Wealth_n^\star)^t\,\mathrm{d}t
=|\widehat\beta_n|\underline f_w(\widehat\beta_n)
\frac{\Wealth_n^\star-1}{\ln\Wealth_n^\star}\,.
\end{align*}
Taking logarithms gives \eqref{eq:oj-lm-regret}. We retain the logarithmic mean, rather than using its square-root lower bound in equation~(17) of the cited proof. The other branch in \eqref{eq:oj-efficient-regret} follows from $\Wealth_n^{F_w}\geq w\Wealth_n^{F_{\mathrm{CO96}}}$ and the regret bound of \citet{CoverO96}. Thus, $\underline\Wealth_n\leq\Wealth_n^{F_w}$, and Ville's inequality proves validity.

For the expected-time upper bound, fix $Q\in\cQ$, $\delta\in(0,1/2)$, $\eta\in(0,1)$, and $\epsilon>0$. Recall that $\beta_Q^\star$ maximizes expected log wealth, with value $\lstar_Q$. Set $\gamma_Q:=(1-\delta)\beta_Q^\star$ and $g_n(\beta):=\ln\Wealth_n(\beta)$. The increments of $g_n(\gamma_Q)$ are i.i.d., lie in $[\ln\delta,\ln\max\{1/\mu_0,1/(1-\mu_0)\}]$, and have mean
\[
\ell_Q(\gamma_Q)\geq(1-\delta)\lstar_Q>0
\]
by concavity. Following the Taylor argument in \citet[Section~V-H, equation~(19)]{OrabonaJ21}, we compare the empirical maximizer with $\gamma_Q$. We claim that whenever $g_n(\gamma_Q)>0$,
\begin{equation}
\label{eq:oj-empirical-bet}
\widehat\beta_n\gamma_Q>0,\qquad
|\widehat\beta_n|\geq c_\delta|\gamma_Q|,\qquad
c_\delta:=\frac{\delta^2}{2(\max\{1/\mu_0,1/(1-\mu_0)\})^2}\,.
\end{equation}
To prove this, write $S_n:=\sum_{i=1}^n(X_i-\mu_0)$ and $V_n:=\sum_{i=1}^n(X_i-\mu_0)^2$. On the segment between zero and $\gamma_Q$,
\[
-\frac{V_n}{\delta^2}\leq g_n''(\beta)
\leq-\frac{V_n}{(\max\{1/\mu_0,1/(1-\mu_0)\})^2}\,.
\]
Suppose $\gamma_Q>0$; the negative case follows analogously. Taylor's theorem and $g_n(\gamma_Q)>0$ give $V_n>0$ and
\[
S_n>\frac{\gamma_Q V_n}{2(\max\{1/\mu_0,1/(1-\mu_0)\})^2}>0\,.
\]
Hence, $\widehat\beta_n>0$. If $\widehat\beta_n\geq\gamma_Q$, the claim is immediate. Otherwise, $g_n'(\widehat\beta_n)=0$, and integrating the lower bound on $g_n''$ gives $S_n\leq\widehat\beta_n V_n/\delta^2$, proving \eqref{eq:oj-empirical-bet}. This argument also allows a population or empirical maximizer at an admissible endpoint.

By \eqref{eq:oj-empirical-bet}, monotonicity of $b\underline f_w(b)$, and $\ln\Wealth_n^\star\geq g_n(\gamma_Q)$, the logarithmic-mean branch gives, whenever $g_n(\gamma_Q)>0$,
\begin{align*}
\underline\Wealth_n
&\geq c_\delta|\gamma_Q|\underline f_w(c_\delta\gamma_Q)
\int_{1-\eta}^1\exp\bigl(tg_n(\gamma_Q)\bigr)\,\mathrm{d}t\\
&\geq\eta c_\delta|\gamma_Q|\underline f_w(c_\delta\gamma_Q)
\exp\bigl((1-\eta)g_n(\gamma_Q)\bigr)\,.
\end{align*}
The index-one behavior of $\underline f_w$ implies that, for some constant $C_{\delta,\epsilon}$ uniform over $Q$,
\begin{equation}
\label{eq:oj-density-penalty}
\ln\frac{1}{c_\delta|\gamma_Q|\underline f_w(c_\delta\gamma_Q)}
\leq(1+\epsilon)\LL(|\beta_Q^\star|)+C_{\delta,\epsilon}\,.
\end{equation}
Indeed, near zero, $\ln(1/(b\underline f_w(b)))/\ln\ln(1/b)\to1$. Also $c_\delta|\gamma_Q|=c_\delta(1-\delta)|\beta_Q^\star|$, and multiplication by this fixed positive constant changes $\LL(|\beta_Q^\star|)$ by a bounded amount near zero. Away from zero, $b\underline f_w(b)$ has a positive lower bound. These observations give \eqref{eq:oj-density-penalty} with the regularized $\LL$ used in the benchmark.

For any admissible nonzero $\beta$, monotonicity and the density minorant give $|\beta|\underline f_w(\beta)\leq1$ by integration over the segment from zero to $\beta$. Thus, the threshold $\ln(1/\alpha)+\ln(1/(\eta c_\delta|\gamma_Q|\underline f_w(c_\delta\gamma_Q)))$ is positive. Apply the bounded-stopping argument of Lemma~\ref{lemma:window} to $(1-\eta)g_n(\gamma_Q)$ at this threshold. Its mean increment is at least $(1-\eta)(1-\delta)\lstar_Q$, and \eqref{eq:oj-density-penalty} gives
\begin{equation}
\label{eq:oj-pointwise-regret-time}
\E_Q[\underline\tau_\alpha]\leq
\frac{\ln(1/\alpha)+(1+\epsilon)\LL(|\beta_Q^\star|)+C_{\mu_0,w,\delta,\eta,\epsilon}}
{(1-\delta)(1-\eta)\lstar_Q}\,,
\end{equation}
with the bounded overshoot absorbed into the constant, which we take to be nonnegative. Dividing by $B_\alpha(Q)$ gives
\begin{align*}
\frac{\E_Q[\underline\tau_\alpha]}{B_\alpha(Q)}
&\leq\frac{1}{(1-\delta)(1-\eta)}
\left(1+\frac{\epsilon\LL(|\beta_Q^\star|)+C_{\mu_0,w,\delta,\eta,\epsilon}}
{\ln(1/\alpha)+\LL(|\beta_Q^\star|)}\right)\\
&\leq\frac{1+\epsilon+C_{\mu_0,w,\delta,\eta,\epsilon}/\ln(1/\alpha)}{(1-\delta)(1-\eta)}~.
\end{align*}
The constant is finite and independent of $Q,\alpha$ for fixed $\delta,\eta,\epsilon$. First take the supremum over $Q$ and the upper limit as $\alpha\downarrow0$, and then send $\epsilon,\delta,\eta\downarrow0$. We obtain
\[
\limsup_{\alpha\downarrow0}\sup_{Q\in\cQ}\frac{\E_Q[\underline\tau_\alpha]}{B_\alpha(Q)}\leq1~.
\]

We now apply the comparison form of Lemma~\ref{lemma:l1-criterion} to the normalized proxy and exact-mixture rejection times:
\[
Z_{\alpha,Q}:=\frac{\underline\tau_\alpha}{B_\alpha(Q)},\qquad
U_{\alpha,Q}:=\frac{\tau_\alpha^{F_w}}{B_\alpha(Q)}~.
\]
Since $\underline\Wealth_n\leq\Wealth_n^{F_w}$, we have $\tau_\alpha^{F_w}\leq\underline\tau_\alpha$ pathwise, and hence $0\leq U_{\alpha,Q}\leq Z_{\alpha,Q}$. Both are integrable by \eqref{eq:oj-pointwise-regret-time}, and the preceding bound gives the required uniform upper limit of $\E_Q[Z_{\alpha,Q}]$. Proposition~\ref{prop:oj-exact-mixture} gives $a_Q(F_w)=1$, so its prior-dependent benchmark equals $B_\alpha(Q)$. Theorem~\ref{thm:pointwise} therefore gives uniform $L^1$ convergence of $U_{\alpha,Q}$ to 1. The lemma now proves \eqref{eq:oj-pointwise-proxy}.

Finally, Jensen's inequality also gives
\[
\left|\frac{\E_Q[\underline\tau_\alpha]}{B_\alpha(Q)}-1\right|
\leq\E_Q\left[\left|\frac{\underline\tau_\alpha}{B_\alpha(Q)}-1\right|\right]~.
\]
Thus, the normalized expectations converge uniformly to one, and Definition~\ref{def:eff} gives $\eff(\underline\tau)=1$.
\end{proof}

\section{A regret bound that does not preserve efficiency}
\label{app:ar-counterexample}

The mixture of \citet{AgrawalR26} has efficiency 1 by our index criterion. Unfortunately, its explicit regret bound can nevertheless lose efficiency when used to lower bound the wealth in a test. We state their bounds at null mean $1/2$ in our notation, including every branch, before giving the counterexample.

\begin{lemma}[Agrawal--Ramdas regret bounds at null mean $1/2$]
\label{lemma:ar-regret-bounds}
Fix $\mu_0:=1/2$ and $c:=6.6e$. Let $F_{\mathrm{AR26}}$ be the prior on $[-2,2]$ with density
\begin{equation}
\label{eq:ar-prior-density}
f_{\mathrm{AR26}}(\beta):=\frac{\ln\ln c}{2|\beta|\ln(2c/|\beta|)[\ln\ln(2c/|\beta|)]^2}
\qquad (0<|\beta|\leq2),
\end{equation}
with value zero outside $[-2,2]$ and arbitrary value at zero. For observations $X_1,X_2,\ldots\in[0,1]$ and $n\geq1$, define
\begin{align*}
\Wealth_n(\beta)&:=\prod_{i=1}^n(1+\beta(X_i-1/2)),\qquad
\Wealth_n^{F_{\mathrm{AR26}}}:=\int_{-2}^2\Wealth_n(\beta)f_{\mathrm{AR26}}(\beta)\,\mathrm{d}\beta,\\
S_n&:=\sum_{i=1}^n(X_i-1/2),\qquad V_n:=\sum_{i=1}^n(X_i-1/2)^2\,.
\end{align*}
Choose $\widehat\beta_n\in\argmax_{\beta\in[-2,2]}\Wealth_n(\beta)$, taking $\widehat\beta_n:=0$ when zero is a maximizer, and define the hindsight-optimal wealth and the actual regret by
\[
\Wealth_n^\star:=\Wealth_n(\widehat\beta_n),\qquad
\Regret_n^{\mathrm{AR26}}:=\ln\Wealth_n^\star-\ln\Wealth_n^{F_{\mathrm{AR26}}}\,.
\]
Define the constants
\[
C:=\ln(4\ln c\,\ln\ln c)+2\ln11,\qquad
C':=1+\ln(16\ln c\,\ln\ln c)\,.
\]
Set $\widehat{\Regret}^{(8)}_n:=0$ if $V_n=0$. For $V_n>0$, define
\begin{equation}
\label{eq:ar-pathwise-regret}
\widehat{\Regret}^{(8)}_n:=
\begin{cases}
\begin{aligned}
&9+\ln\frac{8}{\ln\ln c}+\ln\ln\bigl(14e\sqrt{1+V_n}\bigr)\\
&\quad+2\ln\ln\ln\bigl(14e\sqrt{1+V_n}\bigr),
\end{aligned}
& |S_n|<\sqrt{2V_n},\ |\widehat\beta_n|<2,\\[6pt]
4+\ln(4\ln c\,\ln\ln c),
& |S_n|<\sqrt{2V_n},\ |\widehat\beta_n|=2,\\[6pt]
\begin{aligned}
&1+\ln\frac{4}{\ln\ln c}+\ln\frac{20|S_n|}{3\sqrt{\tfrac43|S_n|+2V_n}}\\
&\quad+\ln\ln\bigl(28e(1+\sqrt{V_n})\bigr)\\
&\quad+2\ln\ln\ln\bigl(28e(1+\sqrt{V_n})\bigr),
\end{aligned}
& \sqrt{2V_n}\leq|S_n|\leq V_n/10,\\[6pt]
\tfrac12\ln\Wealth_n^\star+C,
& \sqrt{2V_n}\leq|S_n|,\ V_n/10<|S_n|\,.
\end{cases}
\end{equation}
Then, for every observation sequence and every $n\geq1$,
\begin{equation}
\label{eq:ar-two-regret-bounds}
\Regret_n^{\mathrm{AR26}}\leq
\min\{\widehat{\Regret}^{(8)}_n,\ln(n+1)+C'\}\,.
\end{equation}
\end{lemma}
\begin{proof}
These are the prior in equation~(7), the four branches of Theorem~4.1, equation~(8), and the logarithmic bound at the end of Remark~1 of \citet{AgrawalR26}, specialized to $\mu_0=1/2$. Their betting fraction is $-\beta$ in our notation, and $f_{\mathrm{AR26}}(-\beta)=f_{\mathrm{AR26}}(\beta)$ at this null mean. The quantities $\min\{\mu_0,1-\mu_0\}$ and $\max\{\mu_0,1-\mu_0\}$ both equal $1/2$, and $f_{\mathrm{AR26}}(2)=f_{\mathrm{AR26}}(-2)=1/(4\ln c\,\ln\ln c)$, giving the displayed constants. If $V_n=0$, all observations equal $1/2$, both wealths equal one, and the regret is zero.
\end{proof}

The last branch of \eqref{eq:ar-pathwise-regret} loses half of the best log wealth. Taking the minimum with the logarithmic bound in \eqref{eq:ar-two-regret-bounds} does not repair the joint efficiency, as the following example shows.

\begin{proposition}[Counterexample for the Agrawal--Ramdas regret bound]
\label{prop:ar-regret-counterexample}
Use the prior, wealths, and regret bounds of Lemma~\ref{lemma:ar-regret-bounds}. Define
\begin{align*}
\widehat{\Regret}^{\mathrm{AR26}}_n&:=\min\{\widehat{\Regret}^{(8)}_n,\ln(n+1)+C'\},\qquad
\underline\Wealth_n^{\mathrm{AR26}}:=\Wealth_n^\star\exp(-\widehat{\Regret}^{\mathrm{AR26}}_n),\\
\tau_\alpha^{F_{\mathrm{AR26}}}&:=\inf\{n\geq1:\Wealth_n^{F_{\mathrm{AR26}}}\geq1/\alpha\},\qquad
\underline\tau_\alpha^{\mathrm{AR26}}:=\inf\{n\geq1:\underline\Wealth_n^{\mathrm{AR26}}\geq1/\alpha\}\,.
\end{align*}
For $0<\Delta<1/2$, let $Q_\Delta$ be the point mass at $1/2+\Delta$, and define
\[
A_\Delta:=\sqrt{\ln(1/\Delta)},\qquad \alpha_\Delta:=e^{-A_\Delta}\,.
\]
Recall that $\lstar_Q=\max_{\beta\in[-2,2]}\E_Q[\ln(1+\beta(X-1/2))]$ is the oracle rate, that $\beta_Q^\star$ is its maximizer, and that
\[
B_\alpha(Q)=\frac{\ln(1/\alpha)+\LL(|\beta_Q^\star|)}{\lstar_Q},\qquad
\LL(x)=\ln\max\{1,\ln(1/x)\}\,.
\]
Then
\begin{equation}
\label{eq:ar-counterexample-ratio}
\frac{\E_{Q_\Delta}[\underline\tau_{\alpha_\Delta}^{\mathrm{AR26}}]}{B_{\alpha_\Delta}(Q_\Delta)}\longrightarrow2
\quad\text{as }\Delta\downarrow0\,.
\end{equation}
Moreover,
\begin{equation}
\label{eq:ar-counterexample-exact-ratio}
\frac{\E_{Q_\Delta}[\underline\tau_{\alpha_\Delta}^{\mathrm{AR26}}]}
{\E_{Q_\Delta}[\tau_{\alpha_\Delta}^{F_{\mathrm{AR26}}}]}\longrightarrow2\,.
\end{equation}
In particular, $\eff(\underline\tau^{\mathrm{AR26}})\geq2$, whereas $\eff(\tau^{F_{\mathrm{AR26}}})=1$.
\end{proposition}
\begin{proof}
By \eqref{eq:ar-two-regret-bounds}, $\underline\Wealth_n^{\mathrm{AR26}}\leq\Wealth_n^{F_{\mathrm{AR26}}}$, and Ville's inequality gives validity.

Under $Q_\Delta$, every observation equals $1/2+\Delta$, so
\[
\lstar_{Q_\Delta}=\ln(1+2\Delta),\qquad
\Wealth_n^\star=(1+2\Delta)^n,\qquad
S_n=n\Delta,\qquad V_n=n\Delta^2\,.
\]
For every $n\geq2$, we have $|S_n|=n\Delta\geq\sqrt{2n\Delta^2}=\sqrt{2V_n}$ and $|S_n|=n\Delta>n\Delta^2/10=V_n/10$. Thus, the last branch of \eqref{eq:ar-pathwise-regret} applies, giving
\[
\widehat{\Regret}^{(8)}_n=\tfrac12 n\lstar_{Q_\Delta}+C\,.
\]
Consequently,
\begin{equation}
\label{eq:ar-counterexample-wealth}
\ln\underline\Wealth_n^{\mathrm{AR26}}
=\max\left\{\tfrac12 n\lstar_{Q_\Delta}-C,
n\lstar_{Q_\Delta}-\ln(n+1)-C'\right\}\,.
\end{equation}
The first branch reaches $A_\Delta$ at $\lceil2(A_\Delta+C)/\lstar_{Q_\Delta}\rceil$. No rejection can occur when $n<A_\Delta/\lstar_{Q_\Delta}$, since the proxy is at most $\Wealth_n^\star$. Throughout the remaining possible crossing window
\[
\frac{A_\Delta}{\lstar_{Q_\Delta}}\leq n\leq
\left\lceil\frac{2(A_\Delta+C)}{\lstar_{Q_\Delta}}\right\rceil,
\]
we have $\ln(n+1)\sim\ln(1/\Delta)\gg A_\Delta$ uniformly, while $n\lstar_{Q_\Delta}\leq2A_\Delta+O(1)$. Therefore, the second branch in \eqref{eq:ar-counterexample-wealth} cannot reach $A_\Delta$. For all sufficiently small $\Delta$, the stopping time is deterministic and equals $\lceil2(A_\Delta+C)/\lstar_{Q_\Delta}\rceil$. Since $\beta_{Q_\Delta}^\star=2$ and $\LL(|\beta_{Q_\Delta}^\star|)=0$, the benchmark is $B_{\alpha_\Delta}(Q_\Delta)=A_\Delta/\lstar_{Q_\Delta}$, proving \eqref{eq:ar-counterexample-ratio}.

The density in \eqref{eq:ar-prior-density} has full support and both logarithmic decay indices equal to one, so Corollary~\ref{thm:main} gives $\eff(\tau^{F_{\mathrm{AR26}}})=1$. Theorem~\ref{thm:pointwise} also gives $\E_{Q_\Delta}[\tau_{\alpha_\Delta}^{F_{\mathrm{AR26}}}]\sim B_{\alpha_\Delta}(Q_\Delta)$, proving \eqref{eq:ar-counterexample-exact-ratio}. The factor 2 comes from the term $\tfrac12\ln\Wealth_n^\star$ in \eqref{eq:ar-pathwise-regret}. In \citet[Lemma~B.1]{AgrawalR26}, this term results from replacing the logarithmic mean $(\Wealth_n^\star-1)/\ln\Wealth_n^\star$ by its smaller lower bound $\sqrt{\Wealth_n^\star}$.
\end{proof}

\section{Details on OPT-FFT}
\label{app:prior-moments}

The algorithm in Section~\ref{sec:computation} divides the allocation interval into regions, computes an allocation within each region, and combines these allocations using their accumulated wealth. The main task inside a region is to integrate products of returns. These products are polynomials, so finitely many prior moments suffice to integrate them exactly up to a prescribed horizon. We first construct the integration rule, then explain how the FFT recursion updates it as observations arrive. The final part describes the regions and numerical moments used by OPT-FFT for CO96 and R70.

\subsection{Integrating on one region}

Recall that $\overline F$ is the prior in the portfolio coordinate $p\in[0,1]$. An observation $X_n$ gives the return
\begin{equation}
\label{eq:regional-return}
r_n(p)=p\frac{X_n}{\mu_0}+(1-p)\frac{1-X_n}{1-\mu_0}~.
\end{equation}
Within a region $I_j$, the next allocation is the average of $p$, weighted by its prior mass and the wealth it has earned so far:
\begin{equation}
\label{eq:regional-allocation}
p_{j,n}:=\frac{\int_{I_j}p\prod_{i=1}^{n-1}r_i(p)\,\mathrm{d}\overline F(p)}
{\int_{I_j}\prod_{i=1}^{n-1}r_i(p)\,\mathrm{d}\overline F(p)}~.
\end{equation}
This is the regional posterior mean whenever the denominator is positive. The two integrands have degrees at most $n$ and $n-1$, respectively. Thus, for a horizon $T\geq1$, it is enough to integrate polynomials through degree $D$, where $D$ is the smallest power of 2 at least $T$.

We use the same local coordinate $z\in[-1,1]$ in every region. For a region $I_j$ with endpoints $a_j<b_j$, the map back to the allocation interval and the conditional Chebyshev moments are
\begin{equation}
\label{eq:regional-map-moments}
p_j(z):=\frac{a_j+b_j}2+\frac{b_j-a_j}2z,\qquad
m_{j,k}:=\frac1{\overline F(I_j)}\int_{I_j}T_k\!\left(\frac{2p-a_j-b_j}{b_j-a_j}\right)\,\mathrm{d}\overline F(p)
\quad(0\leq k\leq D),
\end{equation}
where the Chebyshev polynomials are
\[
T_k(z):=\cos(k\arccos z)\qquad(z\in[-1,1],\ k=0,1,\ldots)~.
\]
The normalization by $\overline F(I_j)$ makes these moments conditional on the region; in particular, $m_{j,0}=1$. We omit zero-mass regions. For priors with atoms, each boundary atom belongs to exactly one region. The same construction works in betting coordinates, with $\beta$ and $F$ in place of $p$ and $\overline F$.

We turn these moments into an integration rule on $D+1$ cosine-spaced points:
\begin{align}
z_\ell&:=\cos\!\left(\frac{\pi(\ell+1/2)}{D+1}\right),\nonumber\\
h_{j,\ell}&:=\frac1{D+1}\left(m_{j,0}+2\sum_{k=1}^D m_{j,k}
\cos\!\left(\frac{\pi k(\ell+1/2)}{D+1}\right)\right)
\qquad(0\leq\ell\leq D)~.\label{eq:regional-rule}
\end{align}
To integrate a polynomial, we evaluate it at these points and take the weighted sum with weights $h_{j,\ell}$. Discrete cosine orthogonality makes this sum equal to the prescribed moment for each $T_k$, $0\leq k\leq D$. By linearity, the rule is exact for every polynomial of degree at most $D$. The weights may be negative: they represent integration against the continuous conditional prior. An inverse cosine transform computes them in $O(D\ln(D+1))$ operations.

\subsection{Moving between grids}

Short blocks of observations need fewer grid points. The recursion therefore moves between grids of different sizes, using the one-dimensional operations of \citet[Appendix~E.3]{LiuPO26}. There are two operations: \texttt{Resample} evaluates a polynomial on another grid, and \texttt{Restrict} transfers integration weights to a smaller grid while preserving the integrals needed there.

A block of $d$ observations uses $d+1$ points. Write the grid vector and the initial integration weights as
\begin{equation}
\label{eq:regional-nodes}
\bz^{(d)}:=\left(\cos\frac{\pi(\ell+1/2)}{d+1}\right)_{\ell=0}^d,
\qquad \bh_j:=(h_{j,\ell})_{\ell=0}^D~.
\end{equation}
Both operations use a discrete cosine transform. We specify its normalization explicitly so that the formulas can be implemented with a standard FFT library. The orthonormal type-II cosine matrix $\bC_d\in\R^{(d+1)\times(d+1)}$ is
\begin{equation}
\label{eq:regional-dct}
(\bC_d)_{k\ell}:=
\begin{cases}
1/\sqrt{d+1},&k=0,\\
\sqrt{2/(d+1)}\cos\!\left(\dfrac{\pi k(\ell+1/2)}{d+1}\right),&1\leq k\leq d,
\end{cases}
\qquad 0\leq\ell\leq d~.
\end{equation}
Multiplication by $\bC_d$ converts grid values to cosine coefficients; multiplication by $\bC_d^\top$ converts them back. To resample, we transform the values, adjust the number of coefficients, and apply the inverse transform on the new grid:
\begin{equation}
\label{eq:regional-resample}
\operatorname{Resample}_{d\to d'}(\bq)
:=\sqrt{\frac{d'+1}{d+1}}\,\bC_{d'}^\top[\bC_d\bq]_{0:d'}~.
\end{equation}
Here, $[\cdot]_{0:d'}$ retains coefficients $0,\ldots,d'$, padding with zeros if the original vector is shorter. The square-root factor accounts for the different transform normalizations. The result gives the same polynomial on the new grid whenever its degree is at most $\min\{d,d'\}$.

For integration weights, we want a different property: the dot product with polynomial values should stay the same after moving to the smaller grid. This requires the transpose of interpolation from that smaller grid:
\begin{equation}
\label{eq:regional-restrict}
\operatorname{Restrict}_{d\to d'}(\bu)
:=\frac{d+1}{d'+1}\operatorname{Resample}_{d\to d'}(\bu)
=\operatorname{Resample}_{d'\to d}^{\top}\bu,
\qquad d'\leq d~.
\end{equation}
In particular, the restricted weights give the same integral for every polynomial of degree at most $d'$. We apply these transforms by FFT without forming the dense matrices in \eqref{eq:regional-dct}. Each operation costs $O(\max\{d,d'\}\ln(\max\{d,d'\}+1))$. In the algorithms below, $\odot$ means componentwise multiplication, and a scalar function applied to a vector is evaluated at each entry.

\subsection{Computing allocations as observations arrive}

Algorithm~\ref{alg:regional-fft} coordinates the regions. On each round, it obtains their allocations, averages them using prior mass times past wealth, and then observes $X_n$. It uses the resulting returns to update wealth and passes $X_n$ to each region to prepare the next allocation. Algorithm~\ref{alg:regional-prefix} performs the computation within one region. It is our fixed-start specialization of the \texttt{Prefix} recursion in \citet[Appendix~E.4]{LiuPO26}.

We write each regional computation as a \emph{generator}: a routine that pauses when it produces an allocation and resumes when the next observation arrives. It keeps all its local variables while paused. The operation $\operatorname{next}$ runs it to its first allocation; $\operatorname{send}(X_n)$ supplies the observation and runs to the next allocation or the end of the block. Thus, $p_{j,n}$ is always available before $X_n$ is observed. If $T<D$, we stop after observation $T$ and leave the unfinished recursions paused.

\begin{algorithm}[t]
\caption{Regional FFT portfolio with a prescribed horizon}
\label{alg:regional-fft}
\begin{algorithmic}[1]
\REQUIRE Horizon $T\geq1$, $\mu_0\in(0,1)$, $\alpha\in(0,1)$, leaf size $d_0$, and regional endpoints, masses, and moments from \eqref{eq:regional-map-moments}
\STATE $D:=2^{\lceil\log_2 T\rceil}$; form $\bz^{(D)}$ using \eqref{eq:regional-nodes}; moments are required through degree $D$
\STATE $\Wealth_0:=1$; $\Wealth_{j,0}:=1$ for $j=1,\ldots,K$
\FOR{$j:=1,\ldots,K$}
  \STATE $\displaystyle\bh_j:=\sqrt{\frac2{D+1}}\,\bC_D^\top(m_{j,0}/\sqrt2,m_{j,1},\ldots,m_{j,D})^\top$
  \STATE $\mathcal A_j:=\operatorname{Prefix}(j,1,D,\bh_j,\bz^{(D)}\odot\bh_j)$, a generator from Algorithm~\ref{alg:regional-prefix}
  \STATE $p_{j,1}:=\operatorname{next}(\mathcal A_j)$
\ENDFOR
\FOR{$n:=1,\ldots,T$}
  \STATE $\displaystyle p_n:=\frac{\sum_{j=1}^K\overline F(I_j)\Wealth_{j,n-1}p_{j,n}}{\sum_{j=1}^K\overline F(I_j)\Wealth_{j,n-1}} \qquad$  (use $p_n:=\mu_0$ if the denominator is zero)
  \STATE Allocate fraction $p_n$ to the first stock, then observe $X_n$
  \STATE $\Wealth_n:=\Wealth_{n-1}r_n(p_n)$; $\Wealth_{j,n}:=\Wealth_{j,n-1}r_n(p_{j,n})$ for every $j$, using \eqref{eq:regional-return}
  \IF{$n<T$}
    \STATE $p_{j,n+1}:=\mathcal A_j.\operatorname{send}(X_n)$ for every $j$
  \ELSE
    \STATE Resume each $\mathcal A_j$ with $X_n$ and discard its next yield or completed block product
  \ENDIF
  \STATE Report $\Wealth_n$ and record the first $n$ with $\Wealth_n\geq1/\alpha$, if any
\ENDFOR
\end{algorithmic}
\end{algorithm}

\paragraph{What a recursive call stores.}
A call $\operatorname{Prefix}(j,L,d,\bu,\bv)$ handles the $d$ observations starting at round $L$. It carries two integration-weight vectors of length $d+1$. The vector $\bu$ computes the denominator of the local posterior mean, and $\bv$ computes its numerator in the coordinate $z$. Both already include the returns from observations before $L$. More precisely, $\bu$ integrates those past returns times any polynomial of degree at most $d$; $\bv$ does the same with an extra factor $z$, for polynomials of degree at most $d-1$. Initially these vectors are $\bh_j$ and $\bz^{(D)}\odot\bh_j$.

In a short block, the vector $\bq$ stores the product of the returns observed within that block, evaluated on its grid. When the denominator is positive, the ratio $\bv^\top\bq/(\bu^\top\bq)$ gives the posterior mean of $z$, which we map back to $I_j$ to obtain \eqref{eq:regional-allocation}. Computing the mean in local coordinates keeps the numerator's scale independent of the region's location. Before each observation, this block product has degree at most $d-1$, so the numerator also fits within the available moments. We process blocks of at most $d_0$ observations this way, where the leaf size $d_0$ is a fixed positive integer, 32 by default.

For a longer block, we process its left half first. Once that half is complete, we multiply its return product into both integration-weight vectors and pass the updated weights to the right half. This is how the right half accounts for all earlier observations. Each call returns the product of all returns in its block, evaluated on its own grid.

\begin{algorithm}[t]
\caption{$\operatorname{Prefix}(j,L,d,\bu,\bv)$: yield regional allocations on rounds $L,\ldots,L+d-1$}
\label{alg:regional-prefix}
\begin{algorithmic}[1]
\REQUIRE Region $j$, block start $L$, power-of-two block length $d$, functional weights $\bu,\bv\in\R^{d+1}$; $\mu_0$ and $d_0$ are shared with Algorithm~\ref{alg:regional-fft}
\STATE Use $p_j$ from \eqref{eq:regional-map-moments}, $\bz^{(d)}$ from \eqref{eq:regional-nodes}, and the operations \eqref{eq:regional-resample}--\eqref{eq:regional-restrict}
\IF{$d\leq d_0$}
  \STATE $\bq:=(1,\ldots,1)\in\R^{d+1}$
  \FOR{$n:=L,\ldots,L+d-1$}
    \STATE $\widehat z:=0$
    \IF{$\bu^\top\bq>0$}
      \STATE $\displaystyle\widehat z:=\min\left\{1,\max\left\{-1,\frac{\bv^\top\bq}{\bu^\top\bq}\right\}\right\}$
    \ENDIF
    \STATE Yield $p_{j,n}:=(a_j+b_j)/2+(b_j-a_j)\widehat z/2$ and receive $X_n$ on resumption
    \STATE $\bq:=\bq\odot r_n(p_j(\bz^{(d)}))$, evaluating \eqref{eq:regional-return} at all $d+1$ grid points
  \ENDFOR
  \RETURN $\bq$
\ELSE
  \STATE $k:=d/2$
  \STATE $\bq_{\mathrm L}:=\operatorname{Prefix}(j,L,k,\operatorname{Restrict}_{d\to k}(\bu),\operatorname{Restrict}_{d\to k}(\bv))$
  \STATE $\bq_{\mathrm L}:=\operatorname{Resample}_{k\to d}(\bq_{\mathrm L})$
  \STATE $\bq_{\mathrm R}:=\operatorname{Prefix}(j,L+k,k,\operatorname{Restrict}_{d\to k}(\bu\odot\bq_{\mathrm L}),\operatorname{Restrict}_{d\to k}(\bv\odot\bq_{\mathrm L}))$
  \RETURN $\bq_{\mathrm L}\odot\operatorname{Resample}_{k\to d}(\bq_{\mathrm R})$
\ENDIF
\end{algorithmic}
\end{algorithm}

The two recursive assignments in Algorithm~\ref{alg:regional-prefix} use the programming convention \texttt{yield from}. While a child call is running, the parent forwards its allocations to Algorithm~\ref{alg:regional-fft} and sends each new observation back to it. The assignment finishes only when the child returns its completed block product. Thus, the right call starts only after all observations in the left block have arrived.

\paragraph{Combining CO96 and R70.}
Run Algorithm~\ref{alg:regional-fft} separately for each component. Section~\ref{app:regional-setup} gives their partitions, masses, and moments. The pseudocode uses portfolio coordinates for both components. For R70, map the region endpoints from $\beta\in[-1,1]$ using $p=\mu_0+\beta\mu_0(1-\mu_0)$; the moments in the local coordinate $z$ remain unchanged. Write $p_n^{\mathrm{CO96}}$ and $p_n^{\mathrm{R70}}$ for the allocations returned by the two runs. For the prior mixture in \eqref{eq:remark2_test}, we again average allocations using their past wealths:
\begin{equation}
\label{eq:regional-component-hedge}
\begin{aligned}
p_n^{F_w}&:=\frac{w\Wealth_{n-1}^{F_{\mathrm{CO96}}}p_n^{\mathrm{CO96}}+(1-w)\Wealth_{n-1}^{F_{\mathrm{R70}}}p_n^{\mathrm{R70}}}
{w\Wealth_{n-1}^{F_{\mathrm{CO96}}}+(1-w)\Wealth_{n-1}^{F_{\mathrm{R70}}}},\\
\Wealth_n^{F_w}&=\Wealth_{n-1}^{F_w}r_n(p_n^{F_w})
=w\Wealth_n^{F_{\mathrm{CO96}}}+(1-w)\Wealth_n^{F_{\mathrm{R70}}}~.
\end{aligned}
\end{equation}
Start all wealths at 1 and use $\mu_0$ as the allocation if the denominator is zero. For the mixture test, compare only the combined wealth, $\Wealth_n^{F_w}$, with $1/\alpha$. The implementation keeps R70 in betting coordinates: if its output is $\beta_n^{\mathrm{R70}}$, convert it by $p_n^{\mathrm{R70}}=\mu_0+\beta_n^{\mathrm{R70}}\mu_0(1-\mu_0)$ before using \eqref{eq:regional-component-hedge}. The return then equals $1+\beta_n^{\mathrm{R70}}(X_n-\mu_0)$.

\paragraph{Rescaling and extending the horizon.}
Long products can overflow or underflow. We can rescale a block product by any positive scalar without changing the ratio that determines the allocation. We can also rescale both integration-weight vectors of a region by the same positive scalar. The pseudocode omits these rescalings; the implementation uses them and stores regional wealth logarithmically. If either inner product is nonfinite, it uses the region's midpoint. The projection onto $[-1,1]$ in Algorithm~\ref{alg:regional-prefix} keeps the allocation inside its region; it is inactive in exact arithmetic whenever regional wealth is positive.

The supplied implementation uses a fixed horizon. To continue after round $D$ in exact arithmetic, prepare new integration rules for horizon $2D$ with the same prior, and replay the first $D$ observations to reconstruct the recursive state. Discard the allocations produced during replay, then resume betting at round $D+1$ with the existing wealth. Exactness makes the combined allocations independent of the horizon and partition whenever total wealth is positive. Summing the work over doubled horizons preserves the arithmetic and memory orders in Proposition~\ref{prop:regional-fft}, using the largest region count reached.

\begin{proposition}[Exactness and arithmetic cost]
\label{prop:regional-fft}
Fix $T\geq1$, a candidate mean $\mu_0\in(0,1)$, and a partition into $K$ positive-mass regions. Given their exact prior masses and conditional Chebyshev moments through degree $D$, Algorithms~\ref{alg:regional-fft} and~\ref{alg:regional-prefix} compute $\Wealth_n^F$ for every $1\leq n\leq T$ exactly in exact arithmetic. With constant leaf size, they use $O(KT\ln^2(T+1))$ operations and $O(KT)$ memory, excluding the computation of the supplied moments.
\end{proposition}
\begin{proof}
Before observation $n$, the prefix product has degree $n-1$ and its product with the local coordinate has degree $n$. Both are at most $D$, so the integration rule returns the exact conditional posterior mean. Every recursive split transfers a functional by the transpose of exact polynomial interpolation. Multiplication by the completed left product gives the functional needed for the right prefixes, with sufficient degree at every node. Thus, the recursion preserves these integrals.

Let $\Wealth_{j,n}$ denote the mixture wealth under the prior conditional on $I_j$. Its posterior-mean allocation has return $\Wealth_{j,n}/\Wealth_{j,n-1}$ whenever the denominator is positive. Weight the regional allocations before observation $n$ by
\[
\frac{\overline F(I_j)\Wealth_{j,n-1}}{\sum_\ell\overline F(I_\ell)\Wealth_{\ell,n-1}}~.
\]
Since returns are affine in the allocation, the combined return telescopes to wealth $\sum_j\overline F(I_j)\Wealth_{j,n}=\Wealth_n^F$. A region with zero wealth stays at zero and receives zero weight; if all regions have zero wealth, any allocation suffices.

At recursion size $d$, the transforms take $O(d\ln(d+1))$ work per region, followed by two recursions of half the size. Summing over the levels gives $O(T\ln^2(T+1))$ work per region. Arrays along a recursion branch have geometrically decreasing sizes, so the workspace, cached initial functionals, and allocation traces use $O(T)$ memory per region. Constructing the root rules and combining the regional strategies fit within these bounds.
\end{proof}

\subsection{Setting up the CO96 and R70 priors}
\label{app:regional-setup}

The algorithms take a partition, its prior masses, and its conditional moments as inputs. We now give the choices used by OPT-FFT: \eqref{eq:regional-spacing}--\eqref{eq:regional-defaults} specify the regions, \eqref{eq:co96-region-mass}--\eqref{eq:r70-region-mass} give their masses, and the numerical procedure below computes the moments in \eqref{eq:regional-map-moments}.

\paragraph{Choosing the regions.}
We place more boundaries near the zero bet and the support endpoints, where small betting scales and concentrated wealth products need resolution. Choose an integer $P\geq1$ before observing data and define
\begin{equation}
\label{eq:regional-spacing}
s_j:=\sin^2\!\left(\frac{\pi j}{2P}\right)\qquad(0\leq j\leq P)~.
\end{equation}
For CO96, use the boundaries $\mu_0(1-s_j)$ and $\mu_0+(1-\mu_0)s_j$. They give $2P$ intervals, clustered around $p=0,\mu_0,1$. For R70, use betting coordinates: take the central interval $[-s_1,s_1]$ and the intervals $[s_j,s_{j+1}]$ and $[-s_{j+1},-s_j]$, $1\leq j<P$. This gives $2P-1$ regions. Together, their conditional priors retain the whole support, including arbitrarily small nonzero bets.

The default number of regions depends on the horizon through
\begin{equation}
\label{eq:regional-defaults}
L:=\lceil\log_2(T+1)\rceil,\qquad
P:=\min\!\left\{2^{\lceil\log_2(L^2)\rceil},\max\{1,D/16\}\right\}~.
\end{equation}
Here, $D$ is the power-of-two horizon used above. The cap limits the number of regions on short streams. For longer streams, the squared logarithm increases resolution while keeping $P=O(\ln^2(T+1))$. With leaf size 32, Proposition~\ref{prop:regional-fft} gives $O(T\ln^4(T+1))$ work and $O(T\ln^2(T+1))$ memory. This choice is empirical: exactness holds for every partition, while numerical stability depends on the partition.

\paragraph{Prior masses.}
These can be computed analytically. For CO96, in portfolio coordinates,
\begin{equation}
\label{eq:co96-region-mass}
\overline F_{\mathrm{CO96}}([a,b])=\frac2\pi\bigl(\arcsin\sqrt b-\arcsin\sqrt a\bigr)~.
\end{equation}
For R70, recall $c=6.6e$ from \eqref{eq:oj-robbins-density}. Its central and positive interval masses are
\begin{equation}
\label{eq:r70-region-mass}
\begin{aligned}
F_{\mathrm{R70}}([-h,h])&=\frac{\ln\ln c}{\ln\ln(c/h)}&& (0<h\leq1),\\
F_{\mathrm{R70}}([a,b])&=\frac{\ln\ln c}{2}\left(\frac1{\ln\ln(c/b)}-\frac1{\ln\ln(c/a)}\right)&& (0<a<b\leq1)~.
\end{aligned}
\end{equation}
The negative interval masses satisfy $F_{\mathrm{R70}}([-b,-a])=F_{\mathrm{R70}}([a,b])$ for $0<a<b\leq1$. For narrow intervals, directly subtracting two nearly equal numbers can lose accuracy. We use reflected coordinates or logarithmic identities when they reduce this cancellation.

\paragraph{Moments away from the R70 singularity.}
The substitution $z=\cos\theta$ turns $T_k(z)$ into $\cos(k\theta)$, so the moments become cosine integrals over $[0,\pi]$. At CO96 support endpoints, the Jacobian cancels the square-root singularity; we perform this cancellation analytically before numerical evaluation. For an accuracy target $0<\epsilon<1$, we divide the angular interval into panels of width at most $4/D$, apply positive Gaussian quadrature of order $\lceil\max\{16,\ln(1/\epsilon)/\ln8\}\rceil$ on each panel, and normalize the weights to sum to one. The default target is $10^{-13}$.

Evaluating each cosine sum separately would cost quadratically in $D$. We instead compute all the moments together by FFT. We round the phases to a uniform grid of at least $4D$ points, expand the residual phase in a Taylor series, and sum each term by FFT. In every Fourier mode we need, the residual phase has magnitude at most $\pi/4$. Thus, $O(\ln(1/\epsilon))$ terms give absolute truncation error at most $\epsilon$ for normalized positive weights.

\paragraph{The central R70 region.}
The region $[-h,h]$ needs special treatment because its density is unbounded at zero. Since $f_{\mathrm{R70}}(-\beta)=f_{\mathrm{R70}}(\beta)$, we can integrate over $(0,h]$ and pair the contributions at $\beta$ and $-\beta$. Use the change of variables
\[
a_h:=\ln(c/h),\qquad \beta:=h\sin\phi\quad(0<\phi\leq\pi/2)~.
\]
The identity $T_k(-z)=(-1)^kT_k(z)$ makes all odd moments zero, while the even moments are
\[
m_{j,2k}=(-1)^k\int_0^{\pi/2}\cos(2k\phi)
\frac{\ln a_h\cos\phi}{\sin\phi\,(a_h-\ln\sin\phi)[\ln(a_h-\ln\sin\phi)]^2}\,\mathrm{d}\phi~.
\]
The density in this integral has total mass one. Near zero, the cosine factor is close to 1, so we can integrate that part using its exact mass. For $0<u\leq\pi/2$, this mass is
\[
\frac{\ln a_h}{\ln(a_h-\ln\sin u)}~.
\]
Choose the cutoff
\[
u:=\frac{\sqrt\epsilon}{2\max\{2,D\}}~,
\]
and replace $\cos(2k\phi)$ by 1 on $(0,u]$. Since $2k\leq D$, the resulting error in each moment is at most
\[
\sup_{0\leq\phi\leq u}|1-\cos(2k\phi)|\leq D^2u^2/2\leq\epsilon/8~.
\]
This step retains the prior mass near zero and approximates only the cosine factor. On the remaining interval $[u,\pi/2]$, we use geometrically spaced panels near zero, followed by panels of width at most $4/D$. The first resolve the density near its singularity; the second resolve the cosine oscillations. The Gaussian quadrature order grows as $O(\ln(D/\epsilon))$. The same FFT summation computes all even moments together, and we set the zeroth moment to 1.

\paragraph{Preprocessing cost and accuracy.}
For fixed quadrature oversampling, computing the moments and integration rules takes $O(K(D+\ln(1/\epsilon))\ln(D/\epsilon)\ln(1/\epsilon))$ operations. Storing the moments and rules takes $O(KD)$ memory. The central R70 quadrature has $O(D+\ln(1/\epsilon))$ panels and order $O(\ln(D/\epsilon))$, and the implementation stores all its nodes. Processing regions one at a time therefore requires an additional $O((D+\ln(1/\epsilon))\ln(D/\epsilon))$ workspace, which is $O(D\ln D)$ at fixed tolerance. With our default region count, total memory remains $O(T\ln^2(T+1))$, and preprocessing takes $\widetilde O(T)$ operations at fixed tolerance.

The bounds above control the approximation near zero and the truncation of the Fourier expansion. We check the remaining Gaussian quadrature error by increasing its order. Accordingly, $\epsilon$ is a numerical target, not a certified bound on all moments or accumulated wealth. The implementation also accepts conditional moments computed by another routine; Proposition~\ref{prop:regional-fft} applies whenever those input moments and masses are exact.

\paragraph{Default hyperparameters.}
Our default partitions use squared-sine spacing on each side of the zero bet, with a polylogarithmic number of intervals rounded to a power of 2. At horizons 512, 2,048, and 8,192, this gives respectively 64, 256, and 512 CO96 regions, and 63, 255, and 511 R70 regions. These are continuous regions, not point masses. The degree within each region grows with the horizon. A leaf size of 32 avoids transforms for short products, and the numerical moment routine uses an accuracy target of $10^{-13}$. The spacing targets the sensitive parts of the domain; the region count and leaf size balance their resolution against transform overhead. Appendix~\ref{app:prior-moments} specifies the choices.

We empirically validated these hyperparameter choices on 15 paths: Bernoulli and bounded reversals, Bernoulli alternatives with means $0.505$ and $0.7$, and $\mathrm{Beta}(4,3)$ observations, each at lengths 512, 2,048, and 8,192, with $\mu_0=1/2$. Against refined positive quadrature, the largest discrepancy in mixture log wealth was $6.7\times10^{-11}$, and no regional allocation required projection. Analytic Bernoulli calculations and independent integration provide additional checks.

\end{document}